\documentclass[11pt,letterpaper]{amsart}

\usepackage{amssymb}
\usepackage{color}

\newtheorem{theorem}{Theorem}[section]
\newtheorem{corollary}[theorem]{Corollary}
\newtheorem{definition}[theorem]{Definition}

\newtheorem{lemma}[theorem]{Lemma}
\newtheorem{proposition}[theorem]{Proposition}
\newtheorem{remark}[theorem]{Remark}

\theoremstyle{definition}
\theoremstyle{remark}
\numberwithin{equation}{section}

\begin{document}
\begin{sloppypar}
	
\title{A General and Unified Method to Prove the Uniqueness of Ground State Solutions and the Existence/Non-existence, and Multiplicity of Normalized Solutions with applications to various NLS}
		
\author{Hichem Hajaiej}
\address[H. Hajaiej]{ Department of Mathematics, California State University at Los Angeles, Los Angeles, CA 90032, USA}
\email{ hhajaie@calstatela.edu}
	
\author{Linjie Song}		
\address[L. Song]{Institute of Mathematics, AMSS, Academia Sinica, Beijing 100190, China}
\address[L. Song]{University of Chinese Academy of Science, Beijing 100049, China}
\email{ songlinjie18@mails.ucas.edu.cn.}
\thanks{Linjie Song is supported by CEMS}

\date{}

		\begin{abstract}
			
			We first give an abstract framework to show the uniqueness of Ground State Solutions (GSS) for a large class of PDEs. To the best of our knowledge, all the existing results in the literature only addressed particular cases. Moreover, our self-contained approach offers a general framework to study the existence/non-existence and multiplicity of normalized solutions. We will exhibit concrete examples to which our method applies, and verify all the assumptions we need. Our approach is applicable to a wide range of operators and domains provided that our hypotheses are verified. Additionally, we prove new results about the non-degeneracy and uniqueness of positive GSS in a general setting. Our findings are applicable to fractional nonlinear Schr\"{o}dinger equations with non-autonomous nonlinearities. In particular, we were able to extend the main results of \cite{FL, FLS} to general non-autonomous and mixed nonlinearities. This does not seem possible by using the approach developed by the authors of the above breakthrough papers.  The orbital stability/instability of the standing waves will be addressed thanks to the non-degeneracy.
			
			
			
			\textbf{Data availability statement:} My manuscript has no associate data.
		\end{abstract}
\maketitle
		
		\section{\textbf{Introduction}}
		
		Since the breakthrough paper of Lions in 1985, \cite{Lio1}, on the existence of normalized solutions in the sub-critical case using the concentration-compactness method, and the one of Cazenave and Lions in 1986 on the orbital stability of standing waves, \cite{CazLio}, there were no additional contributions in these two directions involving a heuristic or general methods until quite recently. In 2017, Bellazzini, Boussaid, Jeanjean, and Visciglia developed a general approach to show the existence of normalized solutions in the super-critical case, \cite{BBJV}. Many mathematicians borrowed their idea to address the existence of normalized solutions and to study the orbital stability of standing waves in the super-critical case. To the best of our knowledge, there are no general approaches addressing the uniqueness of GSS or the existence of normalized solutions in any regime.
		
		The main goal of this work is to address these two questions and to provide a unified and very general method that establishes the uniqueness of GSS and applies to study normalized solutions in the sub-critical, critical or super-critical cases, bounded, exterior or entire space, Laplacian, fractional Laplacian, $p$-Laplacian, $p$-fractional Laplacian, for autonomous or non-autonomous nonlinearities and with or without the presence of a potential. We will exhibit some concrete examples of the latter and discuss all the details. Our abstract framework applies to many other situations and our method is quite flexible as the reader will notice it.

  Before introducing the notation and our main results (Theorem \ref{thmc.1}, Theorem \ref{thmc.4}, and Theorem \ref{athmc.7}), let us point out that we believe that the ideas and ingredients used in this paper are novel. The subtle interplay between the bifurcation theory, spectral theory, and variational methods turned out to be powerful to show the uniqueness of GSS.  Let us also mention that some of the results in Section 3 can seem stated under strong assumptions (non-degeneracy, convergence of (P.S) sequence). We count on the reader's patience to see how all these assumptions were satisfied in our concrete examples. Additionally, these hypotheses are optimal in this framework as it is shown in a ongoing work of the authors. This means that our framework is fruitful as it addresses many difficult challenges. In particular, we extend the breakthrough results of \cite{FL, FLS} to non-autonomous nonlinearities using our self-containing approach. In what follows, we will describe our approach:
		
		Assume that $W$ is a Banach space with the norm $\|\cdot\|$, $\Phi_{\lambda} \in C^{2}(W, \mathbb{R})$ is a Fr\'{e}chet-differentiable functional with derivative $D_{u}\Phi_{\lambda}: W \rightarrow W^{\ast}$, where the space $W^{\ast}$ is the dual space of $W$ and $\lambda \in \mathbb{R}$. Consider a particular form of $\Phi_{\lambda}$, (so that \eqref{eq1.1} has an underlying variational structure):
		\begin{equation} \label{eq1.1}
			D_{u}S(u) + D_{u}G(u) - D_{u}F(u) - \lambda D_{u}Q(u) = 0,\quad\text{i.e.}\ u \in W,
		\end{equation}
		where $\Phi_{\lambda} = S(u) + G(u) - F(u) - \lambda Q(u)$, $S(u)$ and $Q(u)$ are positive for any $u \in W \backslash \{0\}$, $G(u)$ is nonnegative for any $u \in W \backslash \{0\}$. Assume that $D_{u}\Phi_{\lambda}(0) = 0$, ( $0$ is always a trivial solution of (\ref{eq1.1}) for any $\lambda$). In this paper, we are interested in the nontrivial solutions. In previous contributions, there are two different approaches to study the existence of solutions of equation (\ref{eq1.1}) using variational methods: When $\lambda$ is prescribed, we study $\Phi_{\lambda}$ directly and find its critical points $u_{\lambda}$. However, we don't have information about $Q(u_{\lambda})$. When $\rho > 0$ is prescribed, we study $I(u) = S(u) + G(u) - F(u)$ under the constraint $Q(u) = \rho$ and get critical points $u_{\rho}$ of $I_{|\{Q(u) = \rho\}}$. Then $(\lambda_{\rho}, u_{\rho})$ is a solution of (\ref{eq1.1}), where $\lambda_{\rho}$ is determined as a Lagrange multiplier and unknown, $u_{\rho}$ is a normalized solution that satisfies the normalized condition $Q(u_{\rho}) = \rho$. Normalized solutions are the most important ones because they are known to be the most stable. The latter makes them very attractive especially that they also usually have nice symmetry properties. In this paper, we give an abstract framework in which $\Phi_{\lambda}$ is studied directly and the information of $Q(u)$ is prescribed.
		
		Our method is self-contained and novel. It provides a unified and general framework to prove the existence of normalized solutions. To the best of our knowledge, all previous works had completely different approaches for the sub-critical and the super-critical cases. Moreover, our method applies to equations on the entire space, bounded domains, and exterior domains provided that some assumptions are satisfied. Our unified approach is also novel in this sense. Surprisingly, there are few papers studying the existence of normalized solutions on bounded domains. As far as we know, the result in Section \ref{p} is the first result establishing the existence of normalized solutions for $p$-Laplacian equations on bounded domains. Moreover, the literature is silent regarding the existence of normalized solutions on exterior domains. To the best of our knowledge, the first result in this direction has been proven using the method we developed in this paper ( see \cite{Song2,Song3}). Additionally, we will also address the orbital stability of standing waves by exploring non-degeneracy results that we also needed to achieve other purposes - namely the uniqueness of ground state solutions-.		
		
		Our abstract framework contains three parts. They are all independent and complementary. First, let us introduce the abstract framework 2. For $\lambda \in (a, b)$, we assume that $\mathcal{N}_{\lambda} (\subset W)$ is a natural constraint, which means that any critical point of $\Phi_{\lambda}$ restricted to $\mathcal{N}_{\lambda}$ is a critical point of $\Phi_{\lambda}$ in $W$, and $h(\lambda) = \inf_{u \in \mathcal{N}_{\lambda}}\Phi_{\lambda}(u)$ can be achieved by $u_{\lambda}$, which is called a ground state solution. Set
		$$
		\mathcal{K}_{\lambda}:= \{u \in \mathcal{N}_{\lambda}: D_{u}\Phi_{\lambda}(u) = 0, \Phi_{\lambda}(u) = h(\lambda)\}.
		$$
		If $\Phi_{\lambda|\mathcal{N}_{\lambda}}$ satisfies $(PS)_{h(\lambda)}$ condition (i.e. any sequence $\{u_{n}\}_{n = 1}^{+\infty} \subset \mathcal{N}_{\lambda}$ such that $\Phi_{\lambda}(u_{n}) \rightarrow h(\lambda)$ and $D_{u}\Phi_{\lambda}(u_{n}) \rightarrow 0$ has a convergent subsequence in $\mathcal{N}_{\lambda}$), then $\mathcal{K}_{\lambda}$ is compact and for $\lambda \in (a, b)$, we can define
		\[
		d(\lambda):= \inf_{u \in \mathcal{K}_{\lambda}}Q(u) = \min_{u \in \mathcal{K}_{\lambda}}Q(u),
		\]
		\[
		\hat{d}(\lambda):= \sup_{u \in \mathcal{K}_{\lambda}}Q(u) = \max_{u \in \mathcal{K}_{\lambda}}Q(u).
		\]
		When $d \equiv \hat{d}$ in $(a, b)$, we can prove that $d(\lambda)$ is continuous. This property will play a crucial role in our discussion. The behavior of $d(\lambda)$ is studied when $\lambda \rightarrow a$ and $\lambda \rightarrow b$ later. Then for any $\rho \in (\inf_{\lambda \in (a, b)}d(\lambda), \sup_{\lambda \in (a, b)}d(\lambda))$, there exists $u_{\lambda} \in \mathcal{K}_{\lambda}$ such that $Q(u_{\lambda}) = \rho$ and $D_{u}\Phi_{\lambda}(u_{\lambda}) = 0$.
		
		\begin{remark}
			\label{rmk1.1} Usually, the Nehari manifold can be chosen as the natural constraint. Define
			$$
			\mathcal{N}_{\lambda}:= \{u \in W \backslash \{0\}: D_{u}\Phi_{\lambda}(u)(u) = 0\}.
			$$
			Then any solution of (\ref{eq1.1}) belongs to $\mathcal{N}_{\lambda}$. Suppose that $\Phi_{\lambda}(u_{\lambda}) = \inf_{u \in U}\Phi_{\lambda}(u)$ and $\Phi_{\lambda} \in C^{2}(U, \mathbb{R})$ where $U$ is some neighbourhood of $u_{\lambda}$ in $W$. Then by the Lagrange multiplier rule, one has $D_{u}\Phi_{\lambda}(u_{\lambda}) = \mu (D_{u}\Phi_{\lambda}(u_{\lambda}) + D_{uu}\Phi_{\lambda}(u_{\lambda})(u_{\lambda}, \cdot))$ for some $\mu$. Plugging into this equality $u_{\lambda}$, we obtain $\mu D_{uu}\Phi_{\lambda}(u_{\lambda})(u_{\lambda}, u_{\lambda}) = 0$. Hence, if
			$$
			D_{uu}\Phi_{\lambda}(u_{\lambda})(u_{\lambda}, u_{\lambda}) \neq 0,
			$$
			then $\mu = 0$ and therefore $D_{u}\Phi_{\lambda}(u_{\lambda}) = 0$. The Nehari manifold method was introduced in \cite{Ne1, Ne2}, \cite{Il, SW, Wi} and the references therein provide more details.
		\end{remark}

		A key ingredient of the framework 2 is to show that the quantity $d(\lambda)= \hat{d} (\lambda)$, which holds true if $d$ is a continuous function. In some cases, this is ensured by the uniqueness of the ground state solution ( see Section 5.1, Section \ref{p}). In Section \ref{p}.3, we provide some examples in which this uniqueness is known. In other cases, we can show the continuity of $d(\lambda)$ and consequently establish the existence of normalized solutions even if we do not know that the positive ground state is unique ( see Section 5.5). The fNLS with a Hardy potential falls into this category. For such case, we will only provide examples for NLS that can be straightforwardly extended to other PDEs (fNLS, for example).
		
		In the abstract framework 3, $a = -\infty$, $b = 0$, and we are mainly concerned with the situations when $\lambda \rightarrow -\infty$ and $\lambda \rightarrow 0^-$. Under some abstract settings, we can prove that the local uniqueness
		of the ground state solution when $\lambda \rightarrow -\infty$ or $\lambda \rightarrow 0^-$, and obtain the existence of normalized solutions. Some scaling methods are used and we assume that the limit equation has a unique, non-degenerate ground state solution. We will give an application to mixed nonlinearities, this kind of nonlinearity has gained a lot of interest in the last years ( see e.g. \cite{S1,S2,Ste,WW}), and we cover the most attractive cases, i.e, cubic-quintic ( see Section \ref{mn}). Moreover, we will show how our approach applies to mixed fractional Laplacians ( see Section \ref{mfl}). This type of equations arises in various fields ranging from biophysics to population dynamics and has very recently received an increasing interest. Non-autonomous mixed nonlinearities will be considered.
	
	    In the abstract {\bf framework 2}, we assume {\bf the uniqueness of the ground state solution}. Note that the abstract framework 3 does not provide any information when the $\lambda$ is away from $-\infty$ and $0$. {\bf The abstract framework 1 will complete the picture}. More precisely, it will establish a $C^1$ global branch parameterized by $\lambda$ under a non-degeneracy assumption. Furthermore, once we show that the ground state is unique when $\lambda \rightarrow -\infty$, we can prove the uniqueness of the global branch. {\bf Thus, the abstract framework 1 provides a framework to prove the uniqueness of all ground states}. We will exhibit concrete examples showing that our abstract framework plays an important role for non-local operators. Regarding scalar field equations involving a fractional Laplacian, it is more difficult to show the uniqueness since $(-\Delta)^{s}$ is a nonlocal operator. When the nonlinearity is pure power $|u|^{q-2}u$, \cite{FL, FLS} showed the uniqueness of positive ground state, leaving the problem open for general nonlinearities. To the best of our knowledge, the problem is still open when the nonlinearity is non-autonomous. We will provide a class of non-autonomous nonlinearities for which the positive ground state solution is unique and provide two different methods to prove it, one is based on our abstract framework and one is inspired with some ideas developed in \cite{FL}. Furthermore, our general method also applies to the fractional nonlinear Schr\"{o}dinger equations (fNLS) with a potential ( see Section 5.4). To solve the uniqueness problem, we first show the non-degeneracy. As far as we know, these non-degeneracy ( see Lemmas \ref{lem4.13b}, \ref{lem4.26}) and uniqueness ( see Theorems \ref{thm4.9}, \ref{thm4.27}) results are completely new, even in the case of the Laplacian. They present interest by themselves.
		
		The rest of this paper is organized as follows. We present our general method and our main results in Sections \ref{AF1} - \ref{AF3}. We then discuss some applications for scalar field equations involving a fractional Laplacian with non-autonomous nonlinearities in Section \ref{frac}. A fractional Laplacian with mixed nonlinearities is considered in section \ref{mn}, p-laplacian in Section \ref{p}, and mixed fractional Laplacians in Section \ref{mfl}.

		\section{\textbf{Abstract Framework 1}} \label{AF1}
		
		In this section, we give the first part of our abstract framework to establish a $C^1$ global branch parameterized by $\lambda \in (a,b)$ under a non-degeneracy assumption. The readers will see that $a = -\infty$ under our setting in this section. Furthermore, once we show that the ground state is unique when $\lambda \rightarrow -\infty$, we can prove the uniqueness of the global branch. The existence of a $C^1$ global branch is enough to obtain the existence of normalized solutions. Nevertheless, the uniqueness of ground state has its own meaning and importance in various fields. Hence, this abstract framework 1 can also be viewed as a framework to prove the uniqueness of ground states. In particular, it is worth noting that the ODE techniques (which are essential to show the uniqueness in many cases) are not applicable in some important situations, such as non-local operators. We will exhibit concrete examples showing that our abstract framework plays an important role in these settings.
		
		Throughout this section, we assume
		
		$(H_1)$ $W$ is a Hilbert space, $D_{u}Q(u) = u$, the linear operator $D_{u}S + D_{u}G$ is self-adjoint and bounded below on a Hilbert space $E$ with operator domain $\widehat{W}$ and form domain $W$.
		
		Then $S(u) = \frac{1}{2}\left\langle D_{u}S(u),u\right\rangle, G(u) = \frac{1}{2}\left\langle D_{u}G(u),u\right\rangle, Q(u) = \frac{1}{2}\left\langle u,u\right\rangle$. Moreover, if $\lambda < \lambda_1 := \inf \sigma(D_{u}S + D_{u}G)$ where $\sigma(D_{u}S + D_{u}G)$ denotes the spectrum of $D_{u}S + D_{u}G: \widehat{W} \rightarrow E$, then
		$$
		\|u\|_{\lambda} = \sqrt{\left\langle D_{u}S(u) + D_{u}G(u) - \lambda u,u\right\rangle}
		$$
		is equivalent to the norm of $W$.
		
		Let $G$ be a topological group, and assume that the action of $G$ on $E$ is isometric. Define
		\begin{equation}
		E_G = \{u \in E: gu = u, \forall g \in G\}, \widehat{W}_G = \widehat{W} \cap E_G, W_G = W \cap E_G.
		\end{equation}
		We assume the following:
		
		$(H_2)$ $\Phi_{\lambda}(gu) = \Phi_{\lambda}(u), \forall g \in G$.
		
		$(H_3)$ $\Phi_{\lambda}|_{W_G}$ satisfies (PS) condition in $W_G$ for all $\lambda < \lambda_1$.
		
		$(H_4)$ $D_{uu}F(0) = 0$, $\left\langle D_{uu}F(u)v,w \right\rangle = \left\langle v,D_{uu}F(u)w \right\rangle$. And there exists some $p > 2$ such that
		\begin{equation}
		\left\langle D_{uu}F(u)u,u \right\rangle \geq (p-1)\left\langle D_uF(u),u \right\rangle > 0, \forall u \in W \setminus \{0\}.
		\end{equation}
	
	    By the principle of symmetric criticality (see \cite[Theorem 1.28]{Wi}), $(H_2)$ ensures that any critical point of $\Phi_\lambda$ restricted to $\widehat{W}_G$ is a critical point of $\Phi_\lambda$.  To show the existence of solutions, we define
	    $$
	    \mathcal{N}_{G,\lambda} = \{u \in W_G \setminus \{0\}: \left\langle D_uS(u) + D_uG(u) - \lambda u,u \right\rangle = \left\langle D_uF(u) ,u \right\rangle.
	    $$
	    \begin{definition}
	    	Set $h_G(\lambda) = \inf_{u \in \mathcal{N}_{G,\lambda}}\Phi_\lambda(u)$. We say that $u \in W_G \setminus \{0\}$ is a $G$-ground state solution if $u$ solves (\ref{eq1.1}) and achieves $h_G(\lambda)$.
	    \end{definition}
	
	    If $\mathcal{N}_{G,\lambda}$ is a $C^1$ manifold with codimension $1$ in $W_G$, $h_G(\lambda)$ is well-defined and $h_G(\lambda) \neq 0$ for all $\lambda < \lambda_1$, by $(H_3)$, standard arguments yield the existence of a $G$-ground state solution with $G$-Morse index $1$ (the definition will be given by \eqref{eqb.3}). Then, we establish a $C^1$ global branch in $\widehat{W}_G$ under the following assumptions:
	
	    $(H_5)$ If $\Phi_\lambda'(u_\lambda) = 0$ and $\mu_{G}(u_\lambda) = 1$, then $\ker D_{uu}\Phi_\lambda(u_\lambda)|_{E_G} = \{0\}$, where $\mu_{G}(u)$, the $G$-Morse index of $u$, is defined as
	    \begin{equation} \label{eqb.3}
	    \mu_{G}(u):= \sharp \{e < 0: e \ is \ an \ eigenvalue \ of \ D_{uu}\Phi_\lambda(u)|_{E_G}\}.
	    \end{equation}
	
	    $(H_6)$ If $u \in W$ solves (\ref{eq1.1}). Then $u \in \widehat{W}$.
	
	    $(H_7)$ $D_{u}F(u) \in E, \forall u \in W$, and $D_{uu}F(u)$ maps $\widehat{W}$ to $E$ for all $u \in \widehat{W}$.
	
	    Our main result in this section reads as follows:
	
	    \begin{theorem}
	    	\label{thmB.1} Assume that $(H_1)$ - $(H_7)$ hold, and that $\mathcal{N}_{G,\lambda}$ is a $C^1$ manifold with codimension $1$ in $W_G$, $h_G(\lambda)$ is well-defined and $h_G(\lambda) \neq 0$ for all $\lambda < \lambda_1$. Then there exists
	    	$$
	    	\lambda \mapsto u_\lambda \in C^1((-\infty,\lambda_1),W_G \setminus \{0\}),
	    	$$
	    	such that $u_\lambda$ solves (\ref{eq1.1}), and $\mu_{G}(u_\lambda) = 1$.
	    \end{theorem}
	
	    For the readers' convenience, the proof of our main result will be divided into three lemmas.
	
	    \begin{lemma}
	    	\label{lemB.10} Assume that $(H_1)$ and $(H_4)$ hold. If
	    	$$
	    	\lambda \mapsto u_\lambda \in C^1((\lambda^\ast,\lambda_\ast],W \setminus \{0\})
	    	$$
	    	for some $\lambda^\ast < \lambda_\ast < \lambda_1$ where $u_\lambda$ solves (\ref{eq1.1}). Then $u_\lambda$ is uniformly bounded with respect to $\lambda \in (\lambda^\ast,\lambda_\ast]$.
	    \end{lemma}
	
	    \textit{Proof.}  Since $u_{\lambda}$ solves (\ref{eq1.1}), we obtain
	    \begin{equation} \label{eqb.9}
	    \left\langle D_uS(u_\lambda) + D_uG(u_\lambda) - \lambda u_\lambda,u_\lambda \right\rangle = \left\langle D_uF(u_\lambda) ,u_\lambda \right\rangle.
	    \end{equation}
	    Hence, we have
	    \begin{equation}
	    (\lambda_1 - \lambda)\left\langle u_\lambda,u_\lambda \right\rangle \leq \left\langle D_uF(u_\lambda) ,u_\lambda \right\rangle.
	    \end{equation}
	    Differentiating $D_uS(u_\lambda) + D_uG(u_\lambda) = \lambda u_\lambda + D_uF(u_\lambda)$ with respect to $\lambda$, we get
	    \begin{equation}
	    D_uS(\partial_\lambda u_\lambda) + D_uG(\partial_\lambda u_\lambda) = u_\lambda + \lambda \partial_\lambda u_\lambda + D_{uu}F(u_\lambda)\partial_\lambda u_\lambda,
	    \end{equation}
	    i.e.
	    \begin{equation}
	    D_{uu}\Phi(u_\lambda)\partial_\lambda u_\lambda = u_\lambda.
	    \end{equation}
	    Set
	    $$
	    g(\lambda) = \left\langle D_uF(u_\lambda),u_\lambda \right\rangle - 2F(u_\lambda).
	    $$
	    Note that
	    \begin{equation}
	    D_{uu}\Phi(u_\lambda)u_\lambda = D_uF(u_\lambda) - D_{uu}F(u_\lambda)u_\lambda.
	    \end{equation}
	    Then
	    \begin{eqnarray}
	    g'(\lambda) &=& \left\langle D_{uu}F(u_\lambda)\partial_\lambda u_\lambda,u_\lambda \right\rangle - \left\langle D_uF(u_\lambda),\partial_\lambda u_\lambda \right\rangle \nonumber \\
	    &=& \left\langle D_{uu}\Phi(u_\lambda)^{-1}(D_{uu}F(u_\lambda)u_\lambda - D_uF(u_\lambda)),u_\lambda \right\rangle \nonumber \\
	    &=& -\left\langle u_\lambda,u_\lambda \right\rangle \nonumber \\
	    &\geq& -\frac{1}{\lambda_1 - \lambda}\left\langle D_uF(u_\lambda) ,u_\lambda \right\rangle \nonumber \\
	    &\geq& -\frac{p}{p-2}\frac{1}{\lambda_1 - \lambda}g(\lambda).
	    \end{eqnarray}
	    Let
	    $$
	    \widetilde{g}(\lambda) = (\lambda_1 - \lambda)^{-\frac{p}{p-2}}g(\lambda), \lambda \in (\lambda^\ast,\lambda_\ast].
	    $$
	    Then
	    $$
	    \widetilde{g}'(\lambda) \geq 0,
	    $$
	    yielding that $\widetilde{g}(\lambda) \leq \widetilde{g}(\lambda_\ast)$ for any $\lambda \in (\lambda^\ast,\lambda_\ast]$ and hence
	    $$
	    g(\lambda) \leq (\frac{\lambda_1 - \lambda^\ast}{\lambda_1 - \lambda_\ast})^{\frac{p}{p-2}}g(\lambda_\ast).
	    $$
	    By (\ref{eqb.9}) and $(H_4)$, $g(\lambda) \sim \|u_\lambda\|_\lambda \sim \|u_\lambda\|$. Thus $u_\lambda$ is uniformly bounded with respect to $\lambda \in (\lambda^\ast,\lambda_\ast]$.
	    \qed\vskip 5pt
	
	    \begin{lemma}
	    	\label{lemB.11} Assume that $(H_1)$ - $(H_4)$ hold and that
	    	$$
	    	\lambda \mapsto u_\lambda \in C^1((\lambda^\ast,\lambda_\ast],W_G \setminus \{0\})
	    	$$
	    	for some $\lambda^\ast < \lambda_\ast < \lambda_1$ where $u_\lambda$ is a solution of (\ref{eq1.1}). Let $\{\lambda_n\} \subset (\lambda^\ast,\lambda_\ast)$ be a sequence such that $\lambda_n \rightarrow \lambda^\ast$, and $\mu_{G}(u_{\lambda_n}) = 1$. Then, after possibly passing to a subsequence, we have $u_{\lambda_n} \rightarrow u_{\lambda^\ast}$ in $W \setminus \{0\}$ where $\mu_{G}(u_{\lambda^\ast}) = 1$ and $u_{\lambda^\ast}$ satisfies
	    	\begin{equation} \label{eqb.15}
	    	D_u\Phi_{\lambda^\ast}(u_{\lambda^\ast}) = 0.
	    	\end{equation}
	    \end{lemma}
	
	    \textit{Proof.} From Lemma \ref{lemB.10}, we obtain the boundedness of $u_{\lambda_n}$ in $W$ and the boundedness of $\left\langle u_{\lambda_n},u_{\lambda_n} \right\rangle$. Then we have that $\Phi_{\lambda^\ast}(u_{\lambda_n}) = \Phi_{\lambda_n}(u_{\lambda_n}) + \frac{1}{2}(\lambda_n - \lambda^\ast)\left\langle u_{\lambda_n},u_{\lambda_n} \right\rangle$ is bounded and
	    $$
	    D_u\Phi_{\lambda^\ast}(u_{\lambda_n}) = D_u\Phi_{\lambda_n}(u_{\lambda_n}) + (\lambda_n - \lambda^\ast)u_{\lambda_n} \rightarrow 0.
	    $$
	    Since $\Phi_{\lambda^\ast}$ satisfies (PS) condition in $W_G$, up to a subsequence, there exists $u_{\lambda^\ast} \in W_G$ such that $u_{\lambda_n} \rightarrow u_{\lambda^\ast}$ in $W$ and $u_{\lambda^\ast}$ satisfies (\ref{eqb.15}).
	
	    We claim that $u_{\lambda^\ast} \neq 0$. If it is not the case, $D_{uu}\Phi_{\lambda^\ast}(0) = D_uS + D_uG - \lambda^\ast$ is invertible. Since $0$ is always a solution of (\ref{eq1.1}) for any $\lambda$, the implicit function theorem yields that $u_{\lambda_n}$ must be $0$ for $n$ large enough, contradicting the fact that $u_{\lambda_n} \neq 0$.
	
	    Finally, we show that $\mu_{G}(u_{\lambda^\ast}) = 1$. On the one hand, $\mu_{G}(u_{\lambda^\ast}) \geq 1$ since
	    $$
	    \langle D_{uu}\Phi_{\lambda^\ast}(u_{\lambda^\ast})u_{\lambda^\ast},u_{\lambda^\ast}\rangle = \langle D_{u}F(u_{\lambda^\ast}) - D_{uu}F(u_{\lambda^\ast})u_{\lambda^\ast},u_{\lambda^\ast}\rangle < 0.
	    $$
	    One the other hand, since the Morse index is lower semi-continuous with respect to the norm-resolvent topology, we conclude that
	    $$
	    \mu_{G}(u_{\lambda^\ast}) \leq \liminf_{n \rightarrow +\infty}\mu_{G}(u_{\lambda_{n}}) = 1,
	    $$
	    which completes the proof.
	    \qed\vskip 5pt
	    	
	    \begin{lemma}
	    	\label{lemB.12} Assume that $(H_1)$ - $(H_7)$ hold, and that for some $\lambda_\ast <\lambda_1$, there exists $u_{\lambda_\ast} \in W_G \setminus \{0\}$ such that $u_{\lambda_\ast}$ solves (\ref{eq1.1}), and $\mu_{G}(u_{\lambda_\ast}) = 1$. Then there exists
	    	$$
	    	\lambda \mapsto u_\lambda \in C^1((-\infty,\lambda_\ast],W_G \setminus \{0\}),
	    	$$
	    	such that $u_\lambda$ solves (\ref{eq1.1}), and $\mu_{G}(u_\lambda) = 1$.
	    \end{lemma}
	    	
	    \textit{Proof.}  We use an implicit function argument for the map
	    $$
	    \Psi: \widehat{W}_G \times (-\infty,\lambda_1) \rightarrow \widehat{W}_G, \Psi(u,\lambda) = u - (D_{u}S + D_{u}G - \lambda)^{-1}D_{u}F(u).
	    $$
	    It is not difficult to see that $\Psi$ is a well-defined map of class $C^{1}$. By $(H_6)$, $u \in W_G$ solves (\ref{eq1.1}) if and only if $\Psi(u,\lambda) =0$. $u_{\lambda_\ast}$ solves (\ref{eq1.1}), thus $\Psi(u_{\lambda_\ast},\lambda_\ast) = 0$. Next, we consider the Fr\'{e}chet derivative
	    \begin{equation}
	    \Psi_{u}(u_{\lambda_\ast},\lambda_\ast) = 1 - K, K = (D_{u}S + D_{u}G - \lambda_\ast)^{-1}D_{uu}F(u_{\lambda_\ast}).
	    \end{equation}
	    By $(H_7)$, $D_{uu}F(u_{\lambda_\ast})$ maps $\widehat{W}_G$ to $E_{G}$. Hence, $K$ maps $\widehat{W}_G$ to $\widehat{W}_G$. If $v - Kv = 0$ for some $v \in \widehat{W}_G$, then $D_{uu}\Phi_{\lambda_\ast}(u_{\lambda_\ast})v = 0$. By $(H_5)$,
	    $$
	    \ker D_{uu}\Phi_{\lambda_\ast}(u_{\lambda_\ast})|_{E_{G}} = 0.
	    $$
	    Thus we deduce that $v = 0$ and thus $1 - K : \widehat{W}_G \rightarrow \widehat{W}_G$ is injective. By standard arguments, the operator $D_{uu}\Phi_{\lambda_\ast}(u_{\lambda_\ast})$ has a bounded inverse on $E_{G}$ since $D_{uu}\Phi_{\lambda_\ast}(u_{\lambda_\ast})$ has trivial kernel on $E_{G}$. Then for any $w \in \widehat{W}_G$, $D_{uu}\Phi_{\lambda_\ast}(u_{\lambda_\ast})^{-1}(D_{u}S + D_{u}G - \lambda_\ast)w \in \widehat{W}_G$. Therefore,
	    \begin{eqnarray}
	    && (1-K)D_{uu}\Phi_{\lambda_\ast}(u_{\lambda_\ast})^{-1}Bw \nonumber \\
	    &=& B^{-1}D_{uu}\Phi_{\lambda_\ast}(u_{\lambda_\ast})D_{uu}\Phi_{\lambda_\ast}(u_{\lambda_\ast})^{-1}Bw \nonumber \\
	    &=& w,
	    \end{eqnarray}
	    where $B = D_{u}S + D_{u}G - \lambda_\ast$, showing that $1 - K : \widehat{W}_G \rightarrow \widehat{W}_G$ is surjective. Using the implicit function theorem for $\Psi(u,\lambda)$ at $(u_{\lambda_\ast},\lambda_\ast)$, we deduce that for some $\delta > 0$, there exists a map $u(\lambda) \in C^{1}((\lambda_\ast-\delta,\lambda_\ast], \widehat{W}_G \setminus \{0\})$ such that the following holds true:
	
	    $(i)$ $u(\lambda)$ solves $\Psi(u,\lambda) = 0$ for all $\lambda \in (\lambda_\ast-\delta,\lambda_\ast]$;
	
	    $(ii)$ there exists $\epsilon > 0$ such that $u(\lambda)$ is the unique solution of $\Psi(u,\lambda) = 0$ for $\lambda \in (\lambda_\ast-\delta,\lambda_\ast]$ in the neighborhood $\{u \in \widehat{W}_G: \|u - u(\lambda_\ast)\|_{\widehat{W}} < \epsilon\}$; in particular, we have that $u(\lambda_\ast) = u_{\lambda_\ast}$ holds.
	
	    Note that $D_{uu}\Phi_{\lambda}(u_{\lambda}) \rightarrow D_{uu}\Phi_{\lambda_\ast}(u_{\lambda_\ast})$ in the norm-resolvent sense as $\lambda \rightarrow \lambda_\ast$. Then, along the line of the proof of Lemma \ref{lem4.16} $(iii)$ below, we can deduce that $\mu_{G}(u(\lambda)) = 1$.
	
	    We claim that this branch $u(\lambda)$ can be extend to $-\infty$. Otherwise, if this branch $u(\lambda)$ can only be extend to $\lambda^\ast$, where $\mu_{G}(u(\lambda)) = 1$ for $\lambda \in (\lambda^\ast,\lambda_\ast]$. By Lemma \ref{lemB.11}, there exists $u_{\lambda^\ast} \in W_G \setminus \{0\}$ such that $\mu_{G}(u_{\lambda^\ast}) = 1$ and $u_{\lambda^\ast}$ satisfies (\ref{eqb.15}). By $(H_5)$, $D_{uu}\Phi_{\lambda^\ast}(u_{\lambda^\ast})|_{E_{G}} = \{0\}$. Then similar implicit function arguments yields that this branch $u(\lambda)$ can be extend to $\lambda^\ast - \delta < \lambda^\ast$, in a contradiction with the minimality of $\lambda^\ast$. Thus we complete the proof.
	    \qed\vskip 5pt
	
	    \textbf{Completion of the Proof for Theorem \ref{thmB.1}.}  Note that $\mathcal{N}_{G,\lambda}$ is a $C^1$ manifold with codimension $1$ in $W_G$, $h_G(\lambda)$ is well-defined and $h_G(\lambda) \neq 0$ for all $\lambda \neq \lambda_1$. Then there exists a (PS) sequence $\{u_n\} \subset W_G$ at level $h_G(\lambda)$. By $(H_3)$, there exists $u_\lambda \in W_G$ such that $u_n \rightarrow u_\lambda$ in $W$. Hence, $\Phi_\lambda'(u_\lambda) = 0$ and $\Phi_\lambda(u_\lambda) = h_G(\lambda)$, showing the existence of a nontrivial $G$-ground state solution. Next we prove that $\mu_G(u_\lambda) = 1$. On the one hand, from
	    \begin{equation}
	    \left\langle D_{uu}\Phi_\lambda(u_\lambda)u_\lambda,u_\lambda\right\rangle = \left\langle D_{u}F(u_\lambda) - D_{uu}F(u_\lambda)u_\lambda,u_\lambda\right\rangle < 0,
	    \end{equation}
	    we derive that $\mu_{G}(u_\lambda) \geq 1$. On the other hand, noticing that the codimension of $\mathcal{N}_{\lambda}$ is $1$ and $u_\lambda$ achieves the minimum of $h(\lambda)$, we have that $\mu_G(u_\lambda) \leq 1$. Hence, $\mu_G(u_\lambda) = 1$. Then by Lemma \ref{lemB.12}, we can take $\lambda_\ast < \lambda_1$ such that there exists
	    $$
	    \lambda \mapsto u_\lambda \in C^1((-\infty,\lambda_\ast],W_G \setminus \{0\}),
	    $$
	    and $u_\lambda$ solves (\ref{eq1.1}), $\mu_{G}(u_\lambda) = 1$. Noticing that $\lambda_\ast$ can be chosen arbitrarily close to $\lambda_1$, we know that there exists
	    $$
	    \lambda \mapsto u_\lambda \in C^1((-\infty,\lambda_1),W_G \setminus \{0\}),
	    $$
	    such that $u_\lambda$ solves (\ref{eq1.1}), and $\mu_{G}(u_\lambda) = 1$.
	    \qed\vskip 5pt
	
	    \begin{corollary} \label{cor uni}
	    	(Uniqueness)  Under the hypotheses of Theorem \ref{thmB.1}. If we also assume that there exists $\Lambda_1 < \lambda_1$ such that (\ref{eq1.1}) admits a unique solution with $G$-Morse index $1$ for all $\lambda < \Lambda_1$. Then for all $\lambda < \lambda_1$, (\ref{eq1.1}) admits a unique solution with $G$-Morse index $1$.
	    	
	    	In particular, (\ref{eq1.1}) admits a unique $G$-ground state solution for all $\lambda < \lambda_1$.
	    	
	    	Furthermore, if all ground states of (\ref{eq1.1}) belong to $W_G$. Then a $G$-ground state is indeed a ground state, and (\ref{eq1.1}) admits a unique ground state solution for all $\lambda < \lambda_1$.
	    \end{corollary}	
	
	    \textit{Proof.} The existence of a solution with $G$-Morse index $1$ has been established in the proof of Theorem \ref{thmB.1}. Next we will show its uniqueness. Suppose on the contrary that (\ref{eq1.1}) has two solutions $u_{\lambda_0}, \widetilde{u}_{\lambda_0}$ with $G$-Morse index $1$ in $W_G$ for some $\lambda_0 < \lambda_1$. By the proof of Theorem \ref{thmB.1}, we obtain two global branches $u(\lambda), \widetilde{u}(\lambda), \lambda \in [\Lambda_1 - 1,\lambda_0]$ stemming from $u(\lambda_0) = u_{\lambda_0}, \widetilde{u}(\lambda_0) = \widetilde{u}_{\lambda_0}$ respectively. Furthermore, both $u(\Lambda_1-1)$ and $\widetilde{u}(\Lambda_1-1)$ are solutions of (\ref{eq1.1}) with $\lambda = \Lambda_1-1$ and $\mu_{G}(u(\Lambda_1-1)) = \mu_{G}(\widetilde{u}(\Lambda_1-1)) = 1$. The uniqueness when $\lambda < \Lambda_1$ implies that $u(\Lambda_1-1) = \widetilde{u}(\Lambda_1-1)$. By $(H_5)$, $u(\Lambda_1-1)$ is non-degenerate in $W_G$. However, there are two different local branches $u(\lambda)$ and $\widetilde{u}(\lambda)$, $\lambda \in [\Lambda_1-1,\Lambda_1-1 + \delta)$ stemming form $u(\Lambda_1-1)$, contradicting the non-degeneracy of $u(\Lambda_1-1)$. Hence the proof is complete.
	    \qed\vskip 5pt
		
		\section{\textbf{Abstract Framework 2}} \label{AF2}
		
		In Section \ref{AF1}, we established a $C^1$ global branch $u_\lambda$ parameterized by $\lambda \in (-\infty,\lambda_1)$. To get the existence of normalized solutions, we need abstract frameworks to show the behaviors of $Q(u_\lambda)$ as $\lambda \rightarrow -\infty$ and $\lambda \rightarrow \lambda_1$. This section and the next one provide two different approaches. In this section, we merely need that $\Phi_{\lambda} \in C^{1}$.
		
		Throughout this section, we always assume that $\lambda$, $\lambda_{0} \in (a, b)$, $\mathcal{N}_{\lambda}$, $\mathcal{K}_{\lambda}$, $h(\lambda)$, $d(\lambda)$ and $\hat{d}(\lambda)$ are given in Section 1, and we always assume that $h(\lambda) \neq 0$, ensuring that $0 \notin \mathcal{K}_{\lambda}$. Therefore, all elements in $\mathcal{K}_{\lambda}$ are nontrivial solutions of (\ref{eq1.1}). Our goals are to obtain the continuity of $d(\lambda)$ if $d(\lambda) \equiv \hat{d}(\lambda)$ and to study the asymptotical behaviors of $d(\lambda)$ and $\hat{d}(\lambda)$. It is worth noting that the abstract framework 2 in this section is not just a simple continuation of abstract framework 1, but a framework in its own right. For example, classical arguments are invalid to get the existence of normalized solutions of the equation $\Delta u + \lambda u + |u|^{p-2}u = 0$ on $\mathbb{R}^N \setminus B_1$, particularly in the $L^2$-supercritical case. Relying on ODE techniques, the uniqueness of positive radial solution can be proved, ensuring that $d(\lambda) \equiv \hat{d}(\lambda)$ (in fact, the uniqueness is enough to obtain a continuous global branch $u_\lambda$ parameterized by $\lambda$). Therefore, our framework 2 is applicable. It helps to obtain threshold results for the existence, non-existence and multiplicity of normalized solutions for semi-linear elliptic equations in the exterior of a ball, see \cite{Song3}. In some cases, even if we do not know the uniqueness of ground state, the abstract framework 2 can be used, see Section \ref{PII} for an example with a Hardy potential.
		
		First, we study the properties of $h(\lambda)$, $d(\lambda)$ and $\hat{d}(\lambda)$. To do this, we assume
		
		$(F_{1})$ For any $(\lambda, u) \in (a, b) \times W \backslash \{0\}$, there exists an unique function $t: (a, b) \times W \backslash \{0\} \rightarrow (0, \infty)$ such that
		\begin{equation}
			when \ t > 0, \ then \ (\lambda, t u) \in \mathcal{N} \Leftrightarrow t = t(\lambda, u),\nonumber
		\end{equation}
		and $\Phi_{\lambda}(t(\lambda, u)u) = \max \limits_{t > 0}\Phi_{\lambda}(tu)$, $t \in C((a, b) \times W \backslash \{0\}, (0, \infty))$, where $\mathcal{N} := \cup_{a < \lambda < b}\mathcal{N}_{\lambda}$.
		
		$(F_{2})$ $Q(u_{\lambda})$ is bounded in $\mathbb{R}$ and $D_{u}Q(u_{\lambda})$ is bounded in $W^{\ast}$ whenever $u_{\lambda} \in \mathcal{K}_{\lambda}$ and $\lambda \rightarrow \lambda_{0}$.
		
		$(F_{3})$ $\Phi_{\lambda_{0}}$ satisfies $(PS)_{c}$ condition for $c \leq h(\lambda_{0})$, $\forall \lambda_{0} \in (a,b)$.
		
		$(F_{4})$ $\liminf_{\lambda \rightarrow \lambda_{0}}h(\lambda) > 0$, $\forall \lambda_{0} \in (a,b)$.
		
		\begin{theorem}
			\label{thm2.1}  Suppose $(F_{1})$ - $(F_{4})$.
			
			$(i)$ $d(\lambda_{0}) \leq \liminf_{\lambda \rightarrow \lambda_{0}}d(\lambda)$, $\hat{d}(\lambda_{0}) \geq \limsup_{\lambda \rightarrow \lambda_{0}}\hat{d}(\lambda)$.
			
			$(ii)$ $\lim_{\lambda \rightarrow \lambda_{0}}h(\lambda) = h(\lambda_{0})$.
			
			$(iii)$ $h(\lambda)$ is differentiable at almost everywhere $\lambda \in (a, b)$. Furthermore,
			$$
			\lim \limits_{\lambda \rightarrow \lambda_{0}^{+}}\frac{h(\lambda) - h(\lambda_{0})}{\lambda - \lambda_{0}} = -\hat{d}(\lambda_{0}), \lim \limits_{\lambda \rightarrow \lambda_{0}^{-}}\frac{h(\lambda) - h(\lambda_{0})}{\lambda - \lambda_{0}} = -d(\lambda_{0}).
			$$
			In particular, $d(\lambda_{0}) = \hat{d}(\lambda_{0}) \Leftrightarrow h(\lambda)$ is differentiable at $\lambda_{0}$ and $h'(\lambda_{0}) = -d(\lambda_{0})$.
		\end{theorem}
		
		\textit{Proof.  } Step 1: $\limsup_{\lambda \rightarrow \lambda_{0}} h(\lambda) \leq h(\lambda_{0})$.
		
		Take $u_{\lambda_{0}} \in \mathcal{K}_{\lambda_{0}}$, then $\Phi_{\lambda_{0}}(u_{\lambda_{0}}) = h(\lambda_{0})$. Consider $t(\lambda, u_{\lambda_{0}})u_{\lambda_{0}} \in \mathcal{N}_{\lambda}$. By $(F_{1})$,
		$$
		t(\lambda, u_{\lambda_{0}})u_{\lambda_{0}} \rightarrow t(\lambda_{0}, u_{\lambda_{0}})u_{\lambda_{0}} = u_{\lambda_{0}} \ in \ W \ as \ \lambda \rightarrow \lambda_{0},
		$$
		yielding that
		\begin{eqnarray} \label{eq3.1}
			h(\lambda) &\leq& \Phi_{\lambda}(t(\lambda, u_{\lambda_{0}})u_{\lambda_{0}}) \nonumber \\
			&=& \Phi_{\lambda_{0}}(t(\lambda, u_{\lambda_{0}})u_{\lambda_{0}}) + (\lambda_{0} - \lambda)Q(t(\lambda, u_{\lambda_{0}})u_{\lambda_{0}}) \nonumber \\
			&=& \Phi_{\lambda_{0}}(u_{\lambda_{0}}) + o(1).
		\end{eqnarray}
		Letting $\lambda \rightarrow \lambda_{0},$ we have $\limsup_{\lambda \rightarrow \lambda_{0}} h(\lambda) \leq h(\lambda_{0})$.
		
		Step 2: If $\lambda_{n} \rightarrow \lambda_{0}$, $u_{\lambda_{n}} \in \mathcal{K}_{\lambda_{n}}$, then $u_{\lambda_{n}} \rightarrow u_{\lambda_{0}}$ in $W$ up to a subsequence and $u_{\lambda_{0}} \in \mathcal{K}_{\lambda_{0}}$.
		
		By $(F_{2})$, $(\lambda_{n} - \lambda_{0})Q(u_{\lambda_{n}}) \rightarrow 0$. Noticing that
		$$
		\Phi_{\lambda_{0}}(u_{\lambda_{n}}) = \Phi_{\lambda_{n}}(u_{\lambda_{n}}) + (\lambda_{n} - \lambda_{0})Q(u_{\lambda_{n}}) = h(\lambda_{n}) + o(1),
		$$
		we obtain that $\limsup_{n \rightarrow \infty}\Phi_{\lambda_{0}}(u_{\lambda_{n}}) \leq h(\lambda_{0})$. By $(F_{4})$, we may assume that $\Phi_{\lambda_{0}}(u_{\lambda_{n}}) \rightarrow c \in (0,h(\lambda_{0})]$ passing to a subsequence if necessary. Furthermore, since $D_{u}Q(u_{\lambda_{n}})$ is bounded in $W^{\ast}$, $D_{u}\Phi_{\lambda_{0}}(u_{\lambda_{n}}) = D_{u}\Phi_{\lambda_{n}}(u_{\lambda_{n}}) + (\lambda_{n} - \lambda_{0})D_{u}Q(u_{\lambda_{n}}) \rightarrow 0$ as $n \rightarrow \infty$.
		
		By $(F_{3})$, we may assume that $u_{\lambda_{n}} \rightarrow u_{\lambda_{0}}$ up to a subsequence. Using $(F_{4})$, $\Phi_{\lambda_{0}}(u_{\lambda_{0}}) > 0$, implying that $u_{\lambda_{0}} \neq 0$. Since $\Phi_{\lambda_{0}} \in C^{1}$, $D_{u}\Phi_{\lambda_{0}}(u_{\lambda_{0}}) = 0$. Thus $u_{\lambda_{0}} \in N_{\lambda_{0}}$.
		
		Next, we will show that $u_{\lambda_{0}} \in \mathcal{K}_{\lambda_{0}}$. $\forall \epsilon > 0$, $\exists X > 0$, $\forall n > X$, there holds that
		\begin{eqnarray} \label{eq3.2}
			h(\lambda_{0}) &\leq& \Phi_{\lambda_{0}}(u_{\lambda_{0}}) \nonumber \\
			&\leq& \Phi_{\lambda_{0}}(u_{\lambda_{n}}) + \frac{\epsilon}{2} \nonumber \\
			&=& \Phi_{\lambda_{n}}(u_{\lambda_{n}}) + (\lambda_{n} - \lambda_{0})Q(u_{\lambda_{n}}) + \frac{\epsilon}{2} \nonumber \\
			&\leq& \Phi_{\lambda_{n}}(u_{\lambda_{n}}) + \epsilon.
		\end{eqnarray}
		Letting $n \rightarrow \infty,$ we get that $h(\lambda_{0}) \leq \limsup_{\lambda \rightarrow \lambda_{0}}h(\lambda)$. Together with Step 1, we can conclude that $\limsup_{\lambda \rightarrow \lambda_{0}}h(\lambda) = h(\lambda_{0})$. We can derive that $\Phi_{\lambda_{0}}(u_{\lambda_{0}}) = h(\lambda_{0})$ by (\ref{eq3.2}) immediately. Thus $u_{\lambda_{0}} \in \mathcal{K}_{\lambda_{0}}$.
		
		Step 3: The Proof of Theorem \ref{thm2.1} $(i)$ and $(ii)$.
		
		Take $\{\lambda_{n}\}_{n = 1}^{+\infty}$ such that $\lambda_{n} \rightarrow \lambda_{0}$, $\liminf_{\lambda \rightarrow \lambda_{0}}d(\lambda) = \lim_{n \rightarrow +\infty}d(\lambda_{n})$. Since $\mathcal{K}_{\lambda}$ is compact, we can take $u_{\lambda_{n}} \in \mathcal{K}_{\lambda_{n}}$ such that $d(\lambda_{n}) = Q(u_{\lambda_{n}})$. By Step 2, passing to a subsequence if necessary, $u_{\lambda_{n}} \rightarrow u_{\lambda_{0}} \in \mathcal{K}_{\lambda_{0}}$. Therefore,
		\begin{eqnarray}
			d(\lambda_{0}) &\leq& Q(u_{\lambda_{0}}) \nonumber \\
			&=& \lim_{n \rightarrow \infty}Q(u_{\lambda_{n}}) \nonumber \\
			&=& \lim_{n \rightarrow \infty} d(\lambda_{n}) \nonumber \\
			&=& \liminf_{\lambda \rightarrow \lambda_{0}}d(\lambda). \nonumber
		\end{eqnarray}
		
		Similarly, we can take $\{\lambda_{n}\}_{n = 1}^{+\infty}$ such that $\lambda_{n} \rightarrow \lambda_{0}$, $\limsup_{\lambda \rightarrow \lambda_{0}}\hat{d}(\lambda) = \lim_{n \rightarrow \infty}\hat{d}(\lambda_{n})$ and take $u_{\lambda_{n}} \in \mathcal{K}_{\lambda_{n}}$ such that $\hat{d}(\lambda_{n}) = Q(u_{\lambda_{n}})$. Then
		\begin{eqnarray}
			\hat{d}(\lambda_{0}) &\geq& Q(u_{\lambda_{0}}) \nonumber \\
			&=& \lim_{n \rightarrow \infty}Q(u_{\lambda_{n}}) \nonumber \\
			&=& \lim_{n \rightarrow \infty}\hat{d}(\lambda_{n}) \nonumber \\
			&=& \limsup_{\lambda \rightarrow \lambda_{0}}\hat{d}(\lambda). \nonumber
		\end{eqnarray}
		
		Next we prove $(ii)$. According to Step 2, for any $\{\lambda_{n}\}$ with $\lambda_{n} \rightarrow \lambda_{0}$, passing to a subsequence if necessary, there exists $u_{\lambda_{n}} \in K_{\lambda_{n}}$ s.t. $\lim_{n \rightarrow \infty}h(\lambda_{n}) = \lim_{n \rightarrow \infty}\Phi_{\lambda_{n}}(u_{\lambda_{n}}) = h(\lambda_{0})$. Hence, we conclude that $\lim_{\lambda \rightarrow \lambda_{0}}h(\lambda) = h(\lambda_{0}),$ which completes the proof.
		
		Step 4: Suppose that $u_{\lambda} \in \mathcal{K}_{\lambda}$, $u_{\lambda_{0}} \in \mathcal{K}_{\lambda_{0}}$. Let \(\widetilde{u}_{\lambda} = t(\lambda_{0}, u_{\lambda})u_{\lambda}\) and \(\widetilde{u}_{\lambda_{0}} = t(\lambda, u_{\lambda_{0}})u_{\lambda_{0}}\), then
		\begin{equation} \label{eq3.3}
			\Phi_{\lambda}(u_{\lambda}) \leq (\lambda_{0}-\lambda)Q(\widetilde{u}_{\lambda_{0}}) + \Phi_{\lambda_{0}}(u_{\lambda_{0}}),
		\end{equation}
		\begin{equation} \label{eq3.4}
			\Phi_{\lambda_{0}}(u_{\lambda_{0}}) \leq (\lambda-\lambda_{0})Q(\widetilde{u}_{\lambda}) + \Phi_{\lambda}(u_{\lambda}).
		\end{equation}
		
		By $(F_{1})$, $\Phi_{\lambda}(u_{\lambda}) = \max \limits_{t > 0}\Phi_{\lambda}(tu_{\lambda})$ and $\Phi_{\lambda_{0}}(u_{\lambda_{0}}) = \max \limits_{t > 0}\Phi_{\lambda_{0}}(tu_{\lambda_{0}})$. Therefore,
		\begin{equation} \nonumber
			\Phi_{\lambda}(\widetilde{u}_{\lambda}) \leq \Phi_{\lambda}(u_{\lambda}), \Phi_{\lambda_{0}}(\widetilde{u}_{\lambda_{0}}) \leq \Phi_{\lambda_{0}}(u_{\lambda_{0}}).
		\end{equation}
		Since $u_{\lambda} \in \mathcal{K}_{\lambda}$, $\widetilde{u}_{\lambda_{0}} \in \mathcal{N}_{\lambda},$ $u_{\lambda_{0}} \in \mathcal{K}_{\lambda_{0}}$, and $\widetilde{u}_{\lambda} \in N_{\lambda_{0}}$, then
		\begin{equation} \nonumber
			\Phi_{\lambda}(\widetilde{u}_{\lambda_{0}}) \geq \Phi_{\lambda}(u_{\lambda}), \Phi_{\lambda_{0}}(\widetilde{u}_{\lambda}) \geq \Phi_{\lambda_{0}}(u_{\lambda_{0}}).
		\end{equation}
		It follows that:
		\begin{eqnarray} \label{eq3.5}
			\Phi_{\lambda}(u_{\lambda}) &\leq& \Phi_{\lambda}(\widetilde{u}_{\lambda_{0}}) \nonumber \\
			&=& (\lambda_{0} - \lambda)Q(\widetilde{u}_{\lambda_{0}}) + \Phi_{\lambda_{0}}(\widetilde{u}_{\lambda_{0}}) \nonumber \\
			&\leq& (\lambda_{0} - \lambda)Q(\widetilde{u}_{\lambda_{0}}) + \Phi_{\lambda_{0}}(u_{\lambda_{0}}),
		\end{eqnarray}
		\begin{eqnarray} \label{eq3.6}
			\Phi_{\lambda_{0}}(u_{\lambda_{0}}) &\leq& \Phi_{\lambda_{0}}(\widetilde{u}_{\lambda}) \nonumber \\
			&=& (\lambda - \lambda_{0})Q(\widetilde{u}_{\lambda}) + \Phi_{\lambda}(\widetilde{u}_{\lambda}) \nonumber \\
			&\leq& (\lambda - \lambda_{0})Q(\widetilde{u}_{\lambda}) + \Phi_{\lambda}(u_{\lambda}).
		\end{eqnarray}
		
		Step 5: The Proof of Theorem \ref{thm2.1} $(iii)$.
		
		By Step 4,
		\begin{equation} \label{eq3.7}
			(\lambda_{0} - \lambda)Q(\widetilde{u}_{\lambda}) \leq h(\lambda) - h(\lambda_{0}) \leq (\lambda_{0}-\lambda)Q(\widetilde{u}_{\lambda_{0}}),
		\end{equation}
		where $\widetilde{u}_{\lambda}$ and $\widetilde{u}_{\lambda_{0}}$ are given by Step 4. Thus $h(\lambda)$ is monotonically decreasing in $(a, b)$. Then we know that $h(\lambda)$ is differentiable at almost every $\lambda \in (a, b)$.
		
		Next, we want to estimate
		\[
		\frac{h(\lambda)-h(\lambda_{0})}{\lambda - \lambda_{0}}.
		\]
		First, we assume that \(\lambda > \lambda_{0}\). It follows that:
		\begin{equation} \label{eq3.8}
			-Q(\widetilde{u}_{\lambda}) \leq \frac{h(\lambda) - h(\lambda_{0})}{\lambda - \lambda_{0}} \leq -Q(\widetilde{u}_{\lambda_{0}}).
		\end{equation}
		On the one hand,
		\begin{eqnarray} \label{eq3.9}
			\limsup_{\lambda \rightarrow \lambda_{0}^{+}}\frac{h(\lambda) - h(\lambda_{0})}{\lambda - \lambda_{0}} &\leq& -\liminf_{\lambda \rightarrow \lambda_{0}^{+}}Q(\widetilde{u}_{\lambda_{0}}) \nonumber \\
			&=& -\lim_{\lambda \rightarrow \lambda_{0}^{+}} Q(t(\lambda,u_{\lambda_{0}})u_{\lambda_{0}}) \nonumber \\
			&=& -Q(u_{\lambda_{0}}).
		\end{eqnarray}
		Using the arbitrariness of $u_{\lambda_{0}}$, we know that
		\begin{eqnarray} \label{eq3.10}
			\limsup_{\lambda \rightarrow \lambda_{0}^{+}}\frac{h(\lambda) - h(\lambda_{0})}{\lambda - \lambda_{0}} &\leq& \min_{u_{\lambda_{0}} \in \mathcal{K}_{\lambda_{0}}}(-Q(u_{\lambda_{0}})) \nonumber \\
			&=& -\max \limits_{u_{\lambda_{0}} \in \mathcal{K}_{\lambda_{0}}}Q(u_{\lambda_{0}}) \nonumber \\
			&=&  -\hat{d}(\lambda_{0}).
		\end{eqnarray}
		On the other hand,
		\begin{equation} \label{eq3.11}
			\liminf_{\lambda \rightarrow \lambda_{0}^{+}}\frac{h(\lambda) - h(\lambda_{0})}{\lambda - \lambda_{0}} \geq \liminf_{\lambda \rightarrow \lambda_{0}^{+}}(-Q(\widetilde{u}_{\lambda})) = -\limsup_{\lambda \rightarrow \lambda_{0}^{+}}Q(\widetilde{u}_{\lambda}).
		\end{equation}
		Take $\{\lambda_{n}\}_{n = 1}^{+\infty}$ such that $\lambda_{n} \rightarrow \lambda_{0}^{+}$ (i.e. $\lambda_{n} \geq \lambda_{0}$ and $\lim_{n \rightarrow \infty}\lambda_{n} = \lambda_{0}$) and
		$$
		\lim_{n \rightarrow \infty}Q(\widetilde{u}_{\lambda_{n}}) = \limsup_{\lambda \rightarrow \lambda_{0}^{+}}Q(\widetilde{u}_{\lambda}).
		$$
		By Step 2, passing to a subsequence if necessary, there exists some $u_{0} \in \mathcal{K}_{\lambda_{0}}$ s.t. $u_{\lambda_{n}} \rightarrow u_{0}$ in $W$ as $n \rightarrow +\infty$. Thus $\widetilde{u}_{\lambda_{n}} = t(\lambda_{0}, u_{\lambda_{n}})u_{\lambda_{n}} \rightarrow u_{0}$ in $W$ as $n \rightarrow +\infty$ and $\lim_{n \rightarrow \infty}Q(\widetilde{u}_{\lambda_{n}}) = Q(u_{0})$. Therefore,
		\begin{eqnarray} \label{eq3.12}
			\liminf_{\lambda \rightarrow \lambda_{0}^{+}}\frac{h(\lambda) - h(\lambda_{0})}{\lambda - \lambda_{0}} &\geq& -\limsup_{\lambda \rightarrow \lambda_{0}^{+}}Q(\widetilde{u}_{\lambda}) \nonumber \\
			&=& -\lim_{n \rightarrow \infty}Q(\widetilde{u}_{\lambda_{n}}) \nonumber \\
			&=& -Q(u_{0}) \nonumber \\
			&\geq& -\hat{d}(\lambda_{0}).
		\end{eqnarray}
		By (\ref{eq3.10}) and (\ref{eq3.12}), we have
		$$
		\lim_{\lambda \rightarrow \lambda_{0}^{+}}\frac{h(\lambda) - h(\lambda_{0})}{\lambda - \lambda_{0}} = -\hat{d}(\lambda_{0}).
		$$
		
		Similar arguments enable as to conclude that
		$$
		\lim_{\lambda \rightarrow \lambda_{0}^{-}}\frac{h(\lambda) - h(\lambda_{0})}{\lambda - \lambda_{0}} = -d(\lambda_{0}).
		$$
		\qed\vskip 5pt
		
		Next, we study the asymptotical behavior of $d(\lambda)$ and $\hat{d}(\lambda)$. Assume
		
		$(F_{5})$ There exists \(l > 0\) such that $\Phi_{\lambda}(u) \leq -l\lambda Q(u)$ for any $u \in \mathcal{K}_{\lambda}$.
		
		$(F_{6})$ There exists \(k > 0\) such that $\Phi_{\lambda}(u) \geq -k\lambda Q(u)$ for any $u \in \mathcal{K}_{\lambda}$.
		
		\begin{theorem}
			\label{thm2.2}  Suppose $(F_{1})$ - $(F_{4})$.
			
			$(i)$ Assume that $(F_{5})$ holds with $0 < l < 1$. If $h(\lambda)$ is well-defined in $(-\infty, b)$, then $\lim_{\lambda \rightarrow -\infty}d(\lambda) = +\infty$.

			$(ii)$ Assume $(F_{6})$ holds with $k > 1$. If $h(\lambda)$ is well-defined in $(-\infty, b)$, then $\lim_{\lambda \rightarrow -\infty}\hat{d}(\lambda) = 0$.
			
			$(iii)$ Assume $(F_{5})$ and $(F_{6})$ hold with $l > k > 0$ and $h(\lambda)$ is well-defined in $(a, 0)$. If $0 < k < l < 1$, then $\lim_{\lambda \rightarrow 0^{-}}\hat{d}(\lambda) = 0$; if $l > k > 1$, then $\lim_{\lambda \rightarrow 0^{-}}d(\lambda) = +\infty$.
			
			$(iv)$ Assume $(F_{5})$ and $(F_{6})$ hold with $k = l$ and $\lambda < 0$. Then $d(\lambda) \equiv \hat{d}(\lambda)$ and $d(\lambda) = \frac{C}{k}(-\lambda)^{\frac{1}{k} - 1}$ where $C > 0$ is a constant.
		\end{theorem}
		
		\textit{Proof.  } For $\lambda < 0$, we define
		\[
		c(\lambda) = \frac{h(\lambda)}{\lambda},
		\]
		\[
		c'_{-}(\lambda) = \lim_{\Delta \lambda \rightarrow 0^{-}}\frac{c(\lambda + \Delta \lambda) - c(\lambda)}{\Delta \lambda},
		\]
		\[
		c'_{+}(\lambda) = \lim_{\Delta \lambda \rightarrow 0^{+}}\frac{c(\lambda + \Delta \lambda) - c(\lambda)}{\Delta \lambda}.
		\]
		We may assume that $b > -1$ and $\lambda < 0$. By Theorem \ref{thm2.1} $(iii)$, $h'_{+}(\lambda) = -\hat{d}(\lambda)$, $h'_{-}(\lambda) = -d(\lambda)$.
		
		By $(F_{5})$, $h(\lambda) \leq -l\lambda d(\lambda)$. Since $h'_{\pm}(\lambda) = c'_{\pm}(\lambda)\lambda + c(\lambda)$,
		\begin{equation} \label{eq3.13}
			c'_{+}(\lambda)\lambda + c(\lambda) \leq c'_{-}(\lambda)\lambda + c(\lambda) \leq \frac{1}{l}c(\lambda),
		\end{equation}
		i.e.
		\begin{equation} \label{eq3.14}
			c'_{+}(\lambda)\lambda \leq c'_{-}(\lambda)\lambda \leq (\frac{1}{l}-1)c(\lambda).
		\end{equation}
		Set $g(\lambda) = |c(-1)|(-\lambda)^{\frac{1}{l}-1}$ for $\lambda < 0$. Then $g'(\lambda)\lambda = (\frac{1}{l}-1)g(\lambda)$ and $g(\lambda) > 0$. Now, set $\psi(\lambda) = \frac{c(\lambda)}{g(\lambda)}$, then
		\begin{equation} \label{eq3.15}
			\psi'_{\pm}(\lambda) = \frac{g(\lambda)c'_{\pm}(\lambda) - g'(\lambda)c(\lambda)}{g^{2}(\lambda)} \geq 0.
		\end{equation}
		Therefore, $\psi$ is monotonic increasing in $(-\infty, -1]$, then
		$$
		\psi(\lambda) \leq \psi(-1) = \frac{c(-1)}{|c(-1)|},
		$$
		i.e.
		$$
		c(\lambda) \leq c(-1)(-\lambda)^{\frac{1}{l}-1}, \forall \lambda < -1.
		$$
		When $h(\lambda) > 0$ and $\lambda < 0$, $c(\lambda) < 0$. When $l \in (0, 1)$,
		$$
		c(\lambda) \leq c(-1)(-\lambda)^{\frac{1}{l}-1} \rightarrow -\infty \ as \ \lambda \rightarrow -\infty.
		$$
		By $(F_{5})$, $d(\lambda) \geq -\frac{1}{l}c(\lambda)$. Thus $\lim_{\lambda \rightarrow -\infty}d(\lambda) = +\infty$ and we prove $(i)$.
		
		By $(F_{6})$, $h(\lambda) \geq -k\lambda \hat{d}(\lambda)$. A similar argument implies that:
		$$
		c(\lambda) \geq c(-1)(-\lambda)^{\frac{1}{k}-1}, \forall \lambda < -1.
		$$
		By $(F_4)$ and Theorem \ref{thm2.1} $(ii)$, for any $\lambda_{0} \in (a,b)$, $h(\lambda_{0}) = \lim_{\lambda \rightarrow \lambda_{0}}h(\lambda) > 0$, i.e. $h(\lambda) > 0$ on its definition domain. Hence, when $\lambda < 0$, $c(\lambda) < 0$. When $k > 1$,
		$$
		c(\lambda) \geq c(-1)(-\lambda)^{\frac{1}{k}-1} \rightarrow 0 \ as \ \lambda \rightarrow -\infty.
		$$
		By the definition of $\mathcal{K}_\lambda$, $\Phi_{\lambda}(u) = h(\lambda), \forall  u \in \mathcal{K}_\lambda$. Then by $(F_{6})$, $\Phi_{\lambda}(u) \geq -k\lambda Q(u), \forall u \in \mathcal{K}_\lambda$. Therefore, $-k\lambda \hat{d}(\lambda) = -k\lambda \sup_{u \in \mathcal{K}_{\lambda}}Q(u) \leq \sup_{u \in \mathcal{K}_{\lambda}}\Phi_{\lambda}(u) = \sup_{u \in \mathcal{K}_{\lambda}}h(\lambda) = h(\lambda)$, i.e. $\hat{d}(\lambda) \leq -\frac{1}{k}c(\lambda)$. Thus $\lim \limits_{\lambda \rightarrow -\infty}\hat{d}(\lambda) = 0$ and $(ii)$ is proved.
		
		Next we prove $(iii)$. We may assume that $a < -1$. By the proof of $(i)$ and $(ii)$, we know that:
		$$
		c(-1)(-\lambda)^{\frac{1}{l}-1} \leq c(\lambda) \leq c(-1)(-\lambda)^{\frac{1}{k}-1}, \forall \lambda \in (-1, 0).
		$$
		$h(\lambda) \geq -k\lambda Q(u) > 0$, thus $c(\lambda) < 0$. If $0 < k < l < 1$, then
		$$
		0 > c(\lambda) \geq c(-1)(-\lambda)^{\frac{1}{l}-1} \rightarrow 0 \ as \ \lambda \rightarrow 0^{-}.
		$$
		By $(F_{6})$, $\hat{d}(\lambda) \leq -\frac{1}{k}c(\lambda)$. Thus $\lim_{\lambda \rightarrow 0^{-}}\hat{d}(\lambda) = 0$. If $l > k > 1$, then
		$$
		c(\lambda) \leq c(-1)(-\lambda)^{\frac{1}{k}-1} \rightarrow -\infty \ as \ \lambda \rightarrow 0^{-}.
		$$
		By $(F_{5})$, $d(\lambda) \geq -\frac{1}{l}c(\lambda)$. Thus $\lim_{\lambda \rightarrow 0^{-}}d(\lambda) = +\infty$.
		
		Finally, if $(F_{5})$ and $(F_{6})$ hold with $k = l$, then $h(\lambda) = -k\lambda Q(u)$ for any $u \in \mathcal{K}_{\lambda}$. Thus $Q(u)$ is a constant for $u \in \mathcal{K}_{\lambda}$ and we conclude that $d \equiv \hat{d}$. By Theorem \ref{thm2.1} $(iii)$, $h(\lambda)$ is differentiable and $h'(\lambda) = -d(\lambda)$. Therefore, $h(\lambda) = k\lambda h'(\lambda)$. Since $h(\lambda) = -k\lambda Q(u) > 0$ when $\lambda < 0$, $h(\lambda) = C(-\lambda)^{\frac{1}{k}}$ for some $C > 0$, showing $d(\lambda) = \frac{C}{k}(-\lambda)^{\frac{1}{k}-1}$.  \qed\vskip 5pt

		To get stronger results when $\lambda \rightarrow 0^{-}$, we have the following assumptions and discussions.
		
		$(F_{7})$ $u$ is a solution of (\ref{eq1.1}), then $S(u) \leq d F(u)$ and $-\lambda Q(u) \geq e F(u)$ for some $d > 0$, $e > 0$.
		
		$(F_{8})$ For $u \in W$, $F(u) \leq C[S(u)^{\theta_{1}}Q(u)^{\tau_{1}} + S(u)^{\theta_{2}}Q(u)^{\tau_{2}}]$ where $C > 0$ is a constant and $\theta_{i} > 1$, $\tau_{i} > 0$, $i = 1$, $2$.
		
		\begin{theorem}
			\label{thm2.3}  Suppose $(F_{7})$ and $(F_{8})$. If $(\lambda_{n}, u_{n})$ are solutions of (\ref{eq1.1}) with $\lambda_{n} \rightarrow 0^{-}$, then $Q(u_{n}) \rightarrow +\infty$. In particular, $\lim_{\lambda \rightarrow 0^{-}}d(\lambda) = +\infty$.
		\end{theorem}
		
		\textit{proof.  } Arguing by contradiction, let us assume that
  $Q(u_{n})\nrightarrow +\infty$.
  We can then suppose that $\{Q(u_{n})\}_{n = 1}^{+\infty}$ is bounded passing to a subsequence if necessary.
		
		By $(F_{7})$,
		$$
		F(u_{n}) \geq \frac{1}{d} S(u_{n}) > 0, \ \ \  -\lambda_{n} Q(u_{n}) \geq e F(u_{n}) > 0.
		$$
		Thus $\lambda_{n} < 0$. Sending $n$ to $+\infty$, $F(u_{n}) \rightarrow 0$. Thus $S(u_{n}) \rightarrow 0$.
		
		On the other hand, by $(F_{8})$,
		$$
		F(u_{n}) \leq C[S(u_{n})^{\theta_{1}}Q(u_{n})^{\tau_{1}} + S(u_{n})^{\theta_{2}}Q(u_{n})^{\tau_{2}}].
		$$
		Hence,
		$$
		1 \leq dC[S(u_{n})^{\theta_{1}-1}Q(u_{n})^{\tau_{1}} + S(u_{n})^{\theta_{2}-1}Q(u_{n})^{\tau_{2}}].
		$$
		But $\theta_{i}-1 > 0$, $i = 1$, $2$, contradicting the fact $S(u_{n}) \rightarrow 0$ and $\{Q(u_{n})\}_{n = 1}^{+\infty}$ is bounded. \qed\vskip 5pt
		
		\begin{corollary} \label{cor2.4}
			$(i)$ Assume $(F_{1})$ - $(F_{5})$ hold with $l \in (0, 1)$. If $h(\lambda) > 0$ is well-defined and $d(\lambda) \equiv \hat{d}(\lambda)$ in $(-\infty, b)$, then $\exists \rho_{1} > 0$ s.t. for any $\rho \in (\rho_{1}, +\infty)$, (\ref{eq1.1}) admits a solution $(\lambda, u_{\lambda})$ with $Q(u_{\lambda}) = \rho$.
			
			$(ii)$ Assume $(F_{1})$ - $(F_{4})$ and $(F_{6})$ hold with $k > 1$. If $h(\lambda)$ is well-defined and $d(\lambda) \equiv \hat{d}(\lambda)$ in $(-\infty, b)$, then $\exists \rho_{2} > 0$ s.t. for any $\rho \in (0, \rho_{2})$, (\ref{eq1.1}) admits a solution $(\lambda, u_{\lambda})$ with $Q(u_{\lambda}) = \rho$.
			
			$(iii)$ Assume $(F_{7})$ and $(F_{8})$ hold. If $h(\lambda)$ is well-defined and can be achieved and $d(\lambda) \equiv \hat{d}(\lambda)$ in $(a, 0)$, then $\exists \rho_{3} > 0$ s.t. for any $\rho \in (\rho_{3}, +\infty)$, (\ref{eq1.1}) admits a solution $(\lambda, u_{\lambda})$ with $Q(u_{\lambda}) = \rho$.
			
			$(iv)$ Assume $(F_{1})$ - $(F_{6})$ hold with $0 < k < l < 1$ or $l > k > 1$. If $h(\lambda)$ is well-defined and $d(\lambda) \equiv \hat{d}(\lambda)$ in $(-\infty, 0)$, then for any $\rho \in (0, +\infty)$, (\ref{eq1.1}) admits a solution $(\lambda, u_{\lambda})$ with $Q(u_{\lambda}) = \rho$.
			
			$(v)$ Assume $(F_{1})$ - $(F_{6})$ hold with $k = l \neq 1$. If $h(\lambda)$ is well-defined in $(-\infty, 0)$, then for any $\rho \in (0, +\infty)$, (\ref{eq1.1}) admits a solution $(\lambda, u_{\lambda})$ with $Q(u_{\lambda}) = \rho$.
			
			$(vi)$ Assume $(F_{1})$ - $(F_{4})$, $(F_{6})$ - $(F_{8})$ hold with $k > 1$. If $h(\lambda)$ is well-defined and $d(\lambda) \equiv \hat{d}(\lambda)$ in $(-\infty, 0)$, then for any $\rho \in (0, +\infty)$, (\ref{eq1.1}) admits a solution $(\lambda, u_{\lambda})$ with $Q(u_{\lambda}) = \rho$.
		\end{corollary}
		
		\textit{Proof.  } We prove $(i)$. Standard arguments yield that $h(\lambda)$ can be achieved by $u_\lambda$ and thus $\mathcal{K}_\lambda$ is not empty ($u_\lambda \in \mathcal{K}_\lambda$). By Theorem \ref{thm2.1} $(i)$, $d(\lambda)$ is continuous in $(-\infty, b)$ since $d \equiv \hat{d}$. By Theorem \ref{thm2.2} $(i)$, $\lim_{\lambda \rightarrow -\infty}d(\lambda) = +\infty$. Thus for any $\rho \in (\rho_{1}, \infty)$ where $\rho_{1} = \inf_{\lambda < b}d(\lambda)$, there exists some $\lambda < b$ such that $d(\lambda) = \rho$. Since $d(\lambda) \equiv \hat{d}(\lambda)$, $Q(u) = d(\lambda), \forall u \in \mathcal{K}_\lambda$. Take $u_{\lambda} \in \mathcal{K}_{\lambda}$ and we complete the proof of $(i)$. Similar arguments apply to $(ii)$ - $(vi)$. \qed\vskip 5pt
		
		\begin{remark}	\label{rmk2.5}
		If there exists a global branch
		$$
		\lambda \mapsto u_\lambda \in C^1((-\infty,0),W \setminus \{0\}),
		$$
		such that $u_\lambda$ solves (\ref{eq1.1}), it is easy to verify that $\partial_\lambda \Phi_\lambda(u_\lambda) = -Q(u_\lambda)$. Then arguments in Theorems \ref{thm2.2}, \ref{thm2.3} and Corollary \ref{cor2.4} also hold if we replace $d(\lambda), \hat{d}(\lambda)$ with $Q(u_\lambda)$.
		\end{remark}

		\section{\textbf{Abstract framework 3}} \label{AF3}
		
		Throughout this section, we always assume that $h(\lambda)$ is well-defined and achieved by $u_{\lambda}$ in $(-\infty,0)$. Under some abstract settings, we will prove the uniqueness of the ground state solution when $\lambda \rightarrow -\infty$ or $\lambda \rightarrow 0^{-}$ and obtain the existence of normalized solutions. This section contains two subsections. We will need the following assumptions:
		
		$(S_{1})$ $D_{u}S(tu) = tD_{u}S(u)$, $D_{u}G(tu) = tD_{u}G(u)$, $D_{u}Q(tu) = tD_{u}Q(u)$.
		
		$(S_{2})$ There exists $T(\lambda): W \rightarrow W, W^{\ast} \rightarrow W^{\ast}$ such that
		$$
		D_{u}S(T(\lambda)u) = \frac{1}{\lambda}T(\lambda)D_{u}S(u),
		$$
		$$
		D_{u}Q(T(\lambda)u) = T(\lambda)D_{u}Q(u),
		$$
		$$
		D_{u}F(T(\lambda)u) = T(\lambda)D_{u}F(u).
		$$
		
		$(G_{1})$ $T(\lambda)$ is given by $(S_{2})$ and
		$$
		\lim_{\lambda \rightarrow -\infty}\frac{1}{\lambda}T(\lambda)D_{u}G(u) = 0.
		$$
		
		$(Q_{1})$ $T(\lambda)$ is given by $(S_{2})$ and
		$$
		Q(T(\lambda)u) = \lambda^{\gamma}Q(u).
		$$
		
		$(N_{1})$ There exists $p > 2$ such that
		$$
		\lim_{|t| \rightarrow 0}\frac{D_{u}F(tu)}{(tu)^{p-1}} = m_{1} \in (0,+\infty), \lim_{|t| \rightarrow 0}\frac{pF(tu)}{(tu)^{p}} = m_{1} \in (0,+\infty).
		$$
		
		$(N_{2})$ There exists $q > 2$ such that
		$$
		\lim_{|t| \rightarrow +\infty}\frac{D_{u}F(tu)}{(tu)^{q-1}} = m_{2} \in (0,+\infty), \lim_{|t| \rightarrow +\infty}\frac{qF(tu)}{(tu)^{q}} = m_{2} \in (0,+\infty).
		$$
		
		\subsection{The case when $\lambda \rightarrow 0^{-}$}
		
		We always assume that $D_{u}G = 0$ in this subsection. Let $u_{\lambda}$ be the ground state of (1.1) as $\lambda \rightarrow 0^{-}$. We consider
		\begin{equation} \label{eqa.4}
			w_{\lambda} = |\lambda|^{-\frac{1}{p-2}}T(\lambda)u_{\lambda},
		\end{equation}
		which satisfies the following equation
		\begin{equation} \label{eqa.5}
			D_{w}S(w) - |\lambda|^{-\frac{p-1}{p-2}}D_{w}F(|\lambda|^{\frac{1}{p-2}}w) + D_{w}Q(w) = 0 \ in \ W^{\ast}.
		\end{equation}
		Let
		$$
		\widehat{\Phi}_{\lambda}(w) = S(w) - |\lambda|^{-\frac{p}{p-2}}F(|\lambda|^{\frac{1}{p-2}}w) + Q(w),
		$$
		$$
		\widehat{\Phi}_{0}(w) = S(w) - m_{1}F_{0,p}(w) + Q(w), \ where \ D_{w}F_{0,p}(w) = w^{p-1}.
		$$
		Set
		$$
		\widehat{\mathcal{N}}_{\lambda} = |\lambda|^{-\frac{1}{p-2}}T(\lambda)\mathcal{N}_{\lambda}.
		$$
		We assume that $T(\lambda)$ is invertible. Then $\widehat{\mathcal{N}}_{\lambda}$ is a natural constraint of $\widehat{\Phi}_{\lambda}$. Let $\widehat{\mathcal{N}}_{0}$ be a natural constraint of $\widehat{\Phi}_{0}$ and set
		$$
		\widehat{h}(\lambda) = \inf_{w \in \widehat{\mathcal{N}}_{\lambda}}\widehat{\Phi}_{\lambda}(w),
		$$
		$$
		\mathcal{\widehat{K}}_{\lambda} = \{w \in \widehat{\mathcal{N}}_{\lambda}: D_{w}\widehat{\Phi}_{\lambda}(w) = 0, \widehat{\Phi}_{\lambda}(w) = \widehat{h}(\lambda)\}.
		$$
		We have the following assumptions:
		
		$(A_{1})$ There exists $\delta > 0$ such that for any $(\lambda, w) \in (-\delta,0] \times W \backslash \{0\}$, there exists a unique function $\widehat{t}: (-\delta,0] \times W \backslash \{0\} \rightarrow (0, \infty)$ such that
		\begin{equation}
			when \ t > 0, \ then \ (\lambda, t w) \in \mathcal{\widehat{N}} \Leftrightarrow t = \widehat{t}(\lambda, w),
			\nonumber
		\end{equation}
		$\widehat{\Phi}_{\lambda}(\widehat{t}(\lambda, w)w) = \max_{t > 0}\widehat{\Phi}_{\lambda}(tw)$, $\widehat{t} \in C((-\delta,0] \times W \backslash \{0\}, (0, \infty))$, where $\mathcal{\widehat{N}} := \cup_{-\delta < \lambda \leq 0}\mathcal{\widehat{N}}_{\lambda}$.
		
		$(A_{2})$ $w_{\lambda}$ is bounded in $W$ whenever $w_{\lambda} \in \mathcal{\widehat{K}}_{\lambda}$ and $\lambda \rightarrow 0^{-}$.
		
		$(A_{3})$ $D_{w}\widehat{\Phi}_{0}(w) = 0$ admits an unique and non-degenerate ground state $w^{\ast}$.
		
		$(A_{4})$ $\widehat{\Phi}_{0}$ satisfies (PS)$_{c}$ condition in $W$ for $c \leq \widehat{h}(0)$.
		
		$(A_{5})$ $\liminf_{\lambda \rightarrow 0^{-}}\widehat{h}(\lambda) > 0$.
		
		\begin{theorem}
			\label{thmc.1} Assume that $(S_{1})$, $(S_{2})$, $(N_{1})$, $(A_{1})$ - $(A_{5})$ hold true. Then there exists $\Lambda_{1} < 0$ such that when $\lambda \in (\Lambda_{1},0)$, the ground state of (1.1) is unique and non-degenerate.
		\end{theorem}
		
		\textit{Proof.  } Step 1: $\limsup_{\lambda \rightarrow 0^{-}} \widehat{h}(\lambda) \leq \widehat{h}(0)$.
		
		Consider $\widehat{t}(\lambda, w^{\ast})w^{\ast} \in \mathcal{\widehat{N}}_{\lambda}$. By $(A_{1})$, $\widehat{t}(\lambda, w^{\ast})w^{\ast} \rightarrow \widehat{t}(0, w^{\ast})w^{\ast} = w^{\ast}$ in $W$ as $\lambda \rightarrow 0^{-}$. Together with $(N_{1})$, we know that
		\begin{eqnarray}
			&& \widehat{h}(\lambda) \nonumber \\
&\leq& \widehat{\Phi}_{\lambda}(\widehat{t}(\lambda, w^{\ast})w^{\ast}) \nonumber \\
			&=& \widehat{\Phi}_{0}(\widehat{t}(\lambda, w^{\ast})w^{\ast}) + m_{1}F_{0,p}(\widehat{t}(\lambda, w^{\ast})w^{\ast}) - |\lambda|^{-\frac{p}{p-2}}F(|\lambda|^{\frac{1}{p-2}}\widehat{t}(\lambda, w^{\ast})w^{\ast}) \nonumber \\
			&=& \widehat{\Phi}_{0}(w^{\ast}) + o(1) \nonumber \\
			&=& \widehat{h}(0) + o(1). \nonumber
		\end{eqnarray}
		Letting  $\lambda \rightarrow 0^{-},$ we obtain that $\limsup_{\lambda \rightarrow 0^{-}}\widehat{h}(\lambda) \leq \widehat{h}(0)$.
		
		Step 2: Let $u_{\lambda_{n}} = u_{n}$ be the ground state of (1.1) as $\lambda_{n} \rightarrow 0^{-}$ and $w_{n} = w_{\lambda_{n}}$ is given by (\ref{eqa.4}). Then $w_{n} \rightarrow w^{\ast}$ in $W$ up to a subsequence.
		
		Since $\Phi_{\lambda_{n}}(u_{n}) = h(\lambda_{n})$, we know that $\widehat{\Phi}_{\lambda_{n}}(w_{n}) = \widehat{h}(\lambda_{n})$. By $(A_{2})$, $w_{n}$ is bounded in $W$. Then by $(N_{1})$,
		$$
		|\lambda_{n}|^{-\frac{p}{p-2}}F(|\lambda_{n}|^{\frac{1}{p-2}}w_{n}) - m_{1}F_{0,p}(w_{n}) \rightarrow 0.
		$$
		Noticing that
		$$
		\widehat{\Phi}_{0}(w_{n}) = \widehat{\Phi}_{\lambda_{n}}(w_{n}) + |\lambda_{n}|^{-\frac{p}{p-2}}F(|\lambda_{n}|^{\frac{1}{p-2}}w_{n}) - m_{1}F_{0,p}(w_{n}) = \widehat{h}(\lambda_{n}) + o(1),
		$$
		we conclude that $\limsup_{n \rightarrow \infty}\widehat{\Phi}_{0}(w_{n}) \leq \widehat{h}(0)$. This, together with $(A_{5})$, we may assume that $\widehat{\Phi}_{0}(w_{n}) \rightarrow c \in (0,\widehat{h}(0)],$ passing to a subsequence if necessary. Then by $(N_{1})$, as $n \rightarrow \infty$,
		$$
		D_{w}\widehat{\Phi}_{0}(w_{n}) = D_{w}\widehat{\Phi}_{\lambda_{n}}(w_{n}) + |\lambda_{n}|^{-\frac{p-1}{p-2}}D_{w}F(|\lambda_{n}|^{\frac{1}{p-2}}w_{n}) - m_{1}D_{w}F_{0,p}(w_{n}) \rightarrow 0.
		$$
		By $(A_{4})$, up to a subsequence, $w_{n} \rightarrow w_{0}$ in $W$. Since $\widehat{\Phi}_{0} \in C^{1}$, $D_{w}\widehat{\Phi}_{0}(w_{0}) = 0$. Noticing that
		$$
		\widehat{\Phi}_{0}(w_{0}) = \lim_{n \rightarrow +\infty}\widehat{\Phi}_{0}(w_{n}) > 0,
		$$
		we have that $w_{0} \neq 0$. Thus $w_{0} \in \widehat{\mathcal{N}}_{0}$.
		
		Next, we will show that $w_{0} \in \widehat{\mathcal{K}}_{0}$. $\forall \epsilon > 0$, $\exists X > 0$, $\forall n > X$, we have that:
		\begin{eqnarray} \label{eqa.6}
			\widehat{h}(0) &\leq& \widehat{\Phi}_{0}(w_{0}) \nonumber \\
			&\leq& \widehat{\Phi}_{0}(w_{n}) + \frac{\epsilon}{2} \nonumber \\
			&=& \widehat{\Phi}_{\lambda_{n}}(w_{n}) + |\lambda_{n}|^{-\frac{p}{p-2}}F(|\lambda_{n}|^{\frac{1}{p-2}}w_{n}) - m_{1}F_{0,p}(w_{n}) + \frac{\epsilon}{2} \nonumber \\
			&\leq& \widehat{\Phi}_{\lambda_{n}}(w_{n}) + \epsilon.
		\end{eqnarray}
		Letting $n \rightarrow \infty,$ we get $\widehat{h}(0) \leq \limsup_{\lambda \rightarrow 0^{-}}\widehat{h}(\lambda)$. Therefore, $$\limsup_{\lambda \rightarrow 0^{-}}\widehat{h}(\lambda) = \widehat{h}(0).$$ We can derive that $\widehat{\Phi}_{0}(w_{0}) = \widehat{h}(0)$ by (\ref{eqa.6}) immediately. Thus $w_{0} \in \widehat{\mathcal{K}}_{0}$. Finally, using $(A_{3})$, $w_{0} = w^{\ast}$.
		
		Step 3: There exists $\Lambda_{1} < 0$ such that when $\lambda \in (\Lambda_{1},0)$, the positive ground state of (1.1) is unique and non-degenerate.
		
		Since $w^{\ast}$ is non-degenerate, we can use the implicit function theorem. Consider
		$$
		P(w,\lambda) = D_{w}\widehat{\Phi}_{\lambda}(w) = 0.
		$$
		Then $P(w^{\ast},0) = 0$, $P_{w}(w^{\ast},0) = D_{ww}\widehat{\Phi}_{\lambda}(w^{\ast})$ is invertible. By applying the implicit function theorem to the map $P$ at $(w^{\ast},0)$, we derive that there exists $\delta > 0$ small enough such that $w_{\lambda}$ is the unique and non-degenerate solution of $P(w,\lambda) = 0$ in the neighborhood $\{(w,\lambda) \in W \times (-\infty,0]:\|w - w^{\ast}\| + |\lambda| < \delta\}$. In particular, $w_{0} = w^{\ast}$. The by Step 1, for $\lambda < 0$ with $|\lambda|$ small enough, $w_{\lambda}$ is the unique and non-degenerate ground state of $D_{w}\widehat{\Phi}_{\lambda}(w) = 0$. Note that the uniqueness and non-degeneracy of $w_{\lambda}$ is equivalent to the uniqueness and non-degeneracy of $u_{\lambda}$. Therefore, there exists $\Lambda_{1} < 0$ with $|\Lambda_{1}|$ small such that when $\lambda \in (\Lambda_{1},0)$, the positive ground state of (1.1) is unique and non-degenerate. \qed\vskip 5pt
		
		\begin{corollary}
			\label{corc.2} Assume that $(S_{1})$, $(S_{2})$, $(Q_{1})$, $(N_{1})$, $(A_{1})$ - $(A_{5})$ hold true. Let $u_{\lambda}$ be the unique ground state of (1.1). Then
			
			$(i)$ $\lim_{\lambda \rightarrow 0^{-}}Q(u_{\lambda}) = +\infty$ if $\gamma - 2/(p-2) > 0$;
			
			$(ii)$ $\lim_{\lambda \rightarrow 0^{-}}Q(u_{\lambda}) = Q(w^{\ast})$ if $\gamma - 2/(p-2) = 0$;
			
			$(iii)$ $\lim_{\lambda \rightarrow 0^{-}}Q(u_{\lambda}) = 0$ if $\gamma - 2/(p-2) < 0$.
		\end{corollary}
		
		\textit{Proof.  } Using $(S_{1})$, $(Q_{1}),$ and the proof of Theorem \ref{thmc.1}, we have
		$$
		Q(w_{\lambda}) = |\lambda|^{\gamma - \frac{2}{p-2}}Q(u_{\lambda}) \rightarrow Q(w^{\ast}) \ as \ \lambda \rightarrow 0^{-}.
		$$
		Then the claims follow. \qed\vskip 5pt
		
		\begin{corollary}
		\label{corc.3}	Assume that $(S_{1})$, $(S_{2})$, $(Q_{1})$, $(N_{1})$, $(A_{1})$ - $(A_{5})$ hold true. Then
			
			$(i)$ if $\gamma - 2/(p-2) > 0$, there exists $c_{1} > 0$ large such that for any $c > c_{1}$, (1.1) admits a ground state $u_{\lambda}$ with $\lambda < 0$ and $|\lambda|$ small such that $Q(u_{\lambda}) = c$;
			
			$(ii)$ if $\gamma - 2/(p-2) < 0$, there exists $c_{2} > 0$ small such that for any $c < c_{2}$, (1.1) admits a ground state $u_{\lambda}$ with $\lambda < 0$ and $|\lambda|$ small such that $Q(u_{\lambda}) = c$.
		\end{corollary}
		
		\subsection{The case when $\lambda \rightarrow -\infty$} \label{AF3.2}
		
		First, we assume that $D_{u}G = 0$. Let $u_{\lambda}$ be a ground state of (1.1) as $\lambda \rightarrow -\infty$ and $\mu = 1/\lambda$. We consider
		\begin{equation} \label{eqa.2}
			v_{\mu} = |\mu|^{\frac{1}{q-2}}T(\lambda)u_{\lambda},
		\end{equation}
		which satisfies the following equation
		\begin{equation} \label{eqa.1}
			D_{v}S(v) - |\mu|^{\frac{q-1}{q-2}}D_{v}F(|\mu|^{-\frac{1}{q-2}}v) + D_{v}Q(v) = 0 \ in \ W^{\ast}.
		\end{equation}
		Let
		$$
		\widetilde{\Phi}_{\mu}(v) = S(v) - |\mu|^{\frac{q}{q-2}}F(|\mu|^{-\frac{1}{q-2}}v) + Q(v),
		$$
		$$
		\widetilde{\Phi}_{0}(v) = S(v) - m_{2}F_{0,q}(v) + Q(v), \ where \ D_{v}F_{0,q}(v) = v^{q-1}.
		$$
		Set
		$$
		\widetilde{\mathcal{N}}_{\mu} = |\mu|^{\frac{1}{q-2}}T(1/\mu)\mathcal{N}_{1/\mu}.
		$$
		We assume that $T_{\lambda}$ is invertible. Then $\widetilde{\mathcal{N}}_{\mu}$ is a natural constraint of $\widetilde{\Phi}_{\mu}$. Let $\widetilde{\mathcal{N}}_{0}$ be a natural constraint of $\widetilde{\Phi}_{0}$ and set
		$$
		\widetilde{h}(\mu) = \inf_{v \in \widetilde{\mathcal{N}}_{\mu}}\widetilde{\Phi}_{\mu}(v),
		$$
		$$
		\mathcal{\widetilde{K}}_{\mu} = \{v \in \widetilde{\mathcal{N}}_{\mu}: D_{v}\widetilde{\Phi}_{\mu}(v) = 0, \widetilde{\Phi}_{\mu}(v) = \widetilde{h}(\mu)\}.
		$$
		Let us assume the following:
		
		$(B_{1})$ There exists $\delta > 0$ such that for any $(\mu, v) \in (-\delta,0] \times W \backslash \{0\}$, there exists a unique function $\widetilde{t}: (-\delta,0] \times W \backslash \{0\} \rightarrow (0, \infty)$ such that
		\begin{equation}
			when \ t > 0, \ then \ (\mu, t v) \in \mathcal{\widetilde{N}} \Leftrightarrow t = \widetilde{t}(\mu, v),
			\nonumber
		\end{equation}
		$\widetilde{\Phi}_{\mu}(\widetilde{t}(\mu, v)v) = \max_{t > 0}\widetilde{\Phi}_{\mu}(tv)$, $\widetilde{t} \in C((-\delta,0] \times W \backslash \{0\}, (0, \infty))$, where $\mathcal{\widetilde{N}} := \cup_{-\delta < \mu \leq 0}\mathcal{\widetilde{N}}_{\mu}$.
		
		$(B_{2})$ $v_{\mu}$ is bounded in $W$ whenever $v_{\mu} \in \mathcal{\widetilde{K}}_{\mu}$ and $\mu \rightarrow 0^{-}$.
		
		$(B_{3})$ $D_{v}\widetilde{\Phi}_{0}(v) = 0$ admits an unique and non-degenerate ground state $v^{\ast}$.
		
		$(B_{4})$ $\widetilde{\Phi}_{0}$ satisfies (PS)$_{c}$ condition in $W$ for $c \leq \widetilde{h}(0)$.
		
		$(B_{5})$ $\liminf_{\mu \rightarrow 0^{-}}\widetilde{h}(\mu) > 0$.
		
		\begin{theorem}
		\label{thmc.4}	Assume that $(S_{1})$, $(S_{2})$, $(N_{2})$, $(B_{1})$ - $(B_{5})$ hold true. Then there exists $\Lambda_{2} < 0$ such that when $\lambda < \Lambda_{2}$, the ground state of (1.1) is unique and non-degenerate.
		\end{theorem}
		
		\textit{Proof.  } Step 1: $\limsup_{\mu \rightarrow 0^{-}} \widetilde{h}(\mu) \leq \widetilde{h}(0)$.
		
		Consider $\widetilde{t}(\mu, v^{\ast})v^{\ast} \in \mathcal{\widetilde{N}}_{\mu}$. By $(B_{1})$, $\widetilde{t}(\mu, v^{\ast})v^{\ast} \rightarrow \widetilde{t}(0, v^{\ast})v^{\ast} = v^{\ast}$ in $W$ as $\mu \rightarrow 0^{-}$. The latter, together with $(N_{2})$, imply that:
		\begin{eqnarray}
			\widetilde{h}(\mu) &\leq& \widetilde{\Phi}_{\mu}(\widetilde{t}(\mu, v^{\ast})v^{\ast}) \nonumber \\
			&=& \widetilde{\Phi}_{0}(\widetilde{t}(\mu, v^{\ast})v^{\ast}) + m_{2}F_{0,q}(\widetilde{t}(\mu, v^{\ast})v^{\ast}) - |\mu|^{\frac{q}{q-2}}F(|\mu|^{-\frac{1}{q-2}}\widetilde{t}(\mu, v^{\ast})v^{\ast}) \nonumber \\
			&=& \widetilde{\Phi}_{0}(v^{\ast}) + o(1) \nonumber \\
			&=& \widetilde{h}(0) + o(1).
		\end{eqnarray}
		Letting $\mu \rightarrow 0^{-},$ we have $\limsup_{\mu \rightarrow 0^{-}}\widetilde{h}(\mu) \leq \widetilde{h}(0)$.
		
		Step 2: Let $u_{\lambda_{n}} = u_{n}$ be a ground state of (1.1) as $\lambda_{n} \rightarrow -\infty$ and $v_{n} = v_{\mu_{n}}$ is given by (\ref{eqa.2}). Then $v_{n} \rightarrow v^{\ast}$ in $W$ up to a subsequence.
		
		Since $\Phi_{\lambda_{n}}(u_{n}) = h(\lambda_{n})$, we know that $\widetilde{\Phi}_{\mu_{n}}(v_{n}) = \widetilde{h}(\mu_{n})$. By $(B_{2})$, $v_{n}$ is bounded in $W$. Then by $(N_{2})$,
		$$
		|\mu_{n}|^{\frac{q}{q-2}}F(|\mu_{n}|^{-\frac{1}{q-2}}v_{n}) - m_{2}F_{0,q}(v_{n}) \rightarrow 0.
		$$
		Noticing that
		$$
		\widetilde{\Phi}_{0}(v_{n}) = \widetilde{\Phi}_{\mu_{n}}(v_{n}) + |\mu_{n}|^{\frac{q}{q-2}}F(|\mu_{n}|^{-\frac{1}{q-2}}v_{n}) - m_{2}F_{0,q}(v_{n}) = \widetilde{h}(\mu_{n}) + o(1),
		$$
		we derive that $\limsup_{n \rightarrow \infty}\widetilde{\Phi}_{0}(v_{n}) \leq \widetilde{h}(0)$. This, together with $(B_{5})$, enable as to assume that $\widetilde{\Phi}_{0}(v_{n}) \rightarrow c \in (0,\widetilde{h}(0)]$ passing to a subsequence if necessary. Then by $(N_{2})$, as $n \rightarrow \infty$,
		$$
		D_{v}\widetilde{\Phi}_{0}(v_{n}) = D_{v}\widetilde{\Phi}_{\mu_{n}}(v_{n}) + |\mu_{n}|^{\frac{q-1}{q-2}}D_{v}F(|\mu_{n}|^{-\frac{1}{q-2}}v_{n}) - m_{2}D_{v}F_{0,q}(v_{n}) \rightarrow 0.
		$$
		By $(B_{4})$, up to a subsequence, $v_{n} \rightarrow v_{0}$ in $W$. Since $\widetilde{\Phi}_{0} \in C^{1}$, $D_{v}\widetilde{\Phi}_{0}(v_{0}) = 0$. Noticing that:
		$$
		\widetilde{\Phi}_{0}(v_{0}) = \lim_{n \rightarrow +\infty}\widetilde{\Phi}_{0}(v_{n}) > 0,
		$$
		we derive that $v_{0} \neq 0$. Thus $v_{0} \in \widetilde{\mathcal{N}}_{0}$.
		
		Next, we will show that $v_{0} \in \widetilde{\mathcal{K}}_{0}$. $\forall \epsilon > 0$, $\exists X > 0$, $\forall n > X$, there holds that
		\begin{eqnarray} \label{eqa.3}
			\widetilde{h}(0) &\leq& \widetilde{\Phi}_{0}(v_{0}) \nonumber \\
			&\leq& \widetilde{\Phi}_{0}(v_{n}) + \frac{\epsilon}{2} \nonumber \\
			&=& \widetilde{\Phi}_{\mu_{n}}(v_{n}) + |\mu_{n}|^{\frac{q}{q-2}}F(|\mu_{n}|^{-\frac{1}{q-2}}v_{n}) - m_{2}F_{0,q}(v_{n}) + \frac{\epsilon}{2} \nonumber \\
			&\leq& \widetilde{\Phi}_{\mu_{n}}(v_{n}) + \epsilon.
		\end{eqnarray}
		Letting $n \rightarrow \infty,$ we get that $\widetilde{h}(0) \leq \limsup_{\mu \rightarrow 0^{-}}\widetilde{h}(\mu)$. Therefore, $$\limsup_{\mu \rightarrow 0^{-}}\widetilde{h}(\mu) = \widetilde{h}(0).$$ We can derive that $\widetilde{\Phi}_{0}(v_{0}) = \widetilde{h}(0)$ by (\ref{eqa.3}) immediately. Thus $v_{0} \in \widetilde{\mathcal{K}}_{0}$. By $(B_{3})$, $v_{0} = v^{\ast}$.
		
		Step 3: There exists $\Lambda_{2} < 0$ such that when $\lambda < \Lambda_{2}$, the positive ground state of (1.1) is unique and non-degenerate.
		
		Since $v^{\ast}$ is non-degenerate, we can use the implicit function theorem. Consider
		$$
		K(v,\mu) = D_{v}\widetilde{\Phi}_{\mu}(v) = 0.
		$$
		Then $K(v^{\ast},0) = 0$, $K_{v}(v^{\ast},0) = D_{vv}\widetilde{\Phi}_{0}(v^{\ast})$ is invertible. By applying the implicit function theorem to the map $K$ at $(v^{\ast},0)$, we derive that there exists $\delta > 0$ small enough such that $v_{\mu}$ is the unique and non-degenerate solution of $K(v,\mu) = 0$ in the neighborhood $\{(v,\mu) \in W \times (-\infty,0]:\|v - v^{\ast}\| + |\mu| < \delta\}$. In particular, $v_{0} = v^{\ast}$. Therefore, by Step 1, for $\mu < 0$ with $|\mu|$ small enough, $v_{\mu}$ is the unique and non-degenerate ground state of $D_{v}\widetilde{\Phi}_{\mu}(v) = 0$. Note that the uniqueness and non-degeneracy of $v_{\mu}$ is equivalent to the uniqueness and non-degeneracy of $u_{\lambda}$. Therefore, there exists $\Lambda_{2} < 0$ with $|\Lambda_{2}|$ large such that when $\lambda < \Lambda_{2}$, the positive ground state of (1.1) is unique and non-degenerate.
		\qed\vskip 5pt
		
		\begin{corollary}
			Assume that $(S_{1})$, $(S_{2})$, $(Q_{1})$, $(N_{2})$, $(B_{1})$ - $(B_{5})$ hold true. Let $u_{\lambda}$ be the unique ground state of (1.1). Then
			
			$(i)$ $\lim_{\lambda \rightarrow -\infty}Q(u_{\lambda}) = 0$ if $\gamma - 2/(q-2) > 0$;
			
			$(ii)$ $\lim_{\lambda \rightarrow -\infty}Q(u_{\lambda}) = Q(v^{\ast})$ if $\gamma - 2/(q-2) = 0$;
			
			$(iii)$ $\lim_{\lambda \rightarrow -\infty}Q(u_{\lambda}) = +\infty$ if $\gamma - 2/(q-2) < 0$.
		\end{corollary}
		
		\textit{Proof.  } The proof is similar to the one of Corollary \ref{corc.2}. \qed\vskip 5pt
		
		\begin{corollary}
		\label{corc.6}	Assume that $(S_{1})$, $(S_{2})$, $(Q_{1})$, $(N_{2})$, $(B_{1})$ - $(B_{5})$ hold true. Then
			
			$(i)$ if $\gamma - 2/(q-2) > 0$, there exists $c_{3} > 0$ small such that for any $c < c_{3}$, (1.1) admits a ground state $u_{\lambda}$ with $\lambda < 0$ and $|\lambda|$ large such that $Q(u_{\lambda}) = c$;
			
			$(ii)$ if $\gamma - 2/(q-2) < 0$, there exists $c_{4} > 0$ large such that for any $c > c_{4}$, (1.1) admits a ground state $u_{\lambda}$ with $\lambda < 0$ and $|\lambda|$ large such that $Q(u_{\lambda}) = c$.
		\end{corollary}
		
		When $D_{u}G \neq 0$, we also assume $(G_{1})$, quite similarly, we can prove the following results.
		
		\begin{theorem}\label{athmc.7}
			Assume that $(S_{1})$, $(S_{2})$, $(G_{1})$, $(N_{2})$, $(B_{1})$ - $(B_{5})$ hold true. Then there exist $\Lambda_{3} < 0$ such that when $\lambda < \Lambda_{3}$, the positive ground state of (1.1) is unique and non-degenerate.
		\end{theorem}
		
		\begin{corollary}
			Assume that $(S_{1})$, $(S_{2})$, $(G_{1})$, $(Q_{1})$, $(N_{2})$, $(B_{1})$ - $(B_{5})$ hold true. Let $u_{\lambda}$ be the unique ground state of (1.1). Then
			
			$(i)$ $\lim_{\lambda \rightarrow -\infty}Q(u_{\lambda}) = 0$ if $\gamma - 2/(q-2) > 0$;
			
			$(ii)$ $\lim_{\lambda \rightarrow -\infty}Q(u_{\lambda}) = Q(v^{\ast})$ if $\gamma - 2/(q-2) = 0$;
			
			$(iii)$ $\lim_{\lambda \rightarrow -\infty}Q(u_{\lambda}) = +\infty$ if $\gamma - 2/(q-2) < 0$.
		\end{corollary}
		
		\begin{corollary}
			Assume that $(S_{1})$, $(S_{2})$, $(G_{1})$, $(Q_{1})$, $(N_{2})$, $(B_{1})$ - $(B_{5})$ hold true. Then
			
			$(i)$ if $\gamma - 2/(q-2) > 0$, there exists $c_{5} > 0$ small such that for any $c < c_{5}$, (1.1) admits a ground state $u_{\lambda}$ with $\lambda < 0$ and $|\lambda|$ large such that $Q(u_{\lambda}) = c$;
			
			$(ii)$ if $\gamma - 2/(q-2) < 0$, there exists $c_{6} > 0$ large such that for any $c > c_{6}$, (1.1) admits a ground state $u_{\lambda}$ with $\lambda < 0$ and $|\lambda|$ large such that $Q(u_{\lambda}) = c$.
		\end{corollary}

		\section{\textbf{Scalar Field Equations Involving the Fractional Laplacian}} \label{frac}
		
		In this section and the next ones; we provide some applications of our abstract framework. In this section, we study scalar field equations involving the fractional Laplacian. After showing the existence of normalized solutions, we also prove the orbital stability of standing waves for almost every $L^{2}$ mass in the existence range.
		
		\subsection{Scalar field equations involving a fractional Laplacian in the entire space}
		
		Consider the scalar field equations involving a fractional Laplacian of the form
		\begin{equation} \label{eq4.1}
			\left\{
			\begin{array}{cc}
				(-\Delta)^{s} u = \lambda u + f(|x|,u) \ in \ \mathbb{R}^{N}, 0 < s \leq 1, N \geq 2, \lambda < 0,\\
				u \geq 0 \ in \ \mathbb{R}^{N}, u(x) \rightarrow 0 \ as \ |x| \rightarrow +\infty.
			\end{array}
			\right.
		\end{equation}
		When $s = 1$, $(-\Delta)^{1}$ is the usual Laplacian; when $0 < s < 1$, the fractional Laplacian $(-\Delta)^{s}$ is defined via its multiplier $|\xi|^{2s}$ in Fourier space. If $u$ is smooth enough, it can also be computed by the following singular integral:
		\[
		(-\Delta)^{s} u(x) = c_{N,s}P.V. \int_{\mathbb{R}^{N}}\frac{u(x)-u(y)}{|x-y|^{N+2s}}dy,
		\]
		where $P.V.$ is the principal value, and $c_{N,s}$ is a normalization constant.
		
		As an application of our abstract framework 2, we will show the existence of normalized solutions. When $s = 1$ and $f(|x|,u) = |u|^{q - 2}u$. Recall that the nonlinearity is called $L^{2}$-supercritical case or the mass supercritical case when $q > 2 + N/4$, $L^{2}$-critical case when $q = 2 + N/4$, the $L^{2}$-subcritical case when $q < 2 + N/4$. When $s = 1$, C.A. Stuart proved the existence of normalized solutions for (\ref{eq4.1}) in a series of works \cite{Stu1, Stu2, Stu3} in the $L^{2}$-subcritical case. In the non-autonomous case, H. Hajaiej and C. A. Stuart provided a detailed study including a wide range of nonlinearities. They obtained the existence of the normalized solutions and they proved the orbital stability of standing waves \cite{Hajaiej-Stuart}.
		In the mass supercritical case, there are many results showing the existence of normalized solutions when $f$ is autonomous, see, e.g., \cite{BN, BV, HMOW, Jean}. As for non-autonomous nonlinearities, the second author has obtained some results in a recent work \cite{Song}. Here we extend these results to the case of $0 < s < 1$.
		
		We may assume that $f(r,t) = 0$ whenever $t < 0$ since we aim to show the existence of positive solutions. Furthermore, we suppose that $f(r,t) \in C^{1}([0, +\infty) \times [0, +\infty), \mathbb{R})$:
		
		$(f_{1})$ $f(r,0) = 0$ and there exist $2 < \alpha \leq \beta < 2_{s}^{\ast}$ such that
		\begin{equation}
			0 < (\alpha - 1) f(r,t) \leq f_{t}(r,t)t \leq (\beta - 1) f(r,t), \forall t > 0, r > 0,
			\label{eq4.2}
		\end{equation}
		where $2_{s}^{\ast} = 2N/(N - 2s)$ if $N > 2s$ and $2_{s}^{\ast} = +\infty$ if $N \leq 2s$.
		
		$(f_{2})$ $f(r,1)$ is bounded, i.e. there exists $M > 0$ such that $\vert f(r,1) \vert \leq M$, $\forall r > 0$.
		
		$(f_{3})$ $F_{r}(r,t)r \leq C (|t|^{q_{1}} + |t|^{q_{2}})$ for some $C > 0$, $q_{1}$, $q_{2} \in (2,2_{s}^{\ast})$ and there exists $\theta \leq 0$ such that
		\begin{equation}
			F_{r}(r,t)r \geq \theta F(r,t), \forall t > 0, r > 0,
			\nonumber
		\end{equation}
		where $F(r,t) = \int_{0}^{t}f(r,s)ds$.
		
		$(f_{4})$ $F_{r}(r,t)r \geq -C (|t|^{q_{1}} + |t|^{q_{2}})$ for some $C > 0$, $q_{1}$, $q_{2} \in (2,2_{s}^{\ast})$ and there exists $\tau \geq 0$ such that
		\begin{equation}
			F_{r}(r,t)r \leq \tau F(r,t), \forall t > 0, r > 0.
			\nonumber
		\end{equation}
		
		In this case, $W = H_{rad}^{s}(\mathbb{R}^{N})$, the radially symmetric fractional order Sobolev space,\\
		equipped with the norm
		$$
		\|u\| = \sqrt{\int_{\mathbb{R}^{N}}(|\xi|^{2s}\hat{u}^{2} + \hat{u}^{2})d\xi}, \hat{u} = \mathcal{F}(u),
		$$
		$$
		\Phi_{\lambda}(u) = \frac{1}{2}\int_{\mathbb{R}^{N}}(|(-\Delta)^{\frac{s}{2}}u|^{2} - \lambda u^{2})dx - \int_{\mathbb{R}^{N}}F(|x|,u)dx,
		$$
		$$
		S(u) = \frac{1}{2}\int_{\mathbb{R}^{N}}|(-\Delta)^{\frac{s}{2}}u|^{2}dx, G(u) = 0,
		$$
		$$
		F(u) = \int_{\mathbb{R}^{N}}F(|x|,u)dx, Q(u) = \frac{1}{2}\int_{\mathbb{R}^{N}}|u|^{2}dx,
		$$
		$$
		\mathcal{N}_{\lambda} = \{u \in H_{rad}^{s}(\mathbb{R}^{N}): \int_{\mathbb{R}^{N}}(|(-\Delta)^{\frac{s}{2}}u|^{2} - \lambda u^{2})dx - \int_{\mathbb{R}^{N}}f(|x|,u)udx = 0\}.
		$$
		
		\begin{remark}
			$H_{rad}^{s}(\mathbb{R}^{N})$ can be replaced by $H^{s}(\mathbb{R}^{N})$ throughout this section if we assume additional assumptions on the integrant ensuring that all ground states are radial and radially decreasing (Schwarz symmetric). The reader can refer to the following works for a more detailed account \cite{Hajaiej-Stuart1, Hajaiej1, Hajaiej2, Hajaiej3}.
		\end{remark}
		
		All assumptions given in the abstract framework 2 are satisfied and our main result in this subsection reads:
		
		\begin{theorem}
			\label{thm4.5} Assume $(f_{1})$ - $(f_{4})$ with $(N - 2s)\beta < 2(N + \theta)$. Furthermore, we also assume that (\ref{eq4.1}) admits at most one positive ground state in $H_{rad}^{s}(\mathbb{R}^{N})$ for any fixed $\lambda < 0$. If
			$$
			\alpha > \frac{2(N + \tau) + 4s}{N}
			$$
			or
			$$
			2 + \frac{2\tau}{N} < \alpha \leq \beta < \frac{2(N + \theta) + 4s}{N},
			$$
			then for any $\rho > 0$, (\ref{eq4.1}) admits at least one solution $(\lambda, u_{\lambda}) \in (-\infty, 0) \times H^{s}_{rad}(\mathbb{R}^{N})$ such that $\int_{\mathbb{R}^{N}}u_{\lambda}^{2}dx = \rho$.
		\end{theorem}
 
		\begin{remark}
            Using our approach, it seems challenging to show the uniqueness of positive GSS for a general function $f$. In the next section, 5.2, when $f$ takes some special form, namely $f(x,u)=h(|x|)u^{q-2},$ we will prove the uniqueness of GSS. Our result generalizes the breakthrough papers \cite{FL, FLS} by using a totally different method.
        \end{remark}
        
		\begin{lemma}
			\label{lem4.1} $(i)$ Assume $(f_{1})$ and $(f_{2})$, then $(F_{1})$ holds, i.e. for any $(\lambda, u) \in (-\infty, 0) \times H_{rad}^{s}(\mathbb{R}^{N}) \backslash \{0\}$, there exists a unique function $t: (-\infty, 0) \times H_{rad}^{s}(\mathbb{R}^{N}) \backslash \{0\} \rightarrow (0, +\infty)$ such that
			\begin{equation}
				when \ t > 0, \ then \ (\lambda, t u) \in \mathcal{N} \Leftrightarrow t = t(\lambda, u),
				\nonumber
			\end{equation}
			and $\Phi_{\lambda}(t(\lambda, u)u) = \max_{t > 0}\Phi_{\lambda}(tu)$, $t \in C((-\infty, 0) \times H_{rad}^{s}(\mathbb{R}^{N}) \backslash \{0\}, (0, +\infty))$, where $\mathcal{N} := \cup_{\lambda < 0}\mathcal{N}_{\lambda}$.
		\end{lemma}
		
		\textit{Proof.  }  Consider
		\begin{eqnarray}
			\psi: (0, \infty) \times (-\infty, 0) \times H_{rad}^{s}(\mathbb{R}^{N}) \backslash \{0\} \rightarrow \mathbb{R} \nonumber \\
			(t, \lambda, u) \mapsto \Phi_{\lambda}(tu). \nonumber
		\end{eqnarray}
		Then
		$$
		\psi_{t}(t, \lambda, u) = t \int_{\mathbb{R}^{N}}(|(-\Delta)^{\frac{s}{2}}u|^{2} - \lambda u^{2})dx - \int_{\mathbb{R}^{N}}f(|x|,tu)udx.
		$$
		Thus $(\lambda, tu) \in \mathcal{N}$ if and only if $\psi_{t}(t, \lambda, u) = 0$. Note that
		$$
		\psi_{tt}(t, \lambda, u) = \int_{\mathbb{R}^{N}}(|(-\Delta)^{\frac{s}{2}}u|^{2} - \lambda u^{2})dx - \int_{\mathbb{R}^{N}}f_{t}(|x|,tu)u^{2}dx.
		$$
		When $\psi_{t}(t, \lambda, u) = 0$, then
		$$
		\psi_{tt}(t, \lambda, u) = \int_{\mathbb{R}^{N}}(\frac{f(|x|,tu)}{t}u - f_{t}(|x|,tu)u^{2})dx < 0.
		$$
		By the fact that
		$$
		\Phi_{\lambda}(u) \geq \min\{\frac{1}{2}, -\frac{\lambda}{2}\}\|u\|^{2} - C(\|u\|^{\alpha} + \|u\|^{\beta}),
		$$
		it is easy to verify that $\psi(t, \lambda, u) \geq \delta > 0$ for suitable $t$ (since there exist $\rho > 0$, $\delta > 0$ such that $\inf_{\|v\| = \rho}\Phi_{\lambda}(v) \geq \delta$). Furthermore, $\psi(t, \lambda, u) \rightarrow -\infty$ as $t \rightarrow +\infty$ for fixed $\lambda$ and $u$, implying the existence and uniqueness of $t(\lambda, u)$.
		
		We see $\psi_{t}(t, \lambda, u)$ as $\psi_{t}(t, (\lambda, u))$, and $$\psi_{t(\lambda, u)}(t, (\lambda, u)) = (\psi_{t\lambda}(t, \lambda, u), \psi_{tu}(t, \lambda, u))$$ where
		$$
		\psi_{t\lambda}(t, \lambda, u) = -t \int_{\mathbb{R}^{N}}u^{2}dx,
		$$
		\begin{eqnarray}
			\psi_{tu}(t, \lambda, u)(\varphi) &=& 2t\int_{\mathbb{R}^{N}}((-\Delta)^{\frac{s}{2}}u (-\Delta)^{\frac{s}{2}}\varphi - \lambda u \varphi)dx \nonumber \\
			&-& \int_{\mathbb{R}^{N}}(f(|x|,tu) + tf_{t}(|x|,tu)u)\varphi dx, \forall \varphi \in H_{rad}^{s}(\mathbb{R}^{N}). \nonumber
		\end{eqnarray}
		Since $\psi_{tt}(t(\lambda, u), \lambda, u) < 0$, it is an isomorphism of $\mathbb{R}$ onto $\mathbb{R}$. By the Implicit Function Theorem, $t \in C^{1}((-\infty, 0) \times H_{rad}^{s}(\mathbb{R}^{N})\backslash \{0\}, (0, +\infty))$. \qed\vskip 5pt
		
		\begin{lemma}
			\label{lem4.2} Assume that $(f_{1})$ and $(f_{2})$ hold true, then $(F_{2})$ holds.
		\end{lemma}
		
		\textit{Proof.  } By Lemma \ref{lem4.1}, there exists $\delta_{1} = \delta_{1}(\lambda_{0}) > 0$ such that $h(\lambda) < h(\lambda_{0}) + 1$ for all $|\lambda - \lambda_{0}| < \delta_{1}$ and $\lambda < 0$. This, together with $u_{\lambda} \in \mathcal{K}_{\lambda}$, we have
		\begin{equation} \label{eq4.3}
			\left\{
			\begin{array}{cc}
				\Phi_{\lambda}(u_{\lambda}) = \int_{\mathbb{R}^{N}}(\frac{1}{2}(|(-\Delta)^{\frac{s}{2}} u_{\lambda}|^{2} - \lambda u_{\lambda}^{2})-F(|x|,u_{\lambda}))dx < h(\lambda_{0}) + 1, \\
				\langle \Phi'_{\lambda}(u_{\lambda}), u_{\lambda} \rangle = \int_{\mathbb{R}^{N}}(|(-\Delta)^{\frac{s}{2}} u_{\lambda}|^{2} - \lambda u_{\lambda}^{2}- f(|x|,u_{\lambda})u_{\lambda})dx = 0.
			\end{array}
			\right.
		\end{equation}
		Thus we have that
		\[
		\frac{1}{2}\int_{\mathbb{R}^{N}}f(|x|,u_{\lambda})u_{\lambda}dx - \int_{\mathbb{R}^{N}}F(|x|,u_{\lambda})dx < h(\lambda_{0}) + 1.
		\]
		Since $\alpha F(|x|,t) \leq f(|x|,t)t \leq \beta F(|x|,t)$, it follows that
		\[
		\int_{\mathbb{R}^{N}}f(|x|,u_{\lambda})u_{\lambda}dx < \frac{2 \alpha}{\alpha - 2}(h(\lambda_{0}) + 1),
		\]
		i.e.
		\[
		\int_{\mathbb{R}^{N}}(|(-\Delta)^{\frac{s}{2}} u_{\lambda}|^{2} - \lambda u_{\lambda}^{2})dx < \frac{2 \alpha}{\alpha - 2}(h(\lambda_{0}) + 1).
		\]
		Therefore, if $|\lambda - \lambda_{0}| < \min \{\delta_{1}, -\frac{\lambda_{0}}{2}\}$, which implies that
		\begin{eqnarray}
		&& \int_{\mathbb{R}^{N}}(|(-\Delta)^{\frac{s}{2}} u_{\lambda}|^{2} + u_{\lambda}^{2})dx \nonumber \\
		&<& (1 - \frac{1}{\lambda}) \int_{\mathbb{R}^{N}}(|(-\Delta)^{\frac{s}{2}} u_{\lambda}|^{2} - \lambda u_{\lambda}^{2})dx \nonumber \\
		&<& \frac{2 \alpha(\lambda_{0} - 2)}{\lambda_{0}(\alpha - 2)}(h(\lambda_{0}) + 1).
		\end{eqnarray}
		Take $\delta = \min \{\delta_{1}, -\frac{\lambda_{0}}{2}\}$. For all $\lambda$ with $|\lambda - \lambda_{0}| < \delta$, we have that:
		\[
		\|u_{\lambda}\| < \sqrt{\frac{2 \alpha(\lambda_{0} - 2)}{\lambda_{0}(\alpha - 2)}(h(\lambda_{0}) + 1)}.
		\]
		Therefore, $(F_{2})$ holds and the proof is completed. \qed\vskip 5pt
		
		\begin{lemma}
			\label{lem4.3} (Pohozaev identity) Suppose that $(f_{1})$, $(f_{2})$, $(f_{4})$ or $(f_{1})$ - $(f_{3})$ hold. Let $u \in H_{rad}^{s}(\mathbb{R}^{N})$ be a weak solution of (\ref{eq4.1}), then the following integral identity holds true:
			\begin{eqnarray} \label{eq4.4}
				&& (N-2s)\int_{\mathbb{R}^{N}}|(-\Delta)^{\frac{s}{2}} u|^{2}dx \nonumber \\
				&=& N\lambda\int_{\mathbb{R}^{N}}u^{2}dx + 2\int_{\mathbb{R}^{N}}[NF(|x|,u)dx + |x|F_{r}(|x|,u)]dx.
			\end{eqnarray}
		\end{lemma}
		
		\textit{Proof.  } We know that $u \in H_{rad}^{2s+1}(\mathbb{R}^{N})$ if $u \in H_{rad}^{s}(\mathbb{R}^{N})$ solves (\ref{eq4.1}). In fact, \cite[Lemma B.2]{FL} gives the proof when $N = 1$; for $N \geq 2$, the modifications of the proof given there are straightforward (the proof of $u \in L^\infty$ can be found in \cite[AppendixB.2]{FLS} and then the cases $s \geq 1/2$ and $0 < s <1/2$ can be treated separately as the proof of \cite[Lemma B.2]{FL}).  We integrate (\ref{eq4.1}) after a multiplication by $x \cdot \nabla u$. Integration by parts implies that
		\[
		\langle x \cdot \nabla u, (-\Delta)^{s}u\rangle_{L^{2}} = \frac{2s-N}{2}\langle u, (-\Delta)^{s}u\rangle_{L^{2}},
		\langle x \cdot \nabla u, \lambda u\rangle_{L^{2}} = \frac{-N\lambda}{2}\langle u, u\rangle_{L^{2}},
		\]
		\[
		\langle x \cdot \nabla u, f(|x|,u)\rangle_{L^{2}} = -\int_{\mathbb{R}^{N}}[NF(|x|,u) + |x|F_{r}(|x|,u)]dx,
		\]
		where for the first identity we also use that $[\nabla \cdot x, (-\Delta)^{s}] = -2s(-\Delta)^{s}$, which is easily  verified in Fourier space. Hence, we deduce (\ref{eq4.4}). \qed\vskip 5pt
		
		\begin{lemma}
			\label{lem4.4}  Assume $(f_{1})$ - $(f_{4})$ with $(N - 2s)\beta < 2(N + \theta)$. If
			$$
			\alpha > \frac{2(N + \tau) + 4s}{N},
			$$
			then $(F_{5})$ and $(F_{6})$ hold with $l \geq k > 1$. If
			$$
			2 + \frac{2\tau}{N} < \alpha \leq \beta < \frac{2(N + \theta) + 4s}{N},
			$$
			then $(F_{5})$ and $(F_{6})$ hold with $0 < k \leq l < 1$.
		\end{lemma}
		
		\textit{Proof.  } By integrating (\ref{eq4.1}) after a multiplication by $u$, we obtain
		\begin{equation} \label{eq4.5}
			\int_{\mathbb{R}^{N}}|(-\Delta)^{\frac{s}{2}}u|^{2}dx = \lambda\int_{\mathbb{R}^{N}}u^{2}dx + \int_{\mathbb{R}^{N}}f(|x|,u)udx.
		\end{equation}
		$(f_{3})$, (\ref{eq4.4}) and (\ref{eq4.5}) imply that:
		\begin{eqnarray} \label{eq4.6}
			&& -2s\lambda\int_{\mathbb{R}^{N}}u^{2}dx \nonumber \\
			&=& 2\int_{\mathbb{R}^{N}}[NF(|x|,u) + |x|F_{r}(|x|,u)]dx - (N-2s)\int_{\mathbb{R}^{N}}f(|x|,u)udx \nonumber \\
			&\geq& [2(N+\theta)-(N-2s)\beta]\int_{\mathbb{R}^{N}}F(|x|,u)dx,
		\end{eqnarray}
		\begin{eqnarray} \label{eq4.7}
			&& 2s\int_{\mathbb{R}^{N}}|(-\Delta)^{\frac{s}{2}}u|^{2}dx \nonumber \\
			&=& N\int_{\mathbb{R}^{N}}f(|x|,u)udx - 2\int_{\mathbb{R}^{N}}[NF(|x|,u) + |x|F_{r}(|x|,u)]dx \nonumber \\
			&\leq& [N\beta - 2(N+\theta)]\int_{\mathbb{R}^{N}}F(|x|,u)dx.
		\end{eqnarray}
		$(f_{4})$, (\ref{eq4.4}) and (\ref{eq4.5}) imply that:
		\begin{eqnarray} \label{eq4.8}
			&& -2s\lambda\int_{\mathbb{R}^{N}}u^{2}dx \nonumber \\
			&=& 2\int_{\mathbb{R}^{N}}[NF(|x|,u) + |x|F_{r}(|x|,u)]dx - (N-2s)\int_{\mathbb{R}^{N}}f(|x|,u)udx \nonumber \\
			&\leq& [2(N+\tau)-(N-2s)\alpha]\int_{\mathbb{R}^{N}}F(|x|,u)dx,
		\end{eqnarray}
		\begin{eqnarray} \label{eq4.9}
			&& 2s\int_{\mathbb{R}^{N}}|(-\Delta)^{\frac{s}{2}}u|^{2}dx \nonumber \\
			&=& N\int_{\mathbb{R}^{N}}f(|x|,u)udx - 2\int_{\mathbb{R}^{N}}[NF(|x|,u) + |x|F_{r}(|x|,u)]dx \nonumber \\
			&\geq& [N\alpha - 2(N+\tau)]\int_{\mathbb{R}^{N}}F(|x|,u)dx.
		\end{eqnarray}
		Therefore, if $\beta > \frac{2(N + \theta) + 4s}{N}$ and $2(N + \theta) -\beta(N-2s) > 0$, then
		\begin{eqnarray} \label{eq4.10}
			&& \Phi_{\lambda}(u) \nonumber \\
			&=& \frac{1}{2}\int_{\mathbb{R}^{N}}(|(-\Delta)^{\frac{s}{2}} u|^{2} - \lambda u^{2})dx - \int_{\mathbb{R}^{N}}F(|x|,u)dx \nonumber \\
			&\leq& [\frac{N\beta - 2(N + \theta)}{4s}-1]\int_{\mathbb{R}^{N}}F(|x|,u)dx - \frac{\lambda}{2}\int_{\mathbb{R}^{N}}u^{2}dx \nonumber \\
			&\leq& \{[\frac{N\beta - 2(N + \theta)}{4s}-1]\frac{-2s\lambda}{2(N + \theta) -\beta(N-2s)} - \frac{\lambda}{2}\}\int_{\mathbb{R}^{N}}u^{2}dx,
		\end{eqnarray}
		showing that:
		\begin{equation}
			\Phi_{\lambda}(u) \leq l \frac{-\lambda}{2}\int_{\mathbb{R}^{N}}u^{2}dx,
			\nonumber
		\end{equation}
		where \(l = \frac{N\beta - 2(N + \theta) - 4s}{2(N + \theta) -\beta(N-2s)} + 1 > 1\). Thus, $(F_{5})$ holds with $l > 1$. \\
		Similarly, if $2 < \beta < \frac{2(N + \theta) + 4s}{N}$ and $2(N + \tau) - \alpha(N-2s) > 0$, then
		\begin{eqnarray} \label{eq4.11}
			&& \Phi_{\lambda}(u) \nonumber \\
			&\leq& [\frac{N\beta - 2(N + \theta)}{4s}-1]\int_{\mathbb{R}^{N}}F(|x|,u)dx - \frac{\lambda}{2}\int_{\mathbb{R}^{N}}u^{2}dx \nonumber \\
			&\leq& \{[\frac{N\beta - 2(N + \theta)}{4s}-1]\frac{-2s\lambda}{2(N + \tau) - \alpha(N-2s)} - \frac{\lambda}{2}\}\int_{\mathbb{R}^{N}}u^{2}dx,
		\end{eqnarray}
		showing that:
		\begin{equation}
			\Phi_{\lambda}(u) \leq l \frac{-\lambda}{2}\int_{\mathbb{R}^{N}}u^{2}dx,
			\nonumber
		\end{equation}
		where \(l = \frac{N\beta - 2(N + \theta) - 4s}{2(N + \tau) - \alpha(N-2s)} + 1 \in (0,1)\). Thus, $(F_{5})$ holds with $l \in (0,1)$. \\
		If $\alpha > \frac{2(N + \tau) + 4s}{N}$ and $2(N + \tau) - \alpha(N-2s) > 0$, then
		\begin{eqnarray} \label{eq4.12}
			&& \Phi_{\lambda}(u) \nonumber \\
			&=& \frac{1}{2}\int_{\mathbb{R}^{N}}(|(-\Delta)^{\frac{s}{2}} u|^{2} - \lambda u^{2})dx - \int_{\mathbb{R}^{N}}F(|x|,u)dx \nonumber \\
			&\geq& [\frac{N\alpha - 2(N + \tau)}{4s}-1]\int_{\mathbb{R}^{N}}F(|x|,u)dx - \frac{\lambda}{2}\int_{\mathbb{R}^{N}}u^{2}dx \nonumber \\
			&\geq& \{[\frac{N\alpha - 2(N + \tau)}{4s}-1]\frac{-2s\lambda}{2(N + \tau) - \alpha(N-2s)} - \frac{\lambda}{2}\}\int_{\mathbb{R}^{N}}u^{2}dx,
		\end{eqnarray}
		proving that
		\begin{equation}
			\Phi_{\lambda}(u) \geq k \frac{-\lambda}{2}\int_{\mathbb{R}^{N}}u^{2}dx,
			\nonumber
		\end{equation}
		where \(k = \frac{N\alpha - 2(N + \tau) - 4s}{2(N + \tau) - \alpha(N-2s)} + 1 > 1\). Hence, $(F_{6})$ holds with $k > 1$. \\
		If $2 + \frac{2\tau}{N} < \alpha < \frac{2(N + \tau) + 4s}{N}$ and $2(N + \theta) - \beta(N-2s) > 0$, then
		\begin{eqnarray} \label{eq4.13}
			&& \Phi_{\lambda}(u) \nonumber \\
			&\geq& [\frac{N\alpha - 2(N + \tau)}{4s}-1]\int_{\mathbb{R}^{N}}F(|x|,u)dx - \frac{\lambda}{2}\int_{\mathbb{R}^{N}}u^{2}dx \nonumber \\
			&\geq& \{[\frac{N\alpha - 2(N + \tau)}{4s}-1]\frac{-2s\lambda}{2(N + \theta) - \beta(N-2s)} - \frac{\lambda}{2}\}\int_{\mathbb{R}^{N}}u^{2}dx,
		\end{eqnarray}
		which proves that:
		\begin{equation}
			\Phi_{\lambda}(u) \leq k\frac{-\lambda}{2}\int_{\mathbb{R}^{N}}u^{2}dx,
			\nonumber
		\end{equation}
		where \(k = \frac{N\alpha - 2(N + \tau) - 4s}{2(N + \theta) - \beta(N-2s)} + 1 \in (0,1)\). Thus, $(F_{6})$ holds with $k \in (0,1)$.
		
		The proof is complete. \qed\vskip 5pt
		
		\textbf{Proof of Theorem \ref{thm4.5}.  } Lemmas \ref{lem4.1}, \ref{lem4.2} and \ref{lem4.4} yield that $(F_{1})$, $(F_{2})$, $(F_{5})$ and $(F_{6})$ hold. Standard methods yield that $(F_3)$ and $(F_4)$ holds. Since we assume that $f(r,t) \geq 0$ for all $r > 0$, $t \in \mathbb{R}$, nontrivial solutions of $(\ref{eq4.1})$ are strictly positive by the maximum principle, see e.g. \cite[Lemma 4.13 and Remark 4.14]{CS} when $0 < s < 1$. We also assume that (\ref{eq4.1}) admits one positive ground state in $H_{rad}^{s}(\mathbb{R}^{N})$. Therefore, $\mathcal{K}_{\lambda}$ has a unique element $u_{\lambda}$ for any $\lambda < 0$, i.e. $\mathcal{K}_{\lambda} = \{u_{\lambda}\}$. Thus $d \equiv \hat{d}$ in $(-\infty,0)$. Then, apply Corollary \ref{cor2.4} $(iv)$ and the proof is complete. \qed\vskip 5pt
		
		\begin{remark}
			\label{rmk4.6} In Theorem \ref{thm4.5}, we have $0 \leq (N-2s)\tau < N\theta + 4s^{2} \leq 4$ when $\alpha > \frac{2(N + \tau) + 4s}{N}$ and $(N - 2s)\beta < 2(N + \theta)$ since
			\[
			2(N + \theta) > (N - 2s)\beta \geq (N - 2s)\alpha > (N - 2s)\frac{2(N + \tau) + 4s}{N};
			\]
			and $0 \leq \tau < \theta + 2s \leq 2s$ when $2 + \frac{2\tau}{N} < \alpha \leq \beta < \frac{2(N + \theta) + 4s}{N}$.
		\end{remark}
		
		\begin{remark}
			\label{rmk4.7}  When $f_{r}(r,t) \leq 0$ for all $r > 0$ and $t > 0$, applying the moving plane arguments developed in \cite{MZ} directly, we can show that any positive solution of (\ref{eq4.1}) is radially symmetric and monotonically decreasing with respect to $r = |x|$. Since we discuss in $H_{rad}^{s}(\mathbb{R}^{N})$, we do not need the radially symmetric result here, while the monotonic result will be used in the proof of non-degeneracy, see Lemma \ref{lem4.13} below.
		\end{remark}
		
		\begin{remark}
			\label{rmk4.8} In the special case when
			$$
			f(|r|,u) = |u|^{q - 2}u, q \in (2, 2_{s}^{\ast}),
			$$
			we can show the existence of normalized solutions without the hypothesis of uniqueness for the positive ground state solution, see Corollary \ref{cor2.4} $(v)$.
			
			In fact, when $q \neq 2 + 4s/N$, then for any $\rho > 0$, (\ref{eq4.1}) admits a solution $(\lambda, u_{\lambda}) \in (-\infty, 0) \times H^{s}_{rad}(\mathbb{R}^{N})$ such that $\int_{\mathbb{R}^{N}}u_{\lambda}^{2}dx = 2\rho$. Furthermore,
			$$
			\rho = C\frac{2N - (N-2s)q}{2s(q-2)}(-\lambda)^{\frac{4s-N(q-2)}{2s(q-2)}},
			$$
			which can be found in Theorem \ref{thm2.2} $(iv)$. When $q = 2 + 4s/N$, $\int_{\mathbb{R}^{N}}u_{\lambda}^{2}dx \equiv Const.$ for all $\lambda \in (-\infty,0)$ and $u_{\lambda} \in \mathcal{K}_{\lambda}$.
		\end{remark}
		
		\subsection{A non-autonomous case without the hypothesis of uniqueness of the ground state solution} \label{non-aut}
		
		Here we come back to the uniqueness assumption that we made in Theorem \ref{thm4.5}, where we assumed that (\ref{eq4.1}) admits a unique positive ground state in $H_{rad}^{s}(\mathbb{R}^{N})$ for any fixed $\lambda < 0$. When $s = 1$ and $f$ is autonomous, there are many results about the uniqueness of positive solution. When $f$ is non-autonomous, there are quite few results. When $0 < s < 1$ and $f(|x|,u) = |u|^{q-2}u$, the uniqueness of positive ground state has been shown in \cite{FL, FLS}. However, to the best of our knowledge, there are no results when $f$ is non-autonomous. Hence, in this subsection, we will apply our abstract framework 1 to a class of non-autonomous $f$ without the uniqueness assumption. Then as a corollary, we will show the uniqueness of positive ground state can be proved when $0 < s \leq 1$. This uniqueness result is new. More precisely, consider
		\begin{equation} \label{eq4.14}
			\left\{
			\begin{array}{cc}
				(-\Delta)^{s} u = \lambda u + h(|x|)|u|^{q-2}u \ in \ \mathbb{R}^{N},\\
				u(x) \rightarrow 0 \ as \ |x| \rightarrow +\infty,
			\end{array}
			\right.
		\end{equation}
		where $0 < s \leq 1$, $N \geq 2$, $\lambda < 0$, $2 < q < 2_{s}^{\ast}$, and we assume on $h(r)$:
		
		$(h_{1})$ $h(r) \in C^{1}[0,+\infty) \cap L^{\infty}[0,+\infty)$, $h(r) > 0$ in $[0,+\infty)$, $h(r)$ and $\frac{rh'(r)}{h(r)}$ are non-increasing in $(0,+\infty)$, $\theta = \lim_{r \rightarrow +\infty}\frac{rh'(r)}{h(r)} > -\infty$.
		
		Our main result in this subsection is the following theorem:

		\begin{theorem}
			\label{cor4.19}  Assume that $(h_{1})$ holds with $(N-2s)q < 2(N+\theta)$. Then for any $\rho > 0$, (\ref{eq4.14}) admits at least one solution $(\lambda, u_{\lambda}) \in (-\infty, 0) \times H^{s}_{rad}(\mathbb{R}^{N})$ such that $\int_{\mathbb{R}^{N}}u_{\lambda}^{2}dx = 2\rho$ if $q > 2 + 4s/N$ or $2< q < 2 +(2\theta + 4s)/N$.
		\end{theorem}
		
		\begin{remark}
			\label{cor4.20}  Let $h(r) = 1/(1+r^{k})^{l}$, $k \geq 1$, $l > 0$. Then $(h_{1})$ holds with $\theta = -kl$.
		\end{remark}
		
		We would like to use the abstract framework 1. Let $E = L^2(\mathbb{R}^{N}), W = H^{s}(\mathbb{R}^{N}), \widehat{W} = H^{2s}(\mathbb{R}^{N}), G = O(N)$. It is not difficult to verify $(H_1)$ - $(H_4)$, $(H_6)$, $(H_7)$. Moreover, a $G$-ground state is indeed a ground state here. To show $(H_5)$, we prove the following lemma.		
		
		\begin{lemma}
			\label{lem4.13b}  Let $u_{\lambda} \in H_{rad}^{s}(\mathbb{R}^{N})$ be a positive solution of \eqref{eq4.14} and $\mu_{G}(u_{\lambda}) = 1$. Suppose that $(h_{1})$ holds true, then $\ker D_{uu}\Phi_\lambda(u_{\lambda})|_{L^2_{rad}} = \{0\}$.
		\end{lemma}
		
		\textit{Proof.  }  Note that $\inf\sigma_{ess}(D_{uu}\Phi_\lambda(u_{\lambda})) = -\lambda > 0$. Arguing by contradiction, we assume that $0 \in \sigma_{dis}(D_{uu}\Phi_\lambda(u_{\lambda})|L^{2}_{rad})$. Since $\mu_{G}(u_{\lambda}) = 1$, $0$ is the second eigenvalue. Assume that $v \in H^{2s}_{rad}$ is an unit eigenfunction corresponding to $0$, then $v = v(|x|)$ changes its sign exactly once for $|x| = r \in (0,+\infty)$, which is a standard result when $s = 1$ and can be found in \cite[Theorem 2]{FLS} when $0 < s < 1$. We assume that $v(r^{\ast}) = 0$, $v(r) \geq 0$ in $(0,r^{\ast})$ and $v(r) \leq 0$ in $(r^{\ast},+\infty)$.
		
		Direct calculation yields to
		\begin{equation}
		D_{uu}\Phi_\lambda(u_{\lambda})(u_{\lambda}) = (2-q)h(r)u_{\lambda}^{q-1},
		\end{equation}
		\begin{equation}
		D_{uu}\Phi_\lambda(u_{\lambda})(x \cdot \nabla u_{\lambda}) = (rh'(r) + 2sh(r))u_{\lambda}^{q-1} + 2s\lambda u_{\lambda}.
		\end{equation}
		Let
		\begin{eqnarray}
		z(r) &=& (rh'(r) + 2sh(r))u_{\lambda}^{q-1} + \lambda u_{\lambda} + k h(r)u_{\lambda}^{q-1} \nonumber \\
		&=& h(r)u_{\lambda}^{q-1}[k + 2s + \frac{rh'(r)}{h(r)} + \frac{2s\lambda}{h(r)u_{\lambda}^{q-2}}] \nonumber,
		\end{eqnarray}
		where
		$$
		k = -(2s + \frac{r^{\ast}h'(r^{\ast})}{h(r^{\ast})} + \frac{2s\lambda}{h(r^{\ast})u_{\lambda}(r^{\ast})^{q-2}}),
		$$
		ensuring that $z(r^{\ast}) = 0$. By Remark \ref{rmk4.7}, $u_{\lambda}(r)$ is monotonically decreasing. By $(h_{1})$, $h(r)$ and $rh'(r)/h(r)$ are non-increasing. Thus $z(r) > 0$ in $(0,r^{\ast})$ and $z(r) < 0$ in $(r^{\ast},+\infty)$, yielding $\langle z, v\rangle_{L^{2}} > 0$. However, since $z = D_{uu}\Phi_\lambda(u_{\lambda})(\frac{k}{2-q}u_{\lambda} + x \cdot \nabla u_{\lambda})$, we obtain
		\begin{eqnarray}
		\langle z, v\rangle_{L^{2}} &=& \langle D_{uu}\Phi_\lambda(u_{\lambda})(\frac{k}{2-q}u_{\lambda} + x \cdot \nabla u_{\lambda}), v\rangle_{L^{2}} \nonumber \\
		&=& \langle \frac{k}{2-q}u_{\lambda} + x \cdot \nabla u_{\lambda}, D_{uu}\Phi_\lambda(u_{\lambda})(v)\rangle_{L^{2}} = 0,
		\end{eqnarray}
		which is a contradiction, and this completes the proof. \qed\vskip 5pt	
		
		We assume the positivity of $u_{\lambda}$ in Lemma \ref{lem4.13b}. Thus, to apply our abstract framework 1, we also need the following lemma.	
		
	    \begin{lemma}
	    	\label{lemB.11b} There exists
	    	$$
	    	\lambda \mapsto u_\lambda \in C^1([\lambda^\ast,\lambda_\ast],W \setminus \{0\})
	    	$$
	    	for some $\lambda^\ast < \lambda_\ast < 0$ where $u_\lambda$ is a solution of \eqref{eq4.14}.
	    	
	    	$(i)$ If $u_{\lambda_\ast} > 0$, then $u_\lambda > 0$ for $\lambda$ close enough to $\lambda_\ast$.
	    	
	    	$(ii)$ Let $\{\lambda_n\} \subset (\lambda^\ast,\lambda_\ast)$ be a sequence such that $\lambda_n \rightarrow \lambda^\ast$, and $u_{\lambda_n} > 0$. Then $u_{\lambda^\ast} > 0$.
	    \end{lemma}
	
	    \textit{Proof.} Along the lines of Lemma \ref{lem4.16} $(i)$, we can prove $(i)$, $(ii)$. \qed\vskip 5pt
	
	    \textbf{Proof of Theorem \ref{cor4.19}.}  Thanks to Lemmas \ref{lem4.13b} and \ref{lemB.11b}, Theorem \ref{thmB.1} shows that there exists
	    $$
	    \lambda \mapsto u_\lambda \in C^1((-\infty,0), H^{s}_{rad}(\mathbb{R}^{N}) \setminus \{0\}),
	    $$
	    such that $u_\lambda$ solves \eqref{eq4.14}, and $\mu_{G}(u_\lambda) = 1$.
	
        Then by Remark \ref{rmk2.5}, following the proof of Theorem \ref{thm4.5}, we can use abstract framework 2 to obtain the asymptotical behaviors of $\int_{\mathbb{R}^{N}}u_\lambda^{2}dx$ as $\lambda \rightarrow -\infty$ and $\lambda \rightarrow 0^-$ to complete the proof. \qed\vskip 5pt
		
		Finally in this subsection, we show a uniqueness result. The key step is to show a local uniqueness result when $\lambda \rightarrow -\infty$.
		
		\begin{lemma}
			\label{lem4.9b}  Assume $(h_{1})$ holds with $(N-2s)q < 2(N+\theta)$. Then there exists $\Lambda_1 < 0$ such that the positive solution of $(\ref{eq4.14})$ with $G$-Morse index $1$ is unique in $H_{rad}^{s}(\mathbb{R}^{N})$ for all fixed $\lambda < \Lambda_1$.
		\end{lemma}
		
		\textit{Proof.} Assume that for some $\lambda_0 < 0$, there is a positive solution $u_{\lambda_0}$ of $(\ref{eq4.14})$ with $\mu_G(u_{\lambda_0}) = 1$. Similar to the proof of Theorem \ref{cor4.19}, we know that there exists
		$$
		\lambda \mapsto u_\lambda \in C^1((-\infty,\lambda_0], H^{s}_{rad}(\mathbb{R}^{N}) \setminus \{0\}),
		$$
		such that $u_\lambda$ solves \eqref{eq4.14}, and $\mu_{G}(u_\lambda) = 1$. Using ideas in Section \ref{AF3.2}, let $\mu = 1/\lambda$ and consider	
		\begin{equation}
		v_{\mu} = |\mu|^{\frac{1}{p-2}}u_{\lambda}(|\mu|^{\frac{1}{2s}}x),
		\end{equation}
		which satisfies the following equation
		\begin{equation} \label{eqv}
		(-\Delta)^{s} v_{\mu} + v_{\mu} = h(|\mu|^{\frac{1}{2s}}|x|)|v_{\mu}|^{q-2}v_{\mu}.
		\end{equation}
		Moreover, there exists
		$$
		\mu \mapsto v_\mu \in C^1([\frac{1}{\lambda_0},0), H^{s}_{rad}(\mathbb{R}^{N}) \setminus \{0\}),
		$$
		where $\mu_{G}(v_\mu) = 1$. Then along the lines of the proof for Lemma \ref{lem4.15}, we have the following power estimates
		\[
		\int_{\mathbb{R}^{N}}h(|\mu|^{\frac{1}{2s}}|x|)|v_\mu|^{q}dx \thicksim \int_{\mathbb{R}^{N}}|v_\mu|^{2}dx \thicksim \int_{\mathbb{R}^{N}}|(-\Delta)^{\frac{s}{2}}v_\mu|^{2}dx \gtrsim 1,
		\]
		Next, we show that $\int_{\mathbb{R}^{N}}h(|\mu|^{\frac{1}{2s}}|x|)|v_\mu|^{q}dx \lesssim 1$. Differentiating \eqref{eqv} with respect to $\mu$, we get
		\begin{equation}
		L_{\mu}\partial_\mu u_\mu = \frac{1}{2s\mu}h'(|\mu|^{\frac{1}{2s}}|x|)|\mu|^{\frac{1}{2s}}|x||v_{\mu}|^{q-2}v_{\mu},
		\end{equation}
		showing that
		$$
		\partial_\mu v_\mu = L_{\mu}^{-1}(\frac{1}{2s\mu}h'(|\mu|^{\frac{1}{2s}}|x|)|\mu|^{\frac{1}{2s}}|x||v_{\mu}|^{q-2}v_{\mu}),
		$$
		where
		$$
		L_{\mu} = (-\Delta)^{s} + 1 - (q-1) h(|\mu|^{\frac{1}{2s}}|x|)|v_{\mu}|^{q-2}.
		$$
		Since $L_{\mu}v_{\mu} = (2-q)h(|\mu|^{\frac{1}{2s}}|x|)|v_{\mu}|^{q-2}v_{\mu}$,
		\begin{eqnarray}
		&& \langle h(|\mu|^{\frac{1}{2s}}|x|)|v_{\mu}|^{q-2}v_{\mu}, \partial_\mu v_\mu\rangle_{L^{2}} \nonumber \\
		&=& \langle L_{\mu}^{-1}(h(|\mu|^{\frac{1}{2s}}|x|)|v_{\mu}|^{q-2}v_{\mu}), \frac{1}{2s\mu}h'(|\mu|^{\frac{1}{2s}}|x|)|\mu|^{\frac{1}{2s}}|x||v_{\mu}|^{q-2}v_{\mu}\rangle_{L^{2}} \nonumber \\
		&=& \frac{1}{(2-q)2s\mu}\int_{\mathbb{R}^{N}}h'(|\mu|^{\frac{1}{2s}}|x|)|\mu|^{\frac{1}{2s}}|x||v_{\mu}|^{q}dx.
		\end{eqnarray}
		We write $\phi(\mu) = \int_{\mathbb{R}^{N}}h(|\mu|^{\frac{1}{2s}}|x|)|v_\mu|^{q}dx$. Then
		\begin{eqnarray}
		&& \phi'(\mu) \nonumber \\
		&=& \frac{1}{2s\mu}\int_{\mathbb{R}^{N}}h'(|\mu|^{\frac{1}{2s}}|x|)|\mu|^{\frac{1}{2s}}|x||v_{\mu}|^{q}dx \nonumber \\
		&&+ q\langle h(|\mu|^{\frac{1}{2s}}|x|)|v_{\mu}|^{q-2}v_{\mu}, \partial_\mu v_\mu\rangle_{L^{2}} \nonumber \\
		&=& -\frac{\theta}{(q-2)s\mu}\int_{\mathbb{R}^{N}}h'(|\mu|^{\frac{1}{2s}}|x|)|\mu|^{\frac{1}{2s}}|x||v_{\mu}|^{q}dx \nonumber \\
		&\leq& 0.
		\end{eqnarray}
		Thus $\phi(\mu) \leq \phi(1/\lambda_0)$ for all $\mu \in [1/\lambda_0,0)$. Now we show that
		\[
		\int_{\mathbb{R}^{N}}h(|\mu|^{\frac{1}{2s}}|x|)|v_\mu|^{q}dx \thicksim \int_{\mathbb{R}^{N}}|v_\mu|^{2}dx \thicksim \int_{\mathbb{R}^{N}}|(-\Delta)^{\frac{s}{2}}v_\mu|^{2}dx \thicksim 1,
		\]
		implying that $v_\mu$ is uniformly bounded in $H^{s}(\mathbb{R}^{N})$. Then along the lines of the proof for Theorem \ref{thmc.4}, we can complete the proof.
		\qed\vskip 5pt
		
		\begin{theorem}
			\label{thm4.9}  Assume $(h_{1})$ holds with $(N-2s)q < 2(N+\theta)$. Then the positive ground state solution of $(\ref{eq4.14})$ is unique in $H^{s}(\mathbb{R}^{N})$ for any fixed $\lambda < 0$.
		\end{theorem}
		
		\textit{Proof.} By Remark \ref{rmk4.7}, any positive solution of (\ref{eq4.14}) is radially symmetric if $(h_{1})$ holds. Then the proof can be completed by the combination of Corollary \ref{cor uni} and Lemma \ref{lem4.9b}.
		\qed\vskip 5pt
		
		\begin{remark}
			It seems that the results in Theorem \ref{cor4.19} can be extended. We further assume that $h(r) = r^\theta, r > R$ for some $R > 0$. Similar to the proof of Lemma \ref{lem4.9b}, our abstract framework 3 is helpful to study the asymptotical behavior of $u_\lambda$ as $\lambda \rightarrow 0^-$. Then we believe that the following results hold:
			
			$(i)$ if $q > 2 + 4s/N$ or $2< q < 2 +(2\theta + 4s)/N$, for any $\rho > 0$, (\ref{eq4.14}) admits at least one solution $(\lambda, u_{\lambda}) \in (-\infty, 0) \times H^{s}_{rad}(\mathbb{R}^{N})$ such that $\int_{\mathbb{R}^{N}}u_{\lambda}^{2}dx = 2\rho$;
			
			$(ii)$ if $q = 2 + 4s/N$, there exists $\rho_1 > 0$ such that \\
			- for any $\rho > \rho_1$, (\ref{eq4.14}) admits at least one solution $(\lambda, u_{\lambda}) \in (-\infty, 0) \times H^{s}_{rad}(\mathbb{R}^{N})$ such that $\int_{\mathbb{R}^{N}}u_{\lambda}^{2}dx = 2\rho$, \\
			- for any $\rho < \rho_1$, (\ref{eq4.14}) admits no ground state solution $u_{\lambda}$ such that $\int_{\mathbb{R}^{N}}u_{\lambda}^{2}dx = 2\rho$;
			
			$(iii)$ if $2 +(2\theta + 4s)/N < q < 2 + 4s/N$, there exists $\rho_2 > 0$ such that \\
			- (\ref{eq4.14}) admits at least one solution $(\lambda, u_{\lambda}) \in (-\infty, 0) \times H^{s}_{rad}(\mathbb{R}^{N})$ such that $\int_{\mathbb{R}^{N}}u_{\lambda}^{2}dx = 2\rho_2$, \\
			- for any $\rho > \rho_2$, (\ref{eq4.14}) admits at least two solutions $(\lambda, u_{\lambda}), (\widehat{\lambda}, u_{\widehat{\lambda}}) \in (-\infty, 0) \times H^{s}_{rad}(\mathbb{R}^{N})$ such that $\int_{\mathbb{R}^{N}}u_{\lambda}^{2}dx = \int_{\mathbb{R}^{N}}u_{\widehat{\lambda}}^{2}dx = 2\rho$, \\
			- for any $\rho < \rho_2$, (\ref{eq4.14}) admits no ground state solution $u_{\lambda}$ such that $\int_{\mathbb{R}^{N}}u_{\lambda}^{2}dx = 2\rho$;
			
			$(iv)$ if $q = 2 + (2\theta + 4s)/N$, there exists $\rho_3 > 0$ such that \\
			- for any $\rho > \rho_3$, (\ref{eq4.14}) admits at least one solution $(\lambda, u_{\lambda}) \in (-\infty, 0) \times H^{s}_{rad}(\mathbb{R}^{N})$ such that $\int_{\mathbb{R}^{N}}u_{\lambda}^{2}dx = 2\rho$, \\
			- for any $\rho < \rho_3$, (\ref{eq4.14}) admits no ground state solution $u_{\lambda}$ such that $\int_{\mathbb{R}^{N}}u_{\lambda}^{2}dx = 2\rho$.
			
			Since this is a purely computational aspect, we are not going to prove in this quite long paper.
		\end{remark}
		
		\begin{remark}
			In Appendix \ref{Appen1}, we will provide another method to prove Theorem \ref{thm4.9}. The uniqueness problem of fractional nonlinear Schr\"{o}dinger equations is a long standing question. The breakthrough work in \cite{FL, FLS}, gave an affirmative answer for pure nonlinearity $|u|^{q-2}u$. The discussions in Appendix \ref{Appen1} are inspired by \cite{FL} and make important extensions.
		\end{remark}
		
		\subsection{The orbital stability/instability}
		
		If $(\lambda, u_{\lambda})$ is a solution of $(\ref{eq4.1})$, then $e^{-i\lambda t}u_{\lambda}$ is a standing wave solution of the fractional nonlinear Schr\"{o}dinger equation (fNLS)
		
		\begin{equation} \label{eq4.35}
			\left\{
			\begin{array}{cc}
				i\partial_{t} \psi - (-\Delta)^{s} \psi + \tilde{f}(x,\psi) = 0, (t,x) \in [0,+\infty) \times \mathbb{R}^{N}, \\
				\psi(0,x) = \psi_{0}(x),
			\end{array}
			\right.
		\end{equation}
		where $\tilde{f}: \mathbb{R}^{N} \times \mathbb{C} \rightarrow \mathbb{C}$ is defined by $\tilde{f}(x,e^{i\theta}u) = e^{i\theta}f(x,u)$, $u \in \mathbb{R}$. In this subsection, we study the orbital stability or instability of the ground state solutions for (\ref{eq4.1}) which are shown in Theorem \ref{thm4.5}. Recall that such solutions are called orbitally stable if for each $\epsilon > 0$ there exists $\delta > 0$ such that, whenever $\psi_{0} \in H_{rad}^{s}(\mathbb{R}^{N}, \mathbb{C})$ is such that $\|\psi_{0} - u_{\lambda}\|_{H_{rad}^{s}(\mathbb{R}^{N}, \mathbb{C})} < \delta$ and $\psi(t, x)$ is the solution of (\ref{eq4.35}) with $\psi(0, \cdot) = \psi_{0}$ in some interval $[0, t_{0})$, then $\psi(t, \cdot)$ can be continued to a solution in $0 \leq t < \infty$ and
		\[
		\sup_{0 < t < \infty}\inf_{w \in \mathbb{R}} \|\psi(t, x) - e^{-i\lambda w}u_{\lambda}(x)\|_{H_{rad}^{s}(\mathbb{R}^{N}, \mathbb{C})} < \epsilon;
		\]
		otherwise, they are called unstable. To do this, we assume that (\ref{eq4.35}) is locally well-posed in $H_{rad}^{s}(\mathbb{R}^{N},\mathbb{C})$:
		
		$(LWP)$ For each $\psi_{0} \in H_{rad}^{s}(\mathbb{R}^{N}, \mathbb{C})$, there exist $t_{0} > 0$, only depending on $\|\psi_{0}\|_{H^{s}(\mathbb{R}^{N}, \mathbb{C})}$, and a unique solution $\psi(t, x)$ of (\ref{eq4.35}) with initial datum $\psi_{0}$ in the interval $I = [0, t_{0})$. A detailed account on this aspect is given in \cite{HMOW}.
		
		Our main result in this subsection is as follows.
		
		\begin{theorem}
			\label{thm4.21} Under the hypotheses of Theorem \ref{thm4.5} and $(LWP)$ holds and any positive ground state of (\ref{eq4.1}) is non-degenerate in $H_{rad}^{s}(\mathbb{R}^{N})$ for all $\lambda \in (-\infty, 0)$. $(\lambda, u_{\lambda})$ is given by Theorem \ref{thm4.5}. Then for a.e. $\rho \in (0, \infty)$, $e^{-i\lambda t}u_{\lambda}$ is orbitally unstable in $H_{rad}^{s}(\mathbb{R}^{N}, \mathbb{C})$ when $\alpha > 2 + (2\tau + 4s)/N$ while $e^{-i\lambda t}u_{\lambda}$ is orbitally stable in $H_{rad}^{s}(\mathbb{R}^{N}, \mathbb{C})$ when $2 + 2\tau/N < \alpha \leq \beta < 2 + (2\theta + 4s)/N$ where $\int_{\mathbb{R}^{N}}u_{\lambda}^{2}dx = 2\rho$.
		\end{theorem}
		
		\begin{lemma}
			\label{lem4.22} Assume that $(f_{1})$, $(f_{2})$ hold, any positive ground state of (\ref{eq4.1}) is non-degenerate in $H_{rad}^{s}(\mathbb{R}^{N})$ and $\mathcal{K}_{\lambda}$ has only one element $u_{\lambda}$ for any $\lambda < 0$, then $\mathcal{K}$ is a $C^{1}$ curve in $\mathbb{R} \times H_{rad}^{s}(\mathbb{R}^{N})$ where $\mathcal{K} = \{(\lambda,u): u \in \mathcal{K}_{\lambda}\}$.
		\end{lemma}
		
		\textit{Proof.  } With the help of the result shown in the proof of Step 2 of Theorem \ref{thm2.1}, it is not difficult to verify that $u_{\lambda} \rightarrow u_{\lambda_{0}}$ in $H^{s}(\mathbb{R}^{N})$ as $\lambda \rightarrow \lambda_{0}$. In other words, $\mathcal{K}$ is a continuous curve in $\mathbb{R} \times H_{rad}^{s}(\mathbb{R}^{N})$. Then similar to \cite[Theorem 18]{SS}, we can prove that $\mathcal{K}$ is a $C^{1}$ curve in $\mathbb{R} \times H_{rad}^{s}(\mathbb{R}^{N})$ by the implicit function theorem since $u_{\lambda}$ is non-degenerate in $H_{rad}^{s}(\mathbb{R}^{N})$. \qed\vskip 5pt
		
		To study the orbitally stability, we lean on the following result, which expresses in our context the abstract theory developed in \cite{GSS}:
		
		\begin{proposition}
			\label{prop4.23} Under the hypotheses of Lemma \ref{lem4.22} and assume that $(LWP)$ holds. Then if $d'(\lambda) < 0$ (respectively $> 0$), the standing wave $e^{-i\lambda t}u_{\lambda}(x)$ is orbitally stable (respectively unstable), where $d(\lambda) = \frac{1}{2}\int_{\mathbb{R}^{N}}u_{\lambda}^{2}dx$ here.
		\end{proposition}
		
		\textbf{Proof of Theorem \ref{thm4.21}:  } Case 1: $\alpha > \frac{2(N + \tau) + 4s}{N}$.
		
		In this case, Theorem \ref{thm2.2} yields that $\lim_{\lambda \rightarrow -\infty}d(\lambda) = 0$ and $\lim_{\lambda \rightarrow 0^{-}}d(\lambda) = +\infty$. Thus $d^{-1}(\rho) \subset [e_{1},e_{2}]$ for any fixed $\rho \in (0,+\infty)$, where $-\infty < e_{1} < e_{2} < 0$ and $e_{1}$, $e_{2}$ depend on $\rho$. We assume that $\rho$ is a regular value. Then $d^{-1}(\rho)$ consist of finite points, i.e. $d^{-1}(\rho) = \{\lambda_{\rho}^{(1)}, \cdots, \lambda_{\rho}^{(k)}\}$ and $d'(\lambda_{\rho}^{(i)}) \neq 0$, $i = 1, 2, \cdots, k$. Without loss of generality, assume $\lambda_{\rho}^{(1)} < \cdots < \lambda_{\rho}^{(k)}$. Then $d(\lambda) < d(\lambda_{\rho}^{(1)})$ for any $\lambda < \lambda_{\rho}^{(1)}$ by the continuity of $d(\lambda)$, implying $d'(\lambda_{\rho}^{(1)}) \geq 0$. Considering $d'(\lambda_{\rho}^{(1)}) \neq 0$, we conclude that $d'(\lambda_{\rho}^{(1)}) > 0$. According to Proposition \ref{prop4.23}, the standing wave $e^{-i\lambda_{\rho}^{(1)} t}u_{\lambda_{\rho}^{(1)}}(x)$ is orbitally unstable. By Sard's theorem, we know such $\rho$ is almost everywhere in $(-\infty,0)$.
		
		Case 2: $2 + \frac{2\tau}{N} < \alpha \leq \beta < \frac{2(N + \theta) + 4s}{N}$.
		
		In this case, Theorem \ref{thm2.2} yields that $\lim_{\lambda \rightarrow -\infty}d(\lambda) = +\infty$ and $\lim_{\lambda \rightarrow 0^{-}}d(\lambda) = 0$. Then we can complete the proof after quiet similar arguments to case 1. \qed\vskip 5pt
		
		\begin{remark}
			\label{rmk4.24} When $f(|x|,u) = |u|^{q-2}u$, $2 < q < 2 + 4s/N$ or $2 + 4s/N < q < 2_{s}^{\ast}$, we can derive the stability or instability for all $\rho \in (-\infty,0)$ rather than almost everywhere $\rho \in (-\infty,0)$ in Theorem \ref{thm4.21} above.
		\end{remark}
		
		Considering the case when $f(|x|,u) = h(|x|)|u|^{q-2}u$, $2 < q < 2_{s}^{\ast}$, we have the following corollary:
		
		\begin{corollary}
			\label{cor4.25}  Under the hypotheses of Theorem \ref{cor4.19}, and $(LWP)$, $(\lambda, u_{\lambda})$ is given by Theorem \ref{cor4.19}. Then for a.e. $\rho \in (0, \infty)$, $e^{-i\lambda t}u_{\lambda}$ is orbitally unstable in $H_{rad}^{s}(\mathbb{R}^{N}, \mathbb{C})$ when $q > 2 + 4s/N$ while $e^{-i\lambda t}u_{\lambda}$ is orbitally stable in $H_{rad}^{s}(\mathbb{R}^{N}, \mathbb{C})$ when $2 < q < 2 + (2\theta + 4s)/N$ where $\int_{\mathbb{R}^{N}}u_{\lambda}^{2}dx = 2\rho$.
		\end{corollary}
		
		\subsection{Fractional nonlinear Schr\"{o}dinger equations with a potential I}
		
		In this subsection, we consider the fractional nonlinear Schr\"{o}dinger equation with a potential
		\begin{equation} \label{eq4.36}
			\left\{
			\begin{array}{cc}
				(-\Delta)^{s} u + V(|x|)u = \lambda u + h(|x|)|u|^{q-2}u \ in \ \mathbb{R}^{N},\\
				u(x) \rightarrow 0 \ as \ |x| \rightarrow +\infty,
			\end{array}
			\right.
		\end{equation}
		where $0 < s \leq 1$, $N \geq 2$, $2 < q < 2_{s}^{\ast}$, and $h(r)$ satisfies $(h_{1})$ introduced in Subsection 4.2. Recall that
		
		$(h_{1})$ $h(r) \in C^{1}[0,+\infty) \cap L^{\infty}[0,+\infty)$, $h(r) > 0$ in $[0,+\infty)$, $h(r)$ and $\frac{rh'(r)}{h(r)}$ are non-increasing in $(0,+\infty)$, $\theta = \lim_{r \rightarrow +\infty}\frac{rh'(r)}{h(r)} > -\infty$. \\
		For linear potential $V$, we assume
		
		$(V_{1})$ $V \in C^{1}([0,+\infty)) \cap L^{\infty}$, $V(0) = 0$, $V(r^{\ast}) > 0$ for some $r^{\ast} > 0$, $\gamma = \sup_{r > 0}\frac{rV'(r)}{V} < +\infty$, $V$ and $V + rV'(r)$ are non-decreasing on $[0,+\infty)$.
		
		When $(V_{1})$ holds, $(-\Delta)^{s} + V$ acting on $L^{2}_{rad}(\mathbb{R}^{N})$ is a self-adjoint operator with domain $H^{2s}_{rad}$ and form domain $H^{s}_{rad}$. Let $\lambda_{1} = \inf \sigma((-\Delta)^{s} + V)$. By min-max principle, $0 \leq \lambda_{1} \leq \lim_{r \rightarrow +\infty}V(r)$. Furthermore, $\lambda_{1} > 0$. Since $\sigma_{ess}((-\Delta)^{s} + V) = [\lim_{r \rightarrow +\infty}V(r),+\infty)$, $0$ is an eigenvalue if $\lambda_{1} = 0$. Thus there exists $\psi$ with $\|\psi\|_{L_{2}} = 1$ such that $\int_{\mathbb{R}^{N}}[|(-\Delta)^{\frac{s}{2}}\psi|^{2} + V\psi^{2}]dx = 0$, implying $\int_{\mathbb{R}^{N}}|(-\Delta)^{\frac{s}{2}}\psi|^{2}dx = 0$. We can derive $\psi \equiv 0$, contradicting $\|\psi\|_{L_{2}} = 1$.

		Similar to Lemma \ref{lem4.13b}, we have:
		
		\begin{lemma}
			\label{lem4.26} Let $u \in H_{rad}^{s}$ be a positive solution of $(\ref{eq4.36})$ and $\mu_{G}(u) = 1$. Suppose that $(h_{1})$ and $(V_{1})$ hold true, then $\ker D_{uu}\Phi_\lambda(u_{\lambda})|_{L^2_{rad}} = \{0\}$.
		\end{lemma}
		
		\textit{Proof.  }  Let $L_{+,V} = (-\Delta)^{s} + V - \lambda - h(|x|)|u|^{q-2}$. Direct calculation imply that:
		\begin{equation} \label{eq4.37}
			L_{+,V}(u) = (2-q)h(r)u^{q-1},
		\end{equation}
		\begin{equation} \label{eq4.38}
			L_{+,V}(x \cdot \nabla u) = (rh'(r) + 2sh(r))u^{q-1} - 2s(V(r) - \lambda + rV'(r)) u.
		\end{equation}
		Similarly to the proof of Lemma \ref{lem4.13b}, we just need to show that $-\frac{V(r) - \lambda + rV'(r)}{h(r)u(r)^{q-2}}$ is non-increasing in $(0,+\infty)$. When $\lambda \leq 0$, by $(V_{1})$, $V(r) - \lambda + rV'(r) \geq 0$ is non-decreasing. Since $h(r)u(r)^{q-2}$ is non-increasing, we derive that $-\frac{V(r) - \lambda + rV'(r)}{h(r)u(r)^{q-2}}$ is non-increasing. \qed\vskip 5pt
				
		\begin{theorem}
			\label{thm4.27}  Assume $(h_{1})$ and $(V_{1})$ hold with $(N-2s)q < 2(N+\theta)$. Then any positive ground state of $(\ref{eq4.36})$ is unique in $H^{s}$ for any fixed $\lambda \leq 0$.
		\end{theorem}
		
		\textit{Proof. } The proof is quiet similar to the one of Theorem \ref{thm4.9}, using abstract framework 1. Hence we omit the details. We further provide an alternative proof in Appendix \ref{Appen2}. \qed\vskip 5pt	
		
		Next we use the abstract framework 2 to study the existence of normalized solutions. In this case, $W = H_{rad}^{s}(\mathbb{R}^{N})$,
		$$
		\Phi_{\lambda}(u) = \frac{1}{2}\int_{\mathbb{R}^{N}}[|(-\Delta)^{\frac{s}{2}}u|^{2} + (V - \lambda) u^{2}]dx - \frac{1}{q}\int_{\mathbb{R}^{N}}h(|x|)|u|^{q}dx,
		$$
		$$
		S(u) = \frac{1}{2}\int_{\mathbb{R}^{N}}|(-\Delta)^{\frac{s}{2}}u|^{2}dx, G(u) = \frac{1}{2}\int_{\mathbb{R}^{N}}Vu^{2}dx,
		$$
		$$
		F(u) = \frac{1}{q}\int_{\mathbb{R}^{N}}h(|x|)|u|^{q}dx, Q(u) = \frac{1}{2}\int_{\mathbb{R}^{N}}|u|^{2}dx.
		$$
		Similar to Lemmas \ref{lem4.1} and \ref{lem4.2}, $(F_{1})$ and $(F_{2})$ hold when $(h_{1})$ and $(V_{1})$ hold. The following lemma shows $(F_{5})$ holds in mass subcritical case while $(F_{6})$ holds in mass supercritical case.
		
		After a similar discussion to the proof of Lemma \ref{lem4.3}, we can derive the Pohozaev identity for (\ref{eq4.36}).
		
		\begin{lemma}
			\label{lem4.28} (Pohozaev identity) Suppose $(h_{1})$ and $(V_{1})$ hold. Let $u \in H_{rad}^{s}$ be a weak solution of (\ref{eq4.36}), then we have the following integral identity
			\begin{eqnarray} \label{eq4.39}
			&& (N-2s)\int_{\mathbb{R}^{N}}|(-\Delta)^{\frac{s}{2}} u|^{2}dx + N\int_{\mathbb{R}^{N}}Vu^{2}dx + \int_{\mathbb{R}^{N}}|x|V'(|x|)u^{2}dx \nonumber \\
			&& =  N\lambda\int_{\mathbb{R}^{N}}u^{2}dx + \int_{\mathbb{R}^{N}}\frac{2Nh(|x|)+|x|h'(|x|)}{q}|u|^{q}dx.
			\end{eqnarray}
		\end{lemma}		
		
		\begin{lemma}
			\label{lem4.31}  Assume that $(h_{1})$ and $(V_{1})$ hold with $(N - 2s)\beta < 2(N + \theta)$. If $q > \frac{2N + 4s}{N}$, then $(F_{4})$ holds with $k > 1$. If $q < \frac{2(N + \theta) + 4s}{N}$, then $(F_{3})$ holds with $0 < l < 1$ for $\lambda \leq \lambda^{\ast} < -\frac{(2s + \gamma) \|V\|_{L^{\infty}}}{2s}\frac{2s (q-2) - 2\theta}{2(N+\theta) + 4s - qN}$.
		\end{lemma}
		
		\textit{Proof.  } By Lemma \ref{lem4.28} and
		\begin{equation}
		\int_{\mathbb{R}^{N}}(|(-\Delta)^{\frac{s}{2}}u|^{2} + Vu^{2})dx = \lambda\int_{\mathbb{R}^{N}}u^{2}dx + \int_{\mathbb{R}^{N}}h(|x|)|u|^{q})dx.
		\end{equation}		
		we have
		\begin{eqnarray} \label{eq4.52}
			\Phi_{\lambda}(u) &\geq& -\frac{\lambda}{2}[1 + \frac{qN-2N-4s}{2N-q(N-2s)}]\int_{\mathbb{R}^{N}}u^{2}dx \nonumber \\
			&+& \frac{1}{4s}[1 + \frac{qN-2N-4s}{2N-q(N-2s)}]\int_{\mathbb{R}^{N}}(2sV+rV')u^{2}dx \nonumber \\
			&\geq& -k\lambda Q(u),
		\end{eqnarray}
		where $k = 1 + \frac{qN-2N-4s}{2N-q(N-2s)} > 1$, if $q > \frac{2N + 4s}{N}$; and
		\begin{eqnarray} \label{eq4.53}
			\Phi_{\lambda}(u) &\leq& -\frac{\lambda}{2}[1 + \frac{qN-2(N+\theta)-4s}{2N-q(N-2s)}]\int_{\mathbb{R}^{N}}u^{2}dx  \nonumber \\
			&+& \frac{1}{4s}[1 + \frac{qN-2(N+\theta)-4s}{2N-q(N-2s)}]\int_{\mathbb{R}^{N}}(2sV+rV')u^{2}dx \nonumber \\
			&\leq& -\lambda[1 + \frac{qN-2(N+\theta)-4s}{2N-q(N-2s)}](1-\frac{(2s+\gamma)\|V\|_{L^{\infty}}}{2s\lambda}) Q(u) \nonumber \\
			&\leq& -l\lambda Q(u),
		\end{eqnarray}
		where $l = [1 + \frac{qN - 2(N+\theta) - 4s}{2N - q(N-2s)}](1 - \frac{(2s + \gamma) \|V\|_{L^{\infty}}}{2s\lambda^{\ast}}) \in (0,1)$, if $q < \frac{2(N + \theta) + 4s}{N}$ and $\lambda \leq \lambda^{\ast} < -\frac{(2s + \gamma) \|V\|_{L^{\infty}}}{2s}\frac{2s (q-2) - 2\theta}{2(N+\theta) + 4s - qN}$. \qed\vskip 5pt
		
		Theorem \ref{thm4.27} shows that any positive ground state of (\ref{eq4.36}) is unique in $H^{s}_{rad}$ for any fixed $\lambda \leq 0$. Then applying Corollary \ref{cor2.4} $(i)$ and $(ii)$, we have
		
		\begin{theorem}
			\label{thm4.32}  Assume $(h_{1})$ and $(V_{1})$ hold with $(N-2s)q < 2(N+\theta)$. Then there exists $\rho_{1} > 0$ such that for any $\rho \geq \rho_{1}$, (\ref{eq4.36}) admits at least one solution $(\lambda, u_{\lambda}) \in (-\infty, 0] \times H_{rad}^{s}$ such that $\int_{\mathbb{R}^{N}}u_{\lambda}^{2}dx = 2\rho$ if $2< q < \frac{2(N + \theta) + 4s}{N}$; there exists $\rho_{2} > 0$ such that for any $0 < \rho \leq \rho_{2}$, (\ref{eq4.36}) admits at least one solution $(\lambda, u_{\lambda}) \in (-\infty, 0] \times H_{rad}^{s}$ such that $\int_{\mathbb{R}^{N}}u_{\lambda}^{2}dx = 2\rho$ if $q > \frac{2N + 4s}{N}$.
		\end{theorem}
		
		\begin{remark}
			\label{rmk4.30} It is worth noting that results in Theorems \ref{thm4.27}, \ref{thm4.32} can be extended to the case $\lim_{r \rightarrow +\infty}V(r) = +\infty$. The method in Appendix \ref{Appen2} is invalid since we can not prove that $L_{+,\eta} \rightarrow L_{+,\tilde{\eta}}$ (the definition will be found in Appendix \ref{Appen2}) in the norm-resolvent sense as $\eta \rightarrow \tilde{\eta}$. However, the abstract framework 1 also works, and we can show the uniqueness of the ground state solution of (\ref{eq4.36}), and then use the abstract framework 3 to get the existence of normalized solutions. For unbounded potentials, the local uniqueness when $\lambda \rightarrow -\infty$ needs more discussions ( see some ideas used in \cite{Song4}).
		\end{remark}		
		
		Similar to the discussions in Subsection 4.3, if $(\lambda, u_{\lambda})$ is a solution of $(\ref{eq4.36})$, then $e^{-i\lambda t}u_{\lambda}$ is a standing wave solution of the following fractional nonlinear Schr\"{o}dinger equation
		
		\begin{equation}
			\left\{
			\begin{array}{cc}
				i\partial_{t} \psi - (-\Delta)^{s} \psi - V(x)\psi + h(|x|)|\psi|^{q-2}\psi = 0, (t,x) \in [0,+\infty) \times \mathbb{R}^{N}, \\
				\psi(0,x) = \psi_{0}(x),
			\end{array}
			\right.
		\end{equation}
		where $0 < s \leq 1$, $N \geq 2$, $2 < q < 2_{s}^{\ast}$.
		
		Similar arguments to the proof of Theorem \ref{thm4.21} yield
		
		\begin{theorem}
			\label{thm4.33}  Assume $(h_{1})$, $(V_{1})$ and $(LWP)$ hold with $(N-2s)q < 2(N+\theta)$. When $2< q < \frac{2(N + \theta) + 4s}{N}$, for a.e. $\rho \geq \rho_{1}$, $e^{-i\lambda t}u_{\lambda}$ is orbitally stable in $H_{rad}^{s}(\mathbb{R}^{N}, \mathbb{C})$ where $\int_{\mathbb{R}^{N}}u_{\lambda}^{2}dx = 2\rho$. When $q > \frac{2N + 4s}{N}$, for a.e. $\rho \in (0,\rho_{2})$, $e^{-i\lambda t}u_{\lambda}$ is orbitally unstable in $H_{rad}^{s}(\mathbb{R}^{N}, \mathbb{C})$ where $\int_{\mathbb{R}^{N}}u_{\lambda}^{2}dx = 2\rho$.
			$\rho_{1}$, $\rho_{2}$ and $(\lambda, u_{\lambda})$ are given by Theorem \ref{thm4.32}.
		\end{theorem}
		
		\begin{remark}
			\label{rmk4.34}  We can get more information if the non-degeneracy of positive ground state hold when $0 < \lambda < \lambda_{1}$, when $\lambda_{1}$ is the first eigenvalue of $(-\Delta)^{s} + V$. Indeed, we may derive that when $2< q < \frac{2(N + \theta) + 4s}{N}$, (\ref{eq4.36}) admits at least one solution $(\lambda, u_{\lambda}) \in (-\infty, \lambda_{1}) \times H_{rad}^{s}$ such that $\int_{\mathbb{R}^{N}}u_{\lambda}^{2}dx = 2\rho$ for any $\rho > 0$, and for a.e. $\rho > 0$, $e^{-i\lambda t}u_{\lambda}$ is orbitally stable in $H_{rad}^{s}(\mathbb{R}^{N}, \mathbb{C})$. When $q > \frac{2N + 4s}{N}$, there exists $\rho_{3} > 0$ such that \\
			- (\ref{eq4.36}) admits at least two solutions $(\lambda, u_{\lambda})$, $(\tilde{\lambda}, u_{\tilde{\lambda}}) \in (-\infty, \lambda_{1}) \times H_{rad}^{s}$ such that $\int_{\mathbb{R}^{N}}u_{\lambda}^{2}dx = \int_{\mathbb{R}^{N}}u_{\tilde{\lambda}}^{2}dx = 2\rho$ and $\lambda < \tilde{\lambda}$ for any $0 < \rho < \rho_{3}$, and for a.e. $0 < \rho \leq \rho_{3}$, $e^{-i\lambda t}u_{\lambda}$ is orbitally unstable while $e^{-i\tilde{\lambda} t}u_{\tilde{\lambda}}$ is orbitally stable in $H_{rad}^{s}(\mathbb{R}^{N}, \mathbb{C})$; \\
			- (\ref{eq4.36}) admits at least one solution $(\lambda, u_{\lambda}) \in (-\infty, \lambda_{1}) \times H_{rad}^{s}$ such that $\int_{\mathbb{R}^{N}}u_{\lambda}^{2}dx = 2\rho_{3}$; \\
			- (\ref{eq4.36}) admits no ground state $(\lambda, u_{\lambda}) \in (-\infty, \lambda_{1}) \times H_{rad}^{s}$ such that $\int_{\mathbb{R}^{N}}u_{\lambda}^{2}dx = 2\rho$ for any $\rho > \rho_{3}$.
		\end{remark}
		
		\begin{remark}
			\label{rmk4.35} Let $V(r) = 1 - 1/(1+r^{k})^{l}$, $k \geq 1$, $kl \geq 1$. Then $(V_{1})$ holds.
		\end{remark}
		
		\subsection{Nonlinear Schr\"{o}dinger equations with a potential II} \label{PII}
		
		$2sV + rV'(r)$ plays an important role when verifying $(F_{5})$ and $(F_{6})$, see Lemma \ref{lem4.31}. When $V = r^{-2s}$, $2sV + rV'(r) = 0$. For simplicity, we just study the case when $s = 1$ in this subsection and similar results can be extended to the case when $0 < s < 1$. Consider
		\begin{equation} \label{eq4.55}
			\left\{
			\begin{array}{cc}
				-\Delta u + \frac{1}{|x|^{2}}u = \lambda u + |u|^{q-2}u \ in \ \mathbb{R}^{N} \backslash \{0\},\\
				u(x) \rightarrow 0 \ as \ |x| \rightarrow +\infty,
			\end{array}
			\right.
		\end{equation}
		where $N \geq 3$, $2 < q < \frac{2N}{N-2}$. \cite[Theorem X.28]{RS2} yields that $-\Delta + \frac{1}{|x|^{2}}$ is essential self-adjoint on $C_{0}^{\infty}(\mathbb{R}^{N})$. Furthermore, the form domain is
		\[
		H^{1}_{V}:= \{u \in H^{1}(\mathbb{R}^{N}): \int_{\mathbb{R}^{N}}\frac{1}{|x|^{2}}u^{2}dx < +\infty\},
		\]
		with inner product
		\[
		\langle u,v \rangle_{H^{1}_{V}} = \int_{\mathbb{R}^{N}}[\nabla u \nabla v + (\frac{1}{|x|^{2}}+1)uv]dx
		\]
		and norm
		\[
		\|u\|_{H^{1}_{V}} = \sqrt{\int_{\mathbb{R}^{N}}[|\nabla u|^{2} + (\frac{1}{|x|^{2}}+1)u^{2}]dx}.
		\]
		
		Apply the abstract framework 2 where
		$$
		W = H^{1}_{V,rad} = \{u \in H^{1}_{V}: u(x) = u(|x|)\},
		$$
		$$
		\Phi_{\lambda}(u) = \frac{1}{2}\int_{\mathbb{R}^{N}}[|\nabla u|^{2} + (\frac{1}{|x|^{2}} - \lambda) u^{2}]dx - \frac{1}{q}\int_{\mathbb{R}^{N}}|u|^{q}dx,
		$$
		$$
		S(u) = \frac{1}{2}\int_{\mathbb{R}^{N}}|\nabla u|^{2}dx, G(u) = \frac{1}{2}\int_{\mathbb{R}^{N}}\frac{1}{|x|^{2}}u^{2}dx,
		$$
		$$
		F(u) = \frac{1}{q}\int_{\mathbb{R}^{N}}|u|^{q}dx, Q(u) = \frac{1}{2}\int_{\mathbb{R}^{N}}|u|^{2}dx.
		$$
		We don't know whether the positive ground state of (\ref{eq4.55}) is unique or not. However, we can show the existence of normalized solution. Similar to Lemmas \ref{lem4.1} and \ref{lem4.2}, $(F_{1})$ and $(F_{2})$ hold. Moreover, we will see that both $(F_{5})$ and $(F_{6})$ hold with $k = l$. To show this, we introduce the following Pohozaev identity:
		
		\begin{lemma}
			\label{lem4.36} (\cite[Proposition 2.3]{LLT}) Let $u \in H^{1}_{V}$ be a weak solution of (\ref{eq4.55}), then there holds the following integral identity
			\begin{equation} \label{eq4.56}
				(N-2)\int_{\mathbb{R}^{N}}(|\nabla u|^{2} + \frac{1}{|x|^{2}}u^{2})dx = N\lambda\int_{\mathbb{R}^{N}}u^{2}dx + \frac{2N}{q}\int_{\mathbb{R}^{N}}|u|^{q}dx.
			\end{equation}
		\end{lemma}
		
		\begin{lemma}
			\label{lem4.37}  Let $u \in H^{1}_{V}$ be a weak solution of (\ref{eq4.55}), then $\Phi_{\lambda}(u) = (1 + \frac{(q-2)N-4}{2N-(N-2)q})Q(u)$, i.e. $(F_{5})$ and $(F_{6})$ hold with $k = l = 1 + \frac{(q-2)N-4}{2N-(N-2)q}$.
		\end{lemma}
		
		\textit{Proof.  } By integrating (\ref{eq4.55}) with respect to $u$, we obtain
		\begin{equation} \label{eq4.57}
			\int_{\mathbb{R}^{N}}(|\nabla u|^{2} + \frac{1}{|x|^{2}}u^{2})dx = \lambda\int_{\mathbb{R}^{N}}u^{2}dx + \int_{\mathbb{R}^{N}}|u|^{q}dx.
		\end{equation}
		By (\ref{eq4.56}) and (\ref{eq4.57}), we get $\Phi_{\lambda}(u) = (1 + \frac{(q-2)N-4}{2N-(N-2)q})Q(u)$. \qed\vskip 5pt
		
		\begin{theorem}
			\label{thm4.38}  $(i)$ If $q \neq \frac{2N+4}{N}$, then for any $\rho > 0$, (\ref{eq4.55}) a solution $(\lambda,u_{\lambda}) \in (-\infty,0) \times H_{V,rad}^{1}$ with $\frac{1}{2}\int_{\mathbb{R}^{N}}u_{\lambda}^{2}dx = \rho$ and $\rho = \frac{C}{k}(-\lambda)^{\frac{1}{k}-1}$ for some $C > 0$, $k = 1 + \frac{(q-2)N-4}{2N-(N-2)q}$.
			
			$(ii)$ If $q = \frac{2N+4}{N}$, then there exists some $\rho^{\ast} > 0$ such that for any fixed $\lambda < 0$, (\ref{eq4.55}) admits a solution $u_{\lambda} \in H_{V,rad}^{1}$ with $\frac{1}{2}\int_{\mathbb{R}^{N}}u_{\lambda}^{2}dx = \rho^{\ast}$ (implying that (\ref{eq4.55}) admits infinite solutions $(\lambda,u_{\lambda}) \in (-\infty,0) \times H_{V,rad}^{1}$ with $\frac{1}{2}\int_{\mathbb{R}^{N}}u_{\lambda}^{2}dx = \rho^{\ast}$ when $\lambda$ changes) and (\ref{eq4.55}) admits no solution $u_{\lambda} \in H_{V,rad}^{1}$ with $\frac{1}{2}\int_{\mathbb{R}^{N}}u_{\lambda}^{2}dx = \rho$ for $\rho \neq \rho^{\ast}$.
		\end{theorem}
		
		\textit{Proof.  } $(i)$ and $(ii)$ in Theorem \ref{thm4.38} above are direct corollary of Theorem \ref{thm2.2} $(iv)$ and Corollary \ref{cor2.4} $(v)$. \qed\vskip 5pt
		
		\begin{remark}
			\label{rmk4.39} By Theorem \ref{thm2.2} $(iv)$, we also know that all the ground states of (\ref{eq4.55}) in $H_{V,rad}^{1}$ have the same $L^{2}$ norm. For this reason, we conjecture that the positive ground state of (\ref{eq4.55}) is unique in $H_{V,rad}^{1}$ for a fixed $\lambda < 0$.
		\end{remark}
		
		Similar to the discussions in Subsection 4.3, if $(\lambda, u_{\lambda})$ is a solution of $(\ref{eq4.55})$, then $e^{-i\lambda t}u_{\lambda}$ is a standing wave solution of the following nonlinear Schr\"{o}dinger equation
		
		\begin{equation} \label{eq4.58}
			\left\{
			\begin{array}{cc}
				i\partial_{t} \psi + \Delta \psi - \frac{1}{|x|^{2}}\psi + |\psi|^{q-2}\psi = 0, (t,x) \in [0,+\infty) \times \mathbb{R}^{N} \backslash \{0\}, \\
				\psi(0,x) = \psi_{0}(x),
			\end{array}
			\right.
		\end{equation}
		where $N \geq 3$, $2 < q < \frac{2N}{N-2}$.
		
		Since we don't know the non-degeneracy of ground state for (\ref{eq4.55}), we assume it in the following theorem. In this case, we can get the orbital stability or instability for all $\rho > 0$ rather than almost every $\rho > 0$.
		
		\begin{theorem}
			\label{thm4.40}  Assume $(LWP)$ holds and that any ground state of (\ref{eq4.55}) is non-degenerate in $H_{V,rad}^{1}$. When $2< q < \frac{2N+4}{N}$, for any $\rho > 0$, $e^{-i\lambda t}u_{\lambda}$ is orbitally stable in $H_{V,rad}^{1}(\mathbb{R}^{N}, \mathbb{C})$ where $\int_{\mathbb{R}^{N}}u_{\lambda}^{2}dx = 2\rho$. When $\frac{2N + 4}{N} < q < \frac{2N}{N-2}$, for any $\rho > 0$, $e^{-i\lambda t}u_{\lambda}$ is orbitally unstable in $H_{V,rad}^{1}(\mathbb{R}^{N}, \mathbb{C})$ where $\int_{\mathbb{R}^{N}}u_{\lambda}^{2}dx = 2\rho$.
			$(\lambda, u_{\lambda})$ is the solution given by Theorem \ref{thm4.38} and
			\[
			H_{V,rad}^{1}(\mathbb{R}^{N}, \mathbb{C}):= \{u \in H^{1}(\mathbb{R}^{N}, \mathbb{C}): u(x) = u(|x|), \int_{\mathbb{R}^{N}}\frac{1}{|x|^{2}}u\bar{u}dx < +\infty\}.
			\]
		\end{theorem}
		
		\textit{Proof.  } Let $(\lambda, u_{\lambda})$ be the solution given in Theorem \ref{thm4.38}. Then $Q(u_{\lambda}) = \frac{C}{k}(-\lambda)^{\frac{1}{k}-1}$. Thus $\partial_{\lambda}Q(u_{\lambda}) = -(\frac{1}{k}-1)\frac{C}{k}(-\lambda)^{\frac{1}{k}-2}$. When $2< q < \frac{2N+4}{N}$, $0 < k < 1$, implying $\partial_{\lambda}Q(u_{\lambda}) < 0$. When $\frac{2N + 4}{N} < q < \frac{2N}{N-2}$, $k > 1$, implying $\frac{\partial}{\partial \lambda}Q(u_{\lambda}) > 0$. Then the results in Theorem \ref{thm4.40} follow Proposition \ref{prop4.23}. \qed\vskip 5pt
		
		\section{\textbf{Fractional Laplacian with mixed nonlinearities}} \label{mn}
		
		\subsection{The existence of normalized solutions}
		
		Consider
		\begin{equation} \label{eq4.59}
			\left\{
			\begin{array}{cc}
				(-\Delta)^{s} u = \lambda u + |u|^{p-2}u + |u|^{q-2}u \ in \ \mathbb{R}^{N},\\
				u(x) \rightarrow 0 \ as \ |x| \rightarrow +\infty,
			\end{array}
			\right.
		\end{equation}
		where $N \geq 2$, $2 < p < 2 + 4s/N < q < 2_{s}^{\ast}$. Apply the abstract framework 3 where
		$$
		W = H^{s}_{rad}(\mathbb{R}^{N}) = \{u \in H^{s}(\mathbb{R}^{N}): u(x) = u(|x|)\},
		$$
		$$
		\Phi_{\lambda}(u) = \frac{1}{2}\int_{\mathbb{R}^{N}}[|(-\Delta)^{\frac{s}{2}} u|^{2} - \lambda u^{2}]dx - \frac{1}{p}\int_{\mathbb{R}^{N}}|u|^{p}dx - \frac{1}{q}\int_{\mathbb{R}^{N}}|u|^{q}dx,
		$$
		$$
		S(u) = \frac{1}{2}\int_{\mathbb{R}^{N}}|(-\Delta)^{\frac{s}{2}} u|^{2}dx, G(u) = 0,
		$$
		$$
		F(u) = \frac{1}{p}\int_{\mathbb{R}^{N}}|u|^{p}dx + \frac{1}{q}\int_{\mathbb{R}^{N}}|u|^{q}dx, Q(u) = \frac{1}{2}\int_{\mathbb{R}^{N}}|u|^{2}dx.
		$$
		Note that $(S_{1})$, $(S_{2})$, $(Q_{1})$, $(N_{1})$ and $(N_{1})$ hold with $\gamma = N/2s$ and
		$$
		T(\lambda)u(x) = u(\frac{x}{|\lambda|^{1/2s}}).
		$$
		
		When $\lambda \rightarrow 0^{-}$,
		\begin{eqnarray}
		\widehat{\Phi}_{\lambda}(w) &=& \frac{1}{2}\int_{\mathbb{R}^{N}}|(-\Delta)^{\frac{s}{2}} w|^{2}dx\nonumber \\
		&+& \frac{1}{2}\int_{\mathbb{R}^{N}}|w|^{2}dx - \frac{1}{p}\int_{\mathbb{R}^{N}}|w|^{p}dx - \frac{1}{q}|\lambda|^{\frac{q-p}{p-2}}\int_{\mathbb{R}^{N}}|w|^{q}dx,
		\end{eqnarray}
		\begin{eqnarray}
		\widehat{\Phi}_{0}(w) = \frac{1}{2}\int_{\mathbb{R}^{N}}|(-\Delta)^{\frac{s}{2}} w|^{2}dx + \frac{1}{2}\int_{\mathbb{R}^{N}}|w|^{2}dx - \frac{1}{p}\int_{\mathbb{R}^{N}}|w|^{p}dx.
		\end{eqnarray}
		We have
		\begin{lemma}
		\label{lemE.1}	$(A_{1})$ - $(A_{5})$ hold.
		\end{lemma}
		
		\textit{Proof.  } $(A_{1})$ and $(A_{2})$ can be derived in the same way as in the proofs of Lemmas \ref{lem4.1} and \ref{lem4.2}.  $(A_{3})$ is standard when $s = 1$ and has been shown in \cite[Theorems 3, 4]{FLS} when $0 < s < 1$. $(A_{4})$ can be proved by a standard process, c.f. \cite[Lemma 1.20]{Wi}. Finally, we prove $(A_{5})$. For $\lambda \in (-1,0)$ and $r > 0$ small enough, we have
		$$
		\inf_{\|w\| = r}\widehat{\Phi}_{\lambda}(w) \geq \frac{1}{2}r^{2} - C(r^{p} + r^{q}) > 0.
		$$
		Therefore,
		$$
		\liminf_{\lambda \rightarrow 0^{-}}\widehat{h}(\lambda) \geq \liminf_{\lambda \rightarrow 0^{-}}\inf_{\|w\| = r}\widehat{\Phi}_{\lambda}(w) > 0.
		$$
		\qed\vskip 5pt
		
		When $\mu = 1/\lambda \rightarrow 0^{-}$,
		\begin{eqnarray}
		\widetilde{\Phi}_{\mu}(v) &=& \frac{1}{2}\int_{\mathbb{R}^{N}}|(-\Delta)^{\frac{s}{2}} v|^{2}dx \nonumber \\
		&+& \frac{1}{2}\int_{\mathbb{R}^{N}}|v|^{2}dx - \frac{1}{p}|\mu|^{\frac{q-p}{q-2}}\int_{\mathbb{R}^{N}}|v|^{p}dx - \frac{1}{q}\int_{\mathbb{R}^{N}}|v|^{q}dx,
		\end{eqnarray}
		\begin{eqnarray}
		\widetilde{\Phi}_{0}(v) = \frac{1}{2}\int_{\mathbb{R}^{N}}|(-\Delta)^{\frac{s}{2}} v|^{2}dx + \frac{1}{2}\int_{\mathbb{R}^{N}}|v|^{2}dx - \frac{1}{p}\int_{\mathbb{R}^{N}}|v|^{q}dx.
		\end{eqnarray}
		We have
	    \begin{lemma}
		    $(B_{1})$ - $(B_{5})$ hold.
	    \end{lemma}
	
	    \textit{Proof.  } We can complete the proof after similar arguments as in the one of Lemma \ref{lemE.1}. \qed\vskip 5pt
		
		\begin{theorem}
		\label{thmE.3}	There exist $\Lambda _{2} < \Lambda_{1} < 0$ such that when $\lambda < \Lambda _{2}$ or $\Lambda_{1} < \lambda < 0$, the positive ground state of (\ref{eq4.59}) is unique (up to a translation) and non-degenerate in $H^{s}(\mathbb{R}^{N})$.
		\end{theorem}
		
		\textit{Proof.  } Theorems \ref{thmc.1} and \ref{thmc.4} yield that there exist $\Lambda _{2} < \Lambda_{1} < 0$ such that when $\lambda < \Lambda _{2}$ or $\Lambda_{1} < \lambda < 0$, the positive ground state of (\ref{eq4.59}) is unique and non-degenerate in $H^{s}_{rad}(\mathbb{R}^{N})$.
		
		Similar to the proof of \cite[Proposition 3.1]{FLS}, applying the moving plane method in \cite{MZ} for
		nonlocal equations, we can show that up to a translation, any positive ground state of (\ref{eq4.59}) in $H^{s}(\mathbb{R}^{N})$ is radial and strictly decreasing in $r = |x|$. Therefore, we obtain the results in this theorem. \qed\vskip 5pt
		
		\begin{remark}
			Using the method developed in \cite{FL, FLS}, for any $\lambda < 0$, we can obtain the uniqueness of positive ground state for
			$$
			(-\Delta)^{s}u = \lambda u + |u|^{p-2}u, 0 < s < 1,
			$$
			However, the problem is still open for general nonlinearities. In this section, we provide an extension to this result.
		\end{remark}
		
		\begin{theorem}
		\label{thmE.5}	There exists $\widehat{c} > 0$ such that for any $c < \widehat{c}$, (\ref{eq4.59}) admits two ground states $u_{\lambda}$, $u_{\widetilde{\lambda}}$ with $\lambda < \widetilde{\lambda} < 0$ and $\frac{1}{2}\int_{\mathbb{R}^{N}}|u_{\lambda}|^{2}dx = \frac{1}{2}\int_{\mathbb{R}^{N}}|u_{\widetilde{\lambda}}|^{2}dx = c$.
		\end{theorem}
		
		\textit{Proof.  } Note that $\gamma = N/2s$ and $2 < p < 2 + 4s/N < q < 2_{s}^{\ast}$. Therefore, this theorem is a direct corollary of the combination of Corollary \ref{corc.3} $(ii)$ and Corollary \ref{corc.6} $(i)$ with $\widehat{c} = \min\{c_{2},c_{3}\}$. \qed\vskip 5pt
		
		\begin{remark}
			After a slight modification, our method can handle the case of non-autonomous mixed nonlinearities, i.e.
			\begin{equation}
			\left\{
			\begin{array}{cc}
			(-\Delta)^{s} u = \lambda u + a(x)|u|^{p-2}u + b(x)|u|^{q-2}u \ in \ \mathbb{R}^{N},\\
			u(x) \rightarrow 0 \ as \ |x| \rightarrow +\infty,
			\end{array}
			\right.
			\end{equation}
			where $2 < p < q < 2_{s}^{\ast}$, $a(x) = a(|x|) \in C^{1}$, $\lim_{r \rightarrow +\infty}a(r) = a_{\infty} > 0$, $\inf_{r > 0}a(r) > 0$ and $b(x) = b(|x|) \in L^{\infty} \cap C^{1}$, $\inf_{r > 0}b(r) > 0$ (see details in Section 7). Furthermore, if $a(x)$ is a constant and $b(x) \in L^{\infty}$, we can address the case when $\lambda \rightarrow 0^{-}$ and if $b(x)$ is a constant and $a(x) \in L^{\infty}$, we can handle the case when $\lambda \rightarrow -\infty$. It is worth noting that we do not even need to consider whether the (PS) condition holds if we do not consider the uniqueness of the ground state but only obtain the existence of the normalized solutions.
	
		    In autonomous cases, we may obtain normalized solutions along the lines of \cite{S1}. However, the method in \cite{S1} cannot be extended to non-autonomous cases directly.
		\end{remark}
		
		\begin{remark}
			We can get more information about the normalized solutions of (\ref{eq4.59}) if its positive ground state is unique for any $\lambda < 0$. In fact, once $d(\lambda) = \hat{d}(\lambda)$ for $\lambda < 0$, there exists $\tilde{c} > 0$ such that (\ref{eq4.59}) admits at least two ground states with $L^{2}$ norm $c \in (0,\tilde{c})$,  at least one ground state with $L^{2}$ norm $\tilde{c}$, and no ground state with $L^{2}$ norm $c > \tilde{c}$. Due to the length limitation of this paper, we will address this problem in the future work.
		\end{remark}
		
		\subsection{The orbital stability or instability}
		
		If $(\lambda, u_{\lambda})$ is a solution of $(\ref{eq4.59})$, then $e^{-i\lambda t}u_{\lambda}$ is a standing wave solution of the fractional nonlinear Schr\"{o}dinger equation (fNLS)
		
		\begin{equation}
		\left\{
		\begin{array}{cc}
		i\partial_{t} \psi - (-\Delta)^{s} \psi + |\psi|^{p-2}\psi + |\psi|^{q-2}\psi = 0, (t,x) \in [0,+\infty) \times \mathbb{R}^{N}, \\
		\psi(0,x) = \psi_{0}(x),
		\end{array}
		\right.
		\end{equation}
	    where $N \geq 2$, $2 < p < 2 + 4s/N < q < 2_{s}^{\ast}$. In this subsection, we study the orbital stability or instability of the ground state solutions for (\ref{eq4.59}) which are shown in Theorem \ref{thmE.5}.
	
	    When $\Lambda_{1} < \lambda < 0$ where $\Lambda_{1}$ is given by Theorem \ref{thmE.3}, recall that
	    $$
	    w_{\lambda}(x) = |\lambda|^{-\frac{1}{p-2}}u_{\lambda}(\frac{x}{|\lambda|^{1/2s}}) \rightarrow w^{\ast}(x) \ as \ \lambda \rightarrow 0^{-}
	    $$
	    where $w^{\ast}(x)$ is the unique positive ground state in $H^{s}_{rad}(\mathbb{R}^{N})$ of
	    $$
	    (-\Delta)^{s}w + w = |w|^{p-2}w.
	    $$
	    Let $\chi_{\lambda}(x) = \partial_{\lambda}u_{\lambda}$ and
	    $$
	    \widehat{\chi}_{\lambda}(x) = |\lambda|^{\frac{p-3}{p-2}}\chi_{\lambda}(\frac{x}{|\lambda|^{1/2s}}).
	    $$
	    We have
	
	    \begin{lemma}
	    \label{lemE.8}	$(i)$ $\widehat{\chi}_{\lambda}(x)       \rightarrow \frac{1}{2-p}w^{\ast} - \frac{1}{2s}x \cdotp \nabla w^{\ast}$ in $H^{s}(\mathbb{R}^{N})$ as $\lambda \rightarrow 0^{-}$.
	    	
	    	$(ii)$ There exists $\widetilde{\Lambda}_{1} \in [\Lambda_{1},0)$ such that for any $\lambda \in (\widetilde{\Lambda}_{1},0)$,
	    	$$
	    	\partial_{\lambda}\int_{\mathbb{R}^{N}}|u_{\lambda}|^{2}dx < 0.
	    	$$
	    \end{lemma}
	
	    \textit{Proof.  } $(i)$ can be proved along the lines of the proof of \cite[Proposition 11 $(ii)$]{FSK}. Then
	    \begin{eqnarray}
	    && \lim_{\lambda \rightarrow 0^{-}}|\lambda|^{\frac{p-4}{p-2}+\frac{N}{2s}}\partial_{\lambda}\int_{\mathbb{R}^{N}}|u_{\lambda}|^{2}dx \nonumber \\
	    &=& \lim_{\lambda \rightarrow 0^{-}}2|\lambda|^{\frac{p-4}{p-2}+\frac{N}{2s}}\int_{\mathbb{R}^{N}}u_{\lambda}\chi_{\lambda}dx \nonumber \\
	    &=& \lim_{\lambda \rightarrow 0^{-}}2\int_{\mathbb{R}^{N}}w_{\lambda}\widehat{\chi}_{\lambda}dx \nonumber \\
	    &=& 2\int_{\mathbb{R}^{N}}w^{\ast}(\frac{1}{2-p}w^{\ast} - \frac{1}{2s}x \cdotp \nabla w^{\ast})dx \nonumber \\
	    &=& 2(\frac{1}{2-p}+\frac{N}{4s})\int_{\mathbb{R}^{N}}|w^{\ast}|^{2}dx < 0.
	    \end{eqnarray}
	    Therefore, there exists $\widetilde{\Lambda}_{1} \in [\Lambda_{1},0)$ such that for any $\lambda \in (\widetilde{\Lambda}_{1},0)$,
	    $$
	    \partial_{\lambda}\int_{\mathbb{R}^{N}}|u_{\lambda}|^{2}dx < 0.
	    $$
	    \qed\vskip 5pt
	
	    When $\lambda > \Lambda_{2},$ where $\Lambda_{2}$ is given by Theorem \ref{thmE.3}, set $\chi_{\lambda}(x) = \partial_{\lambda}u_{\lambda}$ and
	    $$
	    \widetilde{\chi}_{\lambda}(x) = |\lambda|^{\frac{q-3}{q-2}}\chi_{\lambda}(\frac{x}{|\lambda|^{1/2s}}).
	    $$
	    Recall that $v^{\ast}(x)$ is the unique positive ground state in $H^{s}_{rad}(\mathbb{R}^{N})$ of
	    $$
	    (-\Delta)^{s}v + v = |v|^{q-2}v.
	    $$
	    Similar to the proof of Lemma \ref{lemE.8}, we have
	
	    \begin{lemma}
	    \label{lemE.9}	$(i)$ $\widetilde{\chi}_{\lambda}(x) \rightarrow \frac{1}{2-p}v^{\ast} - \frac{1}{2s}x \cdotp \nabla v^{\ast}$ in $H^{s}(\mathbb{R}^{N})$ as $\lambda \rightarrow -\infty$.
	    	
	    	$(ii)$ There exists $\widetilde{\Lambda}_{2} \leq \Lambda_{2}$ such that for any $\lambda \leq \widetilde{\Lambda}_{2}$,
	    	$$
	    	\partial_{\lambda}\int_{\mathbb{R}^{N}}|u_{\lambda}|^{2}dx > 0.
	    	$$
	    \end{lemma}
	
	    \begin{theorem}
	    	Assume that (LWP) holds. Then $e^{-i\lambda t}u_{\lambda}$ is orbitally stable in $H^{s}(\mathbb{R}^{N}, \mathbb{C})$ when $\lambda \in (\widetilde{\Lambda}_{1},0)$ while $e^{-i\lambda t}u_{\lambda}$ is orbitally unstable in $H^{s}(\mathbb{R}^{N}, \mathbb{C})$ when $\lambda < \widetilde{\Lambda}_{2}$ where $\widetilde{\Lambda}_{1}$, $\widetilde{\Lambda}_{2}$ are given by Lemmas \ref{lemE.8}, \ref{lemE.9} respectively and $u_{\lambda}$ is the unique and non-degenerate positive ground state of (\ref{eq4.59}) given by Theorem \ref{thmE.3}.
	    \end{theorem}
	
	    \textit{Proof.  } This theorem is a corollary of Proposition \ref{prop4.23} and Lemmas \ref{lemE.8} $(ii)$, \ref{lemE.9} $(ii)$. \qed\vskip 5pt
	
	    \begin{remark}
	    	For $c > 0$ small enough, we know that (\ref{eq4.59}) admits at least two positive ground states with $L^{2}$ norm $c$, one is orbitally stable while another is orbitally unstable.
	    \end{remark}

		\section{\textbf{p-laplace equations $-\Delta_{p} u = \lambda u + f(x,u)$}} \label{p}
		
		Our next application is p-laplace equations. We use abstract framework 2 in this section. In subsection 6.1, we study the case on the entire space. In subsections 6.2 and 6.3, we study the case on smooth bound domains.
		
		\subsection{p-laplace equations on entire space}
		
		Consider
		\begin{equation} \label{eq5.1}
			\left\{
			\begin{array}{cc}
				-\Delta_{p} u = \lambda u^{p-1} + f(|x|,u) \ in \ \mathbb{R}^{N},\\
				u(x) \rightarrow 0 \ as \ |x| \rightarrow +\infty,
			\end{array}
			\right.
		\end{equation}
		where $\Delta_{p} u = div(|\nabla u|^{p-2}\nabla u)$, $p > 1$. We may assume $f(r,t) = 0$ whenever $t < 0$ since we aim to show the existence of positive solutions. Furthermore, we assume on the nonlinear term $f(r,t) \in C^{1}([0, +\infty) \times [0, +\infty), \mathbb{R})$:
		
		$(f_{1}')$ $f(r,0) = 0$ and there exist $p < \alpha \leq \beta < 2_{p}^{\ast}$ such that
		\begin{equation}
			0 < (\alpha - 1) f(r,t) \leq f_{t}(r,t)t \leq (\beta - 1) f(r,t), \forall t > 0, r > 0,
		\end{equation}
		where $2_{p}^{\ast} = \frac{pN}{N-p}$ if $N > p$ and $2_{p}^{\ast} = +\infty$ if $N \leq p$.
		
		$(f_{2})$ $f(r,1)$ is bounded, i.e. there exists $M > 0$ such that $\vert f(r,1) \vert \leq M$, $\forall r > 0$.
		
		$(f_{3}')$ $F_{r}(r,t)r \leq C (|t|^{q_{1}} + |t|^{q_{2}})$ for some $C > 0$, $q_{1}$, $q_{2} \in (p,2_{p}^{\ast})$ and there exists $\theta \leq 0$ such that
		\begin{equation}
			F_{r}(r,t)r \geq \theta F(r,t), \forall t > 0, r > 0,
			\nonumber
		\end{equation}
		where $F(r,t) = \int_{0}^{t}f(r,s)ds$.
		
		$(f_{4}')$ $F_{r}(r,t)r \geq -C (|t|^{q_{1}} + |t|^{q_{2}})$ for some $C > 0$, $q_{1}$, $q_{2} \in (p,2_{p}^{\ast})$ and there exists $\tau \geq 0$ such that
		\begin{equation}
			F_{r}(r,t)r \leq \tau F(r,t), \forall t > 0, r > 0.
			\nonumber
		\end{equation}
		
		In this case, $W = W_{rad}^{1,p}(\mathbb{R}^{N})$, the radially symmetric Sobolev space, equipped with the norm
		$$
		\|u\| = (\int_{\mathbb{R}^{N}}(|\nabla u|^{p} + |u|^{p})dx)^{\frac{1}{p}},
		$$
		$$
		\Phi_{\lambda}(u) = \frac{1}{p}\int_{\mathbb{R}^{N}}(|\nabla u|^{p} - \lambda u^{p})dx - \int_{\mathbb{R}^{N}}F(|x|,u)dx,
		$$
		$$
		S(u) = \frac{1}{p}\int_{\mathbb{R}^{N}}|\nabla u|^{p}dx, G(u) = 0,
		$$
		$$
		F(u) = \int_{\mathbb{R}^{N}}F(|x|,u)dx, Q(u) = \frac{1}{p}\int_{\mathbb{R}^{N}}|u|^{p}dx,
		$$
		$$
		\mathcal{N}_{\lambda} = \{u \in W_{rad}^{1,p}(\mathbb{R}^{N}): \int_{\mathbb{R}^{N}}(|\nabla u|^{p} - \lambda u^{p})dx - \int_{\mathbb{R}^{N}}f(|x|,u)udx = 0\}.
		$$
		
		Similar to Lemmas \ref{lem4.1} and \ref{lem4.2}, we have
		
		\begin{lemma}
			\label{lem5.1} Assume $(f_{1}')$ and $(f_{2})$. Then $(F_{1})$ and $(F_{2})$ hold.
		\end{lemma}
		
		A similar argument to \cite[Theorem B.3]{Wi} yields
		
		\begin{lemma}
			\label{lem5.2} (Pohozaev identity) Suppose $(f_{1}')$, $(f_{2})$, $(f_{3}')$ or $(f_{1}')$, $(f_{2})$, $(f_{4}')$. Let $u \in W_{rad}^{1,p}(\mathbb{R}^{N})$ be a weak solution of (\ref{eq5.1}), then there holds the following integral identity
			\begin{eqnarray} \label{eq5.2}
				&& (N-p)\int_{\mathbb{R}^{N}}|\nabla u|^{p}dx \nonumber \\
				&=& N\lambda\int_{\mathbb{R}^{N}}u^{p}dx + pN\int_{\mathbb{R}^{N}}F(|x|,u)dx + p\int_{\mathbb{R}^{N}}|x|F_{r}(|x|,u)dx.
			\end{eqnarray}
		\end{lemma}
		
		Quiet similar to Lemma \ref{lem4.4}, we have
		
		\begin{lemma}
			\label{lem5.3}  Assume that $(f_{1}')$, $(f_{2})$, $(f_{3}')$ and $(f_{4})$ hold with $(N - p)\beta < p(N + \theta)$. If $\alpha > \frac{p(N + \tau) + p^{2}}{N}$, then $(F_{5})$ and $(F_{6})$ hold with $l \geq k > 1$. If $p + \frac{p\tau}{N} < \alpha \leq \beta < \frac{p(N + \theta) + p^{2}}{N}$, then $(F_{5})$ and $(F_{6})$ hold with $0 < k \leq l < 1$.
		\end{lemma}
		
		As an application of Corollary \ref{cor2.4}, we derive
		
		\begin{theorem}
			\label{thm5.4} Assume that $(f_{1}')$, $(f_{2})$, $(f_{3}')$ and $(f_{4})$ hold with $(N - p)\beta < p(N + \theta)$. Furthermore, we also assume that (\ref{eq5.1}) admits at most one positive ground state in $W_{rad}^{1,p}(\mathbb{R}^{N})$ for any fixed $\lambda < 0$. Then for any $\rho > 0$, (\ref{eq5.1}) admits at least one solution $(\lambda, u_{\lambda}) \in (-\infty, 0) \times W_{rad}^{1,p}(\mathbb{R}^{N})$ such that $\int_{\mathbb{R}^{N}}u_{\lambda}^{p}dx = p\rho$ if $\alpha > \frac{p(N + \tau) + p^{2}}{N}$ or $p + \frac{p\tau}{N} < \alpha \leq \beta < \frac{p(N + \theta) + p^{2}}{N}$.
		\end{theorem}
		
		\begin{remark}
			\label{rmk5.5}  When $f$ is autonomous and satisfies some conditions, (\ref{eq5.1}) admits at most one positive solution in $W_{rad}^{1,p}(\mathbb{R}^{N})$, see, e.g., \cite{FLSe, PS3, ST, Tang}. And it is a future work to provide non-autonomous cases.
		\end{remark}
		
		\subsection{p-laplace equations in a ball}
		
		Consider
		\begin{equation} \label{eq5.4}
			\left\{
			\begin{array}{cc}
				-\Delta_{p} u = \lambda u + f(|x|,u), x \in B_{R}, \\
				u_{|\partial B_{R}} = 0,
			\end{array}
			\right.
		\end{equation}
		where $B_{R} = \{x \in \mathbb{R}^{N}: |x| < R\}$, $N \geq 2$, $\Delta_{p} u = div(|\nabla u|^{p-2}\nabla u)$, $p > 1$. We may assume $f(r,t) = -f(r,-t)$ whenever $t < 0$ since we aim to show the existence of positive solutions. Furthermore, we assume on the nonlinear term $f(r,t) \in C^{1}([0, R) \times [0, +\infty), \mathbb{R})$:
		
		$(f_{5})$ $f(r,0) = f_{t}(r,0) = 0$ and there exist $\frac{pN + p^{2}}{N} < \alpha \leq \beta < 2_{p}^{\ast}$
		such that
		\begin{equation} \label{eq5.5}
			0 < (\alpha - 1) f(r,t) \leq f_{t}(r,t)t \leq (\beta - 1)f(r,t), \forall t > 0, r \in (0,R).
		\end{equation}
		
		$(f_{6})$ $f(r,1)$ is bounded, i.e. there exists $M > 0$ such that $\vert f(r,1) \vert \leq M$, $\forall r \in (0,R)$.
		
		$(f_{7})$ $F_{r}(r,t)r \geq -C (|t|^{q_{1}} + |t|^{q_{2}})$ for some $C > 0$, $q_{1}$, $q_{2} \in (p,2_{p}^{\ast})$ and there exists $\tau \geq 0$ such that
		\begin{equation}
			F_{r}(r,t)r \leq \tau F(x,t), \forall t \in \mathbb{R}, r \in (0,R).
			\nonumber
		\end{equation}
		
		In this case, $W = W_{0,rad}^{1,p}(B_{R})$, the radially symmetric Sobolev space, equipped with the norm
		$$
		\|u\| = (\int_{B_{R}}(|\nabla u|^{p} + |u|^{p})dx)^{\frac{1}{p}},
		$$
		$$
		\Phi_{\lambda}(u) = \frac{1}{p}\int_{B_{R}}(|\nabla u|^{p} - \lambda u^{p})dx - \int_{B_{R}}F(|x|,u)dx,
		$$
		$$
		S(u) = \frac{1}{p}\int_{B_{R}}|\nabla u|^{p}dx, G(u) = 0,
		$$
		$$
		F(u) = \int_{B_{R}}F(|x|,u)dx, Q(u) = \frac{1}{p}\int_{B_{R}}|u|^{p}dx,
		$$
		$$
		\mathcal{N}_{\lambda} = \{u \in W_{0,rad}^{1,p}(B_{R}): \int_{B_{R}}(|\nabla u|^{p} - \lambda u^{p})dx - \int_{B_{R}}f(|x|,u)udx = 0\}.
		$$
		
		Similar to Lemmas \ref{lem4.1} and \ref{lem4.2}, we have
		
		\begin{lemma}
			\label{lem5.6} Assume $(f_{5})$ and $(f_{6})$. Then $(F_{1})$ and $(F_{2})$ hold.
		\end{lemma}
		
		\begin{lemma}
			\label{lem5.7} (Pohozaev identity) Suppose $(f_{5})$ - $(f_{7})$. Let $u \in W_{0,rad}^{1,p}(B_{R})$ be a weak solution of (\ref{eq5.4}), then there holds the following integral identity
			\begin{eqnarray} \label{eq5.6}
				&& (N-p)\int_{B_{R}}|\nabla u|^{p}dx + (p - 1)\int_{\partial B_{R}}|\nabla u|^{p} \sigma \cdot \nu d\sigma \nonumber \\
				&=& N\lambda\int_{B_{R}}u^{p}dx + pN\int_{B_{R}}F(|x|,u)dx + p\int_{B_{R}}|x|F_{r}(|x|,u)dx.
			\end{eqnarray}
		\end{lemma}
		
		The autonomous version of Lemma \ref{lem5.7} above can be found in \cite[Proposition 1]{PS2} or \cite[Theorem 1.3]{DI} and the proof of non-autonomous version is similar.
		
		\begin{lemma}
			\label{lem5.8}  Assume that $(f_{5})$ - $(f_{7})$ hold with $(N - p)\alpha < p(N + \tau)$ and $\alpha > \frac{p(N + \tau) + p^{2}}{N}$, then $(F_{6})$ holds with $k > 1$.
		\end{lemma}
		
		\textit{Proof.  } Since $\sigma \cdot \nu = \sigma \cdot \frac{\sigma}{\vert \sigma \vert} = \vert \sigma \vert = R > 0$, (\ref{eq5.6}) yields that
		\begin{eqnarray} \label{eq5.7}
			&& (N-p)\int_{B_{R}}|\nabla u|^{p}dx \nonumber \\
&\leq& N\lambda\int_{B_{R}}u^{p}dx + pN\int_{B_{R}}F(|x|,u)dx + p\int_{B_{R}}|x|F_{r}(|x|,u)dx.
		\end{eqnarray}
		Then along the lines of the proof of Lemma \ref{lem4.4}, we can show that $(F_{6})$ holds with
		\[
		k = \frac{N\alpha - p(N + \tau) - p^{2}}{p(N + \tau) - \alpha(N - p)} + 1 > 1.
		\]
		\qed\vskip 5pt
		
		\begin{theorem}
			\label{thm5.9} Assume that $(f_{5})$ - $(f_{7})$ hold with $(N - p)\alpha < p(N + \tau)$ and $\alpha > \frac{p(N + \tau) + p^{2}}{N}$. Furthermore, we also assume that (\ref{eq5.4}) admits at most one positive solution in $W_{0,rad}^{1,p}(B_{R})$ for any fixed $\lambda < \lambda_{1}$. Then there exists some $b > 0$ such that (\ref{eq5.4}) admits at least one solution $(\lambda, u_{\lambda}) \in (-\infty, \lambda_{1}) \times W_{0,rad}^{1,p}(B_{R})$ where $\lambda < \lambda_{1}$, $u_{\lambda} > 0$ and $\int_{B_{R}}u_{\lambda}^{p}dx = pb$; (\ref{eq5.4}) admits at least two solutions $(\lambda, u_{\lambda}), (\tilde{\lambda}, u_{\tilde{\lambda}}) \in (-\infty, \lambda_{1}) \times W_{0,rad}^{1,p}(B_{R})$ where $\lambda < \tilde{\lambda} < \lambda_{1}$, $u_{\lambda} > 0$, $u_{\tilde{\lambda}} > 0$ and $\int_{B_{R}}u_{\lambda}^{p}dx = \int_{B_{R}}\tilde{u}_{\lambda}^{p}dx = pc$ for any $0 < c < b$; and (\ref{eq5.4}) admits no $(\lambda, u_{\lambda}) \in (-\infty, \lambda_{1}) \times W_{0,rad}^{1,p}(B_{R})$ where $u_{\lambda} > 0$ and $\int_{B_{R}}u_{\lambda}^{p}dx = pc$ for any $c > b$.
		\end{theorem}
		
		\textit{Proof.  } The uniqueness of the positive solution in $W_{0,rad}^{1,p}(B_{R})$ derives that $d \equiv \hat{d}$. By Theorem \ref{thm2.2} $(ii)$, $\lim_{\lambda \rightarrow -\infty}\hat{d}(\lambda) = 0$. Next, we study the behavior of $d(\lambda)$ near $\lambda_{1}$ with the help of a bifurcation result.
		
		By \cite[Theorem 4]{Dr} or \cite[Theorem 1.1]{PM}, $(\lambda_{1},0)$ is a bifurcation point and for $\epsilon > 0$ small enough, we can take $(\lambda(s),u(s)) \in \mathbb{R} \times W_{0,rad}^{1,p}(B_{R})$, $|s| < \epsilon$ such that $\lambda(s) \rightarrow \lambda_{1}$, $u(s) \rightarrow 0$ and $(\lambda(s),u(s))$ solves (\ref{eq5.4}). Moreover, by \cite[Lemma 3.1]{PM}, $u(s)$ does not change sign. We assume that $u(s)$ is positive for any $s \in (0,\epsilon)$.
		
		Claim 1: $\lambda(s) < \lambda_{1}$. Suppose to the contrary that there exist $\lambda^{\ast} \geq \lambda_{1}$, $u^{\ast} > 0$ such that $(\lambda^{\ast}, u^{\ast})$ solves (\ref{eq5.4}). This implies that $\mu = 1$ is the first eigenvalue of the problem
		\begin{equation} \label{eq5.8}
			\left\{
			\begin{array}{cc}
				-\Delta_{p} u = \mu(\lambda^{\ast} + m(x))|u|^{p-1}u \ in \ B_{R}, \\
				u = 0 \ on \ \partial B_{R},
			\end{array}
			\right.
		\end{equation}
		where $m(x) = f(|x|,u^{\ast})/(u^{\ast})^{p} > 0$. Therefore,
		\begin{equation} \label{eq5.9}
			1 = \inf\{\frac{\int_{B_{R}}|\nabla u|^{p}dx}{\int_{B_{R}}(\lambda^{\ast} + m(x))u^{p}dx}: u \in W_{0}^{1,p}(B_{R}), \int_{B_{R}}u^{p}dx = 1\}.
		\end{equation}
		Suppose that $e_{1} > 0$ is the unit eigenfunction with respect to $\lambda_{1}$, i.e.
		\[
		\int_{B_{R}}|\nabla e_{1}|^{p}dx = \lambda_{1}\int_{B_{R}}e_{1}^{p}dx.
		\]
		Then
		\begin{equation} \label{eq5.10}
			\frac{\int_{B_{R}}|\nabla e_{1}|^{p}dx}{\int_{B_{R}}(\lambda^{\ast} + m(x))e_{1}^{p}dx} = \frac{\lambda_{1}\int_{B_{R}}e_{1}^{p}dx}{\int_{B_{R}}(\lambda^{\ast} + m(x))e_{1}^{p}dx} < 1,
		\end{equation}
		contradicting (\ref{eq5.9}) and the claim is proved thereby.
		
		Claim 2: Ground states of (\ref{eq5.4}) do not change sign. Otherwise, let $u_{\lambda}$ be a ground state and $u_{\lambda}^{+} \neq 0$, $u_{\lambda}^{-} \neq 0$, where $u_{\lambda}^{+} = \max\{u_{\lambda},0\}$ and $u_{\lambda}^{-} = \min\{u_{\lambda},0\}$. $u_{\lambda}^{+}$ and $u_{\lambda}^{-}$ are elements in $W_{0,rad}^{p}$. Furthermore, it is not difficult to verify that $u_{\lambda}^{+}$ and $u_{\lambda}^{-}$ are in $\mathcal{N}_{\lambda}$. Thus $\Phi_{\lambda}(u_{\lambda}^{+}) > 0$, $\Phi_{\lambda}(u_{\lambda}^{-}) > 0$ and we derive a self-contradictory inequality
		\begin{equation} \label{eq5.11}
			\Phi_{\lambda}(u_{\lambda}^{+}) < \Phi_{\lambda}(u_{\lambda}^{+}) + \Phi_{\lambda}(u_{\lambda}^{-}) = \Phi_{\lambda}(u_{\lambda}) \leq \Phi_{\lambda}(u_{\lambda}^{+}).
		\end{equation}
		Therefore, any ground state of (\ref{eq5.4}) does not change sign.
		
		Since (\ref{eq5.4}) admits at most one positive solution in $W_{0,rad}^{1,p}(B_{R})$ for any fixed $\lambda < \lambda_{1}$, we know that $u(s)$ is the positive ground state with respect to $\lambda(s) < \lambda_{1}$ for $s \in (0,\epsilon)$. Thus
		\[
		\lim_{\lambda \rightarrow \lambda_{1}^{-}}d(\lambda) = \lim_{s \rightarrow 0^{+}}\frac{1}{p}\int_{B_{R}}u(s)^{p}dx = 0.
		\]
		Then similar to the proof of Corollary \ref{cor2.4}, we can complete the proof with $b = \max_{\lambda < \lambda_{1}}d(\lambda)$. \qed\vskip 5pt
		
		\subsection{Examples in which the uniqueness of radial positive solution in a ball is known}
		
		In this subsection, we will give some examples in which the uniqueness of radial positive solution is known and $(f_{5})$ - $(f_{7})$ hold. Consider the equation
		\begin{equation} \label{eq5.12}
			\left\{
			\begin{array}{cc}
				-\Delta_{p} u = \lambda |u|^{p-1}u + \frac{1}{(1 + |x|^{s})^{k}}|u|^{q - 2}u \ in \ B_{R}, \\
				u = 0 \ on \ \partial B_{R},
			\end{array}
			\right.
		\end{equation}
		where $B_{R} = \{x \in \mathbb{R}^{N}: |x| < R\}$, $N \geq 2$, $p + p^{2}/N < q < 2_{p}^{\ast}$, $s \geq 1$ and $k \geq 0$ (more conditions will be given later). Obviously, $f \in C^{1}([0,R) \times \mathbb{R}, \mathbb{R})$ and $(f_{5})$ - $(f_{7})$ hold with $\alpha = \beta = q$ and $\tau = 0$. To show the uniqueness of radial positive solution, we recall some results in \cite{NT, Kor}.
		
		\begin{lemma}
			\label{lem5.10} (\cite[Theorems 1.2, 1.3]{NT}) $g:[0, R] \times [0, \infty) \rightarrow [0,\infty)$ is $C^{1}$, $g(r,0) = 0$, $\forall r \geq 0$, and $g(r,u) > 0$, $\forall r > 0$, $u > 0$.
			
			$(i)$ Assume that $N > p$ and the following conditions are satisfied for $u > 0$ and $0 < r \leq R$:
			\begin{equation} \label{eq5.13}
				(p - 1)g(r,u) - g_{u}(r,u)u < 0,
			\end{equation}
			\begin{eqnarray} \label{eq5.14}
				&& -\frac{N - p^{2}}{N - p}g(r,u) - g_{u}(r,u)u + \frac{p - 1}{N - p}rg_{r}(r,u) \nonumber \\
				&\leq& (2 - p)g(r,u)\chi_{(2,+\infty)}(p),
			\end{eqnarray}
			\begin{equation} \label{eq5.15}
				\frac{N(p - 1)}{N - p}g(r,u) - g_{u}(r,u)u + \frac{p - 1}{N - p}rg_{r}(r,u) \geq 0.
			\end{equation}
			Then the equation $\Delta_{p} u + g(|x|,u) = 0$ in $B_{R}$, $u_{|\partial B_{R}} = 0$ has at most one radial positive solution.
			
			$(ii)$ Assume that $N = p$ and there exists a real $\delta > 0$ such that, for all $u > 0$ and $0 < r \leq R$:
			\begin{equation} \label{eq5.16}
				(p - 1)g(r,u) - g_{u}(r,u)u < 0,
			\end{equation}
			\begin{equation} \label{eq5.17}
				\frac{(2 - \delta)p - \delta}{\delta}g(r,u) - g_{u}(r,u)u + \frac{2}{\delta}rg_{r}(r,u) \leq (2 - p)g(r,u)\chi_{(2,+\infty)}(p),
			\end{equation}
			\begin{equation} \label{eq5.18}
				\frac{(2 + \delta)p - \delta}{\delta}g(r,u) - g_{u}(r,u)u + \frac{2}{\delta}rg_{r}(r,u) \geq 0.
			\end{equation}
			Then the equation $\Delta_{p} u + g(|x|,u) = 0$ in $B_{R}$, $u_{|\partial B_{R}} = 0$ has at most one radial positive solution.
		\end{lemma}
		
		\begin{lemma}
			\label{lem5.11} (\cite[Theorem 2.5]{Kor}) Assume that
			\[
			\min\{2,p\} < q < \frac{pN}{N - p},
			\]
			and
			\[
			a(r), b(r) \in C^{1}([0,R]), a(r) \geq 0, b(r) > 0, a'(r) \geq 0, b'(r) < 0 \ for \ r \in (0,R).
			\]
			Define
			\[
			A(r) = pa(r) + ra'(r), B(r) = (\frac{pN}{q} - (N - p))b(r) + \frac{prb'(r)}{q}.
			\]
			Assume also that the function $A(r)$ is positive and nondecreasing while the function $B(r)$ is positive in $(0, R)$, and that the functions $rb'(r)/b(r)$ and $rb'(r)$ are nonincreasing in $(0, R)$. Then the equation $\Delta_{p} u - a(|x|)u^{p - 1} + b(|x|)u^{q - 1} = 0$ in $B_{R}$, $u_{|\partial B_{R}} = 0$ has at most one radial positive solution.
		\end{lemma}
		
		\begin{theorem}
			\label{thm5.12} $(i)$ Assume $p < N$, $0 < \lambda < \lambda_{1}$ and
			\[
			2 < p \leq \frac{\sqrt{1 + 8N} - 1}{2}, ks \leq 2p, p - 1 + \frac{N - p^{2}}{N - p} \leq q \leq \frac{Np + (1 - p)ks - p}{N - p},
			\]
			or
			\[
			p \leq \min\{2,\frac{N + 1}{2}\}, ks \leq \frac{p(N - p)}{p - 1}, \frac{p(p - 1)}{N - p} \leq q \leq \frac{Np + (1 - p)ks - p}{N - p}.
			\]
			Then (\ref{eq5.12}) admits an unique radial positive solution.
			
			$(ii)$ Assume $p = N$, $0 < \lambda < \lambda_{1}$, $ks < p$ and $2 < q \leq 2p + 2 - (1 + \frac{2}{p})ks$. Then (\ref{eq5.12}) admits an unique radial positive solution.
			
			$(iii)$ Assume $\lambda < 0$, $ks \leq N - \frac{q(N - p)}{p}$ and $k \leq \frac{1}{R}$. Then (\ref{eq5.12}) admits an unique radial positive solution.
			
			$(iv)$ Assume
			\[
			p = N, ks < \frac{2p}{p + 2}, ks \leq N, k \leq \frac{1}{R}, 2p < q \leq 2p + 2 - (1 + \frac{2}{p})ks,
			\]
			or
			\begin{eqnarray}
				&& p < N, \max\{2,\sqrt{N}\} < p \leq \frac{\sqrt{1 + 8N} - 1}{2},  \nonumber \\
				&& ks \leq \min\{2p,N - \frac{q(N - p)}{p}\}, ks < \frac{p(p^{2} - N)}{N(p - 1)}, \nonumber \\
				&& k \leq \frac{1}{R}, p - 1 + \frac{N - p^{2}}{N - p} \leq q \leq \frac{Np + (1 - p)ks - p}{N - p}, q > \frac{pN + p^{2}}{N}, \nonumber
			\end{eqnarray}
			or
			\begin{eqnarray}
				&& p < N, \sqrt{N} < p \leq \min\{2,\frac{N + 1}{2}\},  \nonumber \\
				&& ks \leq \min\{\frac{p(N - p)}{p - 1},N - \frac{q(N - p)}{p}\}, ks < \frac{p(p^{2} - N)}{N(p - 1)}, \nonumber \\
				&& k \leq \frac{1}{R}, \frac{p(p - 1)}{N - p} \leq q \leq \frac{Np + (1 - p)ks - p}{N - p}, q > \frac{pN + p^{2}}{N}. \nonumber
			\end{eqnarray}
			Then there exists some $b > 0$ such that (\ref{eq5.12}) admits at least one solution $(\lambda, u_{\lambda}) \in (-\infty, \lambda_{1}) \times W_{0,rad}^{1,p}(B_{R})$ where $\lambda < \lambda_{1}$, $u_{\lambda} > 0$ and $\int_{B_{R}}u_{\lambda}^{p}dx = pb$; (\ref{eq5.12}) admits at least two solutions $(\lambda, u_{\lambda}), (\tilde{\lambda}, u_{\tilde{\lambda}}) \in (-\infty, \lambda_{1}) \times W_{0,rad}^{1,p}(B_{R})$ where $\lambda < \tilde{\lambda} < \lambda_{1}$, $u_{\lambda} > 0$, $u_{\tilde{\lambda}} > 0$ and $\int_{B_{R}}u_{\lambda}^{p}dx = \int_{B_{R}}\tilde{u}_{\lambda}^{p}dx = pc$ for any $0 < c < b$; and (\ref{eq5.12}) admits no $(\lambda, u_{\lambda}) \in (-\infty, \lambda_{1}) \times W_{0,rad}^{1,p}(B_{R})$ where $u_{\lambda} > 0$ and $\int_{B_{R}}u_{\lambda}^{p}dx = pc$ for any $c > b$.
		\end{theorem}
		
		\textit{Proof.  } First we prove $(i)$. The existence of a radial positive solution can be shown by the Nehari manifold method or by the mountain-pass theorem. To show the uniqueness, we apply Lemma \ref{lem5.10} $(i)$ with
		\[
		g(r,u) = \lambda u^{p - 1} + \frac{1}{(1 + r^{s})^{k}}u^{q - 1}, u > 0, 0 < r \leq R.
		\]
		We verify (\ref{eq5.13}) - (\ref{eq5.15}):
		\begin{equation} \label{eq5.19}
			(p - 1)g(r,u) - g_{u}(r,u)u = (p - q) \frac{1}{(1 + r^{s})^{k}}u^{q - 1} < 0,
		\end{equation}
		\begin{eqnarray} \label{eq5.20}
			&& -\frac{N - p^{2}}{N - p}g(r,u) - g_{u}(r,u)u + \frac{p - 1}{N - p}rg_{r}(r,u) \nonumber \\
			&=& \frac{2p^{2} - (N + 1)p}{N - p}\lambda u^{p - 1} - (\frac{N - p^{2}}{N - p} + q - 1)\frac{1}{(1 + r^{s})^{k}}u^{q - 1} \nonumber \\
			&& - ks\frac{p - 1}{N - p}\frac{r^{s}}{(1 + r^{s})^{k+1}}u^{q - 1} \nonumber \\
			&\leq& (2 - p)g(r,u)\chi_{(2,+\infty)}(p),
		\end{eqnarray}
		\begin{eqnarray} \label{eq5.21}
			&& \frac{N(p - 1)}{N - p}g(r,u) - g_{u}(r,u)u + \frac{p - 1}{N - p}rg_{r}(r,u) \nonumber \\
			&=& \frac{N(p - 1)}{N - p}\lambda u^{p - 1} + (\frac{(N - ks)(p - 1)}{N - p} - (q - 1))\frac{1}{(1 + r^{s})^{k}}u^{q - 1} \nonumber \\
			&& + ks\frac{p - 1}{N - p}\frac{1}{(1 + r^{s})^{k+1}}u^{q - 1} \nonumber \\
			&\geq& 0.
		\end{eqnarray}
		Thus the proof of $(i)$ is completed.
		
		To prove $(ii)$, we apply Lemma \ref{lem5.10} $(ii)$ with
		\[
		g(r,u) = \lambda u + \frac{1}{(1 + r^{s})^{k}}u^{p - 1}, u > 0, 0 < r \leq R.
		\]
		We verify (\ref{eq5.16}) - (\ref{eq5.18}) for $\delta \in [\frac{2p}{p + 2},\frac{2(p - ks)}{q - p}]$:
		\begin{equation} \label{eq5.22}
			(p - 1)g(r,u) - g_{u}(r,u)u = (p - q) \frac{1}{(1 + r^{s})^{k}}u^{q - 1} < 0, \nonumber
		\end{equation}
		\begin{eqnarray} \label{eq5.23}
			&& \frac{(2 - \delta)p - \delta}{\delta}g(r,u) - g_{u}(r,u)u + \frac{2}{\delta}rg_{r}(r,u) \nonumber \\
			&=& (\frac{(2 - \delta)p - \delta}{\delta} - (p - 1))\lambda u^{p - 1} + (\frac{(2 - \delta)p - \delta}{\delta} - (q - 1))\frac{1}{(1 + r^{s})^{k}}u^{q - 1} \nonumber \\
			&& - \frac{2}{\delta}ks\frac{r^{s}}{(1 + r^{s})^{k + 1}}u^{q - 1} \nonumber \\
			&\leq& (2 - p)g(r,u)\chi_{(2,+\infty)}(p), \nonumber
		\end{eqnarray}
		\begin{eqnarray} \label{eq5.24}
			&& \frac{(2 + \delta)p - \delta}{\delta}g(r,u) - g_{u}(r,u)u + \frac{2}{\delta}rg_{r}(r,u) \nonumber \\
			&=& (\frac{(2 + \delta)p - \delta}{\delta} - (p - 1))\lambda u + (\frac{(2 - \delta)p - \delta - 2ks}{\delta} - (q - 1))\frac{1}{(1 + r^{s})^{k}}u^{q - 1} \nonumber \\
			&& + 2ks \frac{1}{(1 + r^{s})^{k+1}}u^{q - 1} \geq 0. \nonumber
		\end{eqnarray}
		The proof of $(ii)$ is completed.
		
		Applying Lemma \ref{lem5.11} with $a(r) = -\lambda$, $b(r) = \frac{1}{(1 + r^{s})^{k}}u^{p - 1}$, $u > 0$, $0 < r \leq R$, we can prove $(iii)$ since it is not difficult to verify the conditions in Lemma \ref{lem5.11} under the hypotheses that $\lambda \leq 0$, $ks \leq N - \frac{q(N - p)}{p}$ and $k \leq \frac{1}{R}$.
		
		Since we have shown the uniqueness of radial positive solution of (\ref{eq5.12}), $(iv)$ follows Theorem \ref{thm5.9}.
		\qed\vskip 5pt

		\section{\textbf{Mixed fractional Laplacians}} \label{mfl}
		
		The purpose of this section is to apply our abstract framework 3 to a class of nonlinear equations involving the mixed fractional Laplacians. This type of equations arises in various fields ranging from biophysics to population dynamics and has very recently
		received an increasing interest. The Cauchy problems involving mixed fractional laplacians have been recently addressed in \cite{ HH2, HH1}. The normalized solutions were tacked in \cite{HH5, HH3, HH4}. In this section, we study other aspects of this operator.
		
		We are interested in the following fractional equation
		\begin{equation} \label{eq7.1}
			\left\{
			\begin{array}{cc}
				(-\Delta)^{s_{1}} u + (-\Delta)^{s_{2}} u = \lambda u + h_{2}(|x|)|u|^{p-2}u + h_{3}(|x|)|u|^{q-2}u \ in \ \mathbb{R}^{N},\\
				u(x) \rightarrow 0 \ as \ |x| \rightarrow +\infty,
			\end{array}
				\right.
		\end{equation}
		where $0 < s_{1} < s_{2} \leq 1$, $\max\{2,4s_{1}\} \leq N < 2s_{1}s_{2}/(s_{2} - s_{1})$, $2 < p < q < 2_{s_{2}}^{\ast}$ and we assume that
		
		$(h_{2})$ $h_{2}(r) \in C^{1}([0,+\infty))$, $\lim_{r \rightarrow +\infty}h_{2}(r) = h_{2}^{\infty} > 0$, $\inf_{r > 0}h_{2}(r) > 0$.
		
		$(h_{3})$ $h_{3}(r) \in L^{\infty}([0,+\infty)) \cap C^{1}([0,+\infty))$, $\inf_{r > 0}h_{3}(r) > 0$.
		
		Note that the abstract framework 2 cannot be applied directly since $(S_{2})$ does not hold. Fortunately, we can study this problem after some modifications. More precisely, let
		$$
		W = H^{s_{1}}_{rad}(\mathbb{R}^{N}) \cap H^{s_{2}}_{rad}(\mathbb{R}^{N}),
		$$
		\begin{eqnarray}
		\Phi_{\lambda}(u) &=& \frac{1}{2}\int_{\mathbb{R}^{N}}[|(-\Delta)^{\frac{s_{1}}{2}} u|^{2} + |(-\Delta)^{\frac{s_{2}}{2}} u|^{2} - \lambda u^{2}]dx \nonumber \\
		&-& \int_{\mathbb{R}^{N}}[\frac{1}{p}h_{2}(|x|)|u|^{p} + \frac{1}{q}h_{3}(|x|)|u|^{q}]dx,
		\end{eqnarray}
        \begin{eqnarray}
		\mathcal{N}_{\lambda} = \{u \in W: && \int_{\mathbb{R}^{N}}(|(-\Delta)^{\frac{s_{1}}{2}}u|^{2} + |(-\Delta)^{\frac{s_{2}}{2}}u|^{2})dx \nonumber \\
		&=& \int_{\mathbb{R}^{N}}[\lambda u^{2} + h_{2}(|x|)|u|^{p} + h_{3}(|x|)|u|^{q}]dx\}.
		\end{eqnarray}
		
		Similar to the case of a single fractional Laplacian, we can prove that $\Phi_{\lambda}$ satisfies (PS) condition in $W$ and $h(\lambda) = \inf_{u \in \mathcal{N}_{\lambda}}\Phi_{\lambda}(u)$ can be achieved by a positive ground state $u_{\lambda}$. When $\lambda \rightarrow 0^{-}$, we set
		$$
		w_{\lambda} = |\lambda|^{-\frac{1}{p-2}}u_{\lambda}(\frac{x}{|\lambda|^{1/2s_{1}}}),
		$$
    	which satisfies the following equation
    	\begin{eqnarray} \label{eq7.2}
	    && (-\Delta)^{s_{1}}w + |\lambda|^{\frac{s_{2}-s_{1}}{s_{1}}}(-\Delta)^{s_{2}}w \nonumber \\
	    &=& -w + h_{2}(\frac{|x|}{|\lambda|^{1/2s_{1}}})|w|^{p-2}w + |\lambda|^{\frac{q-p}{p-2}}h_{3}(\frac{|x|}{|\lambda|^{1/2s_{1}}})|w|^{q-2}w.
	    \end{eqnarray}
	    Let
	    \begin{eqnarray}
	    \widehat{\Phi}_{\lambda}(w) &=& \frac{1}{2}\int_{\mathbb{R}^{N}}[|(-\Delta)^{\frac{s_{1}}{2}} w|^{2} + |\lambda|^{\frac{s_{2}-s_{1}}{s_{1}}}|(-\Delta)^{\frac{s_{2}}{2}} w|^{2} + w^{2}]dx \nonumber \\
	    && - \int_{\mathbb{R}^{N}}[\frac{1}{p}h_{2}(\frac{|x|}{|\lambda|^{1/2s_{1}}})|w|^{p} + \frac{1}{q}|\lambda|^{\frac{q-p}{p-2}}h_{3}(\frac{|x|}{|\lambda|^{1/2s_{1}}})|w|^{q}]dx, \nonumber
    	\end{eqnarray}
	    $$
	    \widehat{\Phi}_{0}(w) =  \frac{1}{2}\int_{\mathbb{R}^{N}}[|(-\Delta)^{\frac{s_{1}}{2}} w|^{2} + w^{2}]dx - \frac{h_{2}^{\infty}}{p}\int_{\mathbb{R}^{N}}|w|^{p}dx,
	    $$
	    \begin{eqnarray}
	    \widehat{\mathcal{N}}_{\lambda} = \{w \in W: && \int_{\mathbb{R}^{N}}[|(-\Delta)^{\frac{s_{1}}{2}}w|^{2} + |\lambda|^{\frac{s_{2}-s_{1}}{s_{1}}}|(-\Delta)^{\frac{s_{2}}{2}}w|^{2} + w^{2}]dx \nonumber \\
	    &=&   \int_{\mathbb{R}^{N}}[h_{2}(\frac{|x|}{|\lambda|^{1/2s_{1}}})|w|^{p} - |\lambda|^{\frac{q-p}{p-2}}h_{3}(\frac{|x|}{|\lambda|^{1/2s_{1}}})|w|^{q}]dx\},
	    \end{eqnarray}
	    $$
	    \widehat{\mathcal{N}}_{0} = \{w \in H^{s_{1}}_{rad}(\mathbb{R}^{N}): \int_{\mathbb{R}^{N}}[|(-\Delta)^{\frac{s_{1}}{2}}w|^{2} + w^{2} - h_{2}^{\infty}|w|^{p}]dx = 0\}.
	    $$
		
		\begin{lemma}
			\label{lem7.1}  Assume that $2 < p < 2_{s_{1}}^{\ast}$ and $h_{2}^{\infty} > 0$. Let $w^{\ast}$ be the unique and non-degenerate positive ground state in $H^{s_{1}}_{rad}(\mathbb{R}^{N})$ of
			$$
			(-\Delta)^{s_{1}}w = -w + h_{2}^{\infty}|w|^{p-2}w.
			$$
			Then $w^{\ast} \in H^{s_{2}}_{rad}(\mathbb{R}^{N})$.
		\end{lemma}
		
		\textit{Proof.  }  The uniqueness and non-degeneracy of $w^{\ast}$ have been shown by \cite{FLS}. Then we can simply follow the arguments in \cite[Lemma B.2]{FL} to show that $w^{\ast} \in H^{2s_{1}+1}_{rad}(\mathbb{R}^{N})$. Note that $s_{2} \leq 1 < 2s_{1}+1$. Thus $w^{\ast} \in H^{s_{2}}_{rad}(\mathbb{R}^{N})$.
		\qed\vskip 5pt
		
		\begin{lemma}
			\label{lem7.2}  Assume that $2 < p < q < 2_{s_{1}}^{\ast}$ and $(h_{2})$, $(h_{3})$ hold. Then there exists $\delta > 0$ such that for any $(\lambda, w) \in (-\delta,0] \times W \backslash \{0\}$, there exists a unique function $\widehat{t}: (-\delta,0] \times W \backslash \{0\} \rightarrow (0, \infty)$ such that
			\begin{equation}
			when \ t > 0, \ then \ (\lambda, t w) \in \mathcal{\widehat{N}} \Leftrightarrow t = \widehat{t}(\lambda, w),
			\nonumber
			\end{equation}
			$\widehat{\Phi}_{\lambda}(\widehat{t}(\lambda, w)w) = \max_{t > 0}\widehat{\Phi}_{\lambda}(tw)$, $\widehat{t} \in C((-\delta,0] \times W \backslash \{0\}, (0, \infty))$, where $\mathcal{\widehat{N}} := \cup_{-\delta < \lambda \leq 0}\mathcal{\widehat{N}}_{\lambda}$.
		\end{lemma}
		
		\textit{Proof.  } The proof is similar to the one of Lemma \ref{lem4.1}.		
		\qed\vskip 5pt
		
		\begin{lemma}
			\label{lem7.3}  Assume that $2 < p < q < 2_{s_{1}}^{\ast}$ and $(h_{2})$, $(h_{3})$ hold. Then $w_{\lambda}$ is bounded in $H^{s_{1}}(\mathbb{R}^{N})$ and $|\lambda|^{\frac{s_{2}-s_{1}}{s_{1}}}\int_{\mathbb{R}^{N}}|(-\Delta)^{\frac{s_{2}}{2}}w_{\lambda}|^{2}dx$ is bounded when $\lambda \rightarrow 0^{-}$.
		\end{lemma}
		
		\textit{Proof.  } Lemma \ref{lem7.1} yields that $w^{\ast} \in W$. Then thanks to Lemma \ref{lem7.2},  we can follow the arguments in Step 1 of Theorem \ref{thmc.1} to show that
		$$
		\limsup_{\lambda \rightarrow 0^{-}} \widehat{h}(\lambda) \leq \widehat{h}(0).
		$$
		Mimicking the proof of Lemma \ref{lem4.2}, we arrive at the statement.
		\qed\vskip 5pt
		
		\begin{theorem}
			\label{thm7.4}  Assume that $2 < p < q < 2_{s_{1}}^{\ast}$ and $(h_{2})$, $(h_{3})$ hold. Then $w_{\lambda} \rightarrow w^{\ast}$ in $H^{s_{1}}(\mathbb{R}^{N})$ as $\lambda \rightarrow 0^{-}$ where $w^{\ast}$ be the unique and positive ground state in $H^{s_{1}}_{rad}(\mathbb{R}^{N})$ of
			$$
			(-\Delta)^{s_{1}}w = -w + h_{2}^{\infty}|w|^{p-2}w.
			$$
		\end{theorem}
		
		\textit{Proof.  } Let $w_{n} = w_{\lambda_{n}}$ when $\lambda_{n} \rightarrow 0^{-}$. By Lemma \ref{lem7.3}, $w_{n}$ is bounded in $H^{s_{1}}(\mathbb{R}^{N})$. Up to a subsequence, we assume that $w_{n} \rightharpoonup w_{0}$ in $H^{s_{1}}(\mathbb{R}^{N})$. Noticing $2 < p < q < 2_{s_{1}}^{\ast}$ and thus the embeddings $H^{s_{1}}_{rad}(\mathbb{R}^{N}) \hookrightarrow L^{p}(\mathbb{R}^{N})$, $H^{s_{1}}_{rad}(\mathbb{R}^{N}) \hookrightarrow L^{q}(\mathbb{R}^{N})$ are compact, we have $w_{n} \rightarrow w_{0}$ in $L^{p}(\mathbb{R}^{N})$ and $L^{q}(\mathbb{R}^{N})$.
		
		Step 1: $w_{n} \rightarrow w_{0}$ in $H^{s_{1}}(\mathbb{R}^{N})$.
		
		For any $\phi \in C_{0}^{\infty}(\mathbb{R}^{N})$, we have
		\begin{equation} \label{eq7.4}
		\lim_{n \rightarrow +\infty}\int_{\mathbb{R}^{N}}[h_{2}(\frac{|x|}{|\lambda_{n}|^{1/2s_{1}}})-h_{2}^{\infty}]|w_{n}|^{p-2}w_{n}\phi dx = 0
		\end{equation}
		by Lebesgue convergent theorem since $w_{n}$ is convergent in $L^{p}(\mathbb{R}^{N})$;
		\begin{equation} \label{eq7.5}
		\lim_{n \rightarrow +\infty}\int_{\mathbb{R}^{N}}|\lambda_{n}|^{\frac{q-p}{p-2}}h_{3}(\frac{|x|}{|\lambda_{n}|^{1/2s_{1}}})|w_{n}|^{q-2}w_{n}\phi dx = 0
		\end{equation}
		since $w_{n}$ is bounded in $L^{q}(\mathbb{R}^{N})$;
		\begin{eqnarray} \label{eq7.6}
		&& \lim_{n \rightarrow +\infty}|\int_{\mathbb{R}^{N}}|\lambda_{n}|^{\frac{s_{2}-s_{1}}{s_{1}}}(-\Delta)^{\frac{s_{2}}{2}}w_{n}(-\Delta)^{\frac{s_{2}}{2}}\phi dx| \nonumber \\
		&\leq& \lim_{\lambda \rightarrow 0^{-}}(\int_{\mathbb{R}^{N}}|\lambda_{n}|^{\frac{s_{2}-s_{1}}{s_{1}}}|(-\Delta)^{\frac{s_{2}}{2}}w_{n}|^{2}dx)^{\frac{1}{2}}(\int_{\mathbb{R}^{N}}|\lambda_{n}|^{\frac{s_{2}-s_{1}}{s_{1}}}|(-\Delta)^{\frac{s_{2}}{2}}\phi|^{2}dx)^{\frac{1}{2}} \nonumber \\
		&=& 0 \nonumber
		\end{eqnarray}
		since $\int_{\mathbb{R}^{N}}|\lambda_{n}|^{\frac{s_{2}-s_{1}}{s_{1}}}|(-\Delta)^{\frac{s_{2}}{2}}w_{n}|^{2}dx$ is bounded. Therefore, for any $\phi \in C_{0}^{\infty}(\mathbb{R}^{N})$
		\begin{eqnarray} \label{eq7.7}
		D_{w}\widehat{\Phi}_{0}(w_{n})(\phi) &=& D_{w}\widehat{\Phi}_{\lambda_{n}}(w_{n})(\phi) - \int_{\mathbb{R}^{N}}|\lambda_{n}|^{\frac{s_{2}-s_{1}}{s_{1}}}(-\Delta)^{\frac{s_{2}}{2}}w_{n}(-\Delta)^{\frac{s_{2}}{2}}\phi dx \nonumber \\
		&+& \int_{\mathbb{R}^{N}}[h_{2}(\frac{|x|}{|\lambda_{n}|^{1/2s_{1}}})-h_{2}^{\infty}]|w_{n}|^{p-2}w_{n}\phi dx \nonumber \\
		&+& \int_{\mathbb{R}^{N}}|\lambda_{n}|^{\frac{q-p}{p-2}}h_{3}(\frac{|x|}{|\lambda_{n}|^{1/2s_{1}}})|w_{n}|^{q-2}w_{n}\phi dx \nonumber \\
		&\rightarrow& 0.
		\end{eqnarray}
		Note that $C_{0}^{\infty}(\mathbb{R}^{N})$ is dense in $H^{s_{1}}(\mathbb{R}^{N})$. We know
		$$
		D_{w}\widehat{\Phi}_{0}(w_{n})(\phi) \rightarrow 0, \forall \phi \in H^{s_{1}}(\mathbb{R}^{N}).
		$$
		Then by a standard argument, we can prove that $w_{n} \rightarrow w_{0}$ in $H^{s_{1}}(\mathbb{R}^{N})$.
		
		Step 2: $w_{0} \neq 0$.
		
		Let us argue by contradiction. Suppose to the contrary that $w_{0} = 0$. Notice that $D_{ww}\widehat{\Phi}_{0}(0) = (-\Delta)^{s_{1}} + 1$ is invertible. Then applying implicit function arguments, we can show that $w_{\lambda} \equiv 0$ is a neighborhood of $0$ since $0$ is always a solution of $D_{w}\widehat{\Phi}_{\lambda}(w) = 0$, which is a contradiction.
		
		Step 3: $w_{0} = w^{\ast}$.
		
		Since $w_{0} \neq 0$ and $D_{w}\widehat{\Phi}_{0}(w_{0}) = 0$, $w_{0} \in \widehat{\mathcal{N}}_{0}$. We will show that $\widehat{\Phi}_{0}(w_{0}) \leq \widehat{h}(0)$ to complete the proof. Indeed,
		\begin{eqnarray}
		&& \widehat{\Phi}_{0}(w_{n}) \nonumber \\
		&=& \widehat{\Phi}_{\lambda_{n}}(w_{n}) - \frac{1}{2}\int_{\mathbb{R}^{N}}|\lambda_{n}|^{\frac{s_{2}-s_{1}}{s_{1}}}|(-\Delta)^{\frac{s_{1}}{2}}w_{n}|^{2} dx
		+ \nonumber \\
		&& \frac{1}{p}\int_{\mathbb{R}^{N}}[h_{2}(\frac{|x|}{|\lambda_{n}|^{1/2s_{1}}})-h_{2}^{\infty}]|w_{n}|^{p} dx + \frac{1}{q}\int_{\mathbb{R}^{N}}|\lambda_{n}|^{\frac{q-p}{p-2}}h_{3}(\frac{|x|}{|\lambda_{n}|^{1/2s_{1}}})|w_{n}|^{q} dx \nonumber \\
		&=& \widehat{\Phi}_{\lambda_{n}}(w_{n}) - \frac{1}{2}\int_{\mathbb{R}^{N}}|\lambda_{n}|^{\frac{s_{2}-s_{1}}{s_{1}}}|(-\Delta)^{\frac{s_{1}}{2}}w_{n}|^{2} dx + o_{n}(1) \nonumber \\
		&\leq& \widehat{h}(\lambda_{n}) + o_{n}(1). \nonumber
		\end{eqnarray}
		Thus
		$$
		\widehat{\Phi}_{0}(w_{0}) \leq \limsup_{\lambda \rightarrow 0^{-}} \widehat{h}(\lambda) \leq \widehat{h}(0),
		$$
		and then we know that $\widehat{\Phi}_{0}(w_{0}) = \widehat{h}(0)$. By the uniqueness of $w^{\ast}$, $w_{0} = w^{\ast}$.
		\qed\vskip 5pt
		
		\begin{theorem}
			\label{thm7.5}  Assume that $2 < p < q < 2_{s_{1}}^{\ast}$ and $(h_{2})$, $(h_{3})$ hold. Then there exists $\Lambda_{1} < 0$ such that when $\lambda \in (\Lambda_{1},0)$, the ground state of (\ref{eq7.1}) is unique and non-degenerate. Furthermore,
			
			$(i)$ if $p > 2 + 4s_{1}/N$, then
			$$
			\lim_{\lambda \rightarrow 0^{-}}\int_{\mathbb{R}^{N}}u_{\lambda}^{2}dx \rightarrow +\infty,
			$$
			and thus there exists $c_{1} > 0$ large such that for any $c > c_{1}$, (\ref{eq7.1}) admits a ground state $u_{\lambda}$ with $\lambda < 0$ and $|\lambda|$ small such that $\int_{\mathbb{R}^{N}}u_{\lambda}^{2}dx = c$;
			
			$(ii)$ if $2 < p < 2 + 4s_{1}/N$, then
			$$
			\lim_{\lambda \rightarrow 0^{-}}\int_{\mathbb{R}^{N}}u_{\lambda}^{2}dx \rightarrow 0,
			$$
			and thus there exists $c_{2} > 0$ small such that for any $c \in (0,c_{2}$), (\ref{eq7.1}) admits a ground state $u_{\lambda}$ with $\lambda < 0$ and $|\lambda|$ small such that $\int_{\mathbb{R}^{N}}u_{\lambda}^{2}dx = c$;
			
			$(iii)$  if $p = 2 + 4s_{1}/N$, then
			$$
			\lim_{\lambda \rightarrow 0^{-}}\int_{\mathbb{R}^{N}}u_{\lambda}^{2}dx \rightarrow \int_{\mathbb{R}^{N}}(w^{\ast})^{2}dx.
			$$
		\end{theorem}
		
		\textit{Proof.  } Along the lines of Step 3 in the proof of Theorem \ref{thmc.1}, we can obtain the uniqueness and non-degeneracy. Then $(i)$ - $(iii)$ are similar to Corollaries \ref{corc.2} and \ref{corc.3}.
		\qed\vskip 5pt
		
		 When $\lambda \rightarrow -\infty$, we set $\mu = 1/\lambda \rightarrow 0^{-}$ and
		 $$
		 v_{\mu} = |\mu|^{\frac{1}{q-2}}u_{\lambda}(|\mu|^{\frac{1}{2s_{2}}}|x|),
		 $$
		 which satisfies the following equation
		 \begin{eqnarray} \label{eq7.9}
		 && |\mu|^{\frac{s_{2}-s_{1}}{s_{2}}}(-\Delta)^{s_{1}}v + (-\Delta)^{s_{2}}v \nonumber \\
		 &=& -v + |\mu|^{\frac{q-p}{q-2}}h_{2}(|\mu|^{\frac{1}{2s_{2}}}|x|)|v|^{p-2}v + h_{3}(|\mu|^{\frac{1}{2s_{2}}}|x|)|v|^{q-2}v.
		 \end{eqnarray}
		 Let
		 \begin{eqnarray}
		 \widetilde{\Phi}_{\mu}(v) &=& \frac{1}{2}\int_{\mathbb{R}^{N}}[|\mu|^{\frac{s_{2}-s_{1}}{s_{2}}}|(-\Delta)^{\frac{s_{1}}{2}} v|^{2} + |(-\Delta)^{\frac{s_{2}}{2}} v|^{2} + v^{2}]dx \nonumber \\
		 && - \int_{\mathbb{R}^{N}}[\frac{1}{p}|\mu|^{\frac{q-p}{q-2}}h_{2}(|\mu|^{\frac{1}{2s_{2}}}|x|)|v|^{p} + \frac{1}{q}h_{3}(|\mu|^{\frac{1}{2s_{2}}}|x|)|v|^{q}]dx, \nonumber
		 \end{eqnarray}
		 $$
		 \widetilde{\Phi}_{0}(v) =  \frac{1}{2}\int_{\mathbb{R}^{N}}[|(-\Delta)^{\frac{s_{2}}{2}} v|^{2} + v^{2}]dx - \frac{h_{3}(0)}{q}\int_{\mathbb{R}^{N}}|v|^{q}dx,
		 $$
		 \begin{eqnarray}
		 \widetilde{\mathcal{N}}_{\mu} = \{v \in W: && \int_{\mathbb{R}^{N}}[|\mu|^{\frac{s_{2}-s_{1}}{s_{2}}}||(-\Delta)^{\frac{s_{1}}{2}}v|^{2} + (-\Delta)^{\frac{s_{2}}{2}}v|^{2} + v^{2}]dx \nonumber \\
		 &=&  \int_{\mathbb{R}^{N}}[|\mu|^{\frac{q-p}{q-2}}h_{2}(|\mu|^{\frac{1}{2s_{2}}}|x|)|v|^{p} - h_{3}(|\mu|^{\frac{1}{2s_{2}}}|x|)|v|^{q}]dx\},
		 \end{eqnarray}
		 $$
		 \widetilde{\mathcal{N}}_{0} = \{v \in H^{s_{2}}_{rad}(\mathbb{R}^{N}): \int_{\mathbb{R}^{N}}[|(-\Delta)^{\frac{s_{2}}{2}}v|^{2} + v^{2} - h_{3}(0)|v|^{q}]dx = 0\}.
		 $$
		
		 \begin{lemma}
		 	\label{lem7.6}  Assume that $2 < p < q < 2_{s_{2}}^{\ast}$ and $(h_{2})$, $(h_{3})$ hold. Then there exists $\delta > 0$ such that for any $(\mu, v) \in (-\delta,0] \times W \backslash \{0\}$, there exists a unique function $\widehat{t}: (-\delta,0] \times W \backslash \{0\} \rightarrow (0, \infty)$ such that
		 	\begin{equation}
		 	when \ t > 0, \ then \ (\lambda, t w) \in \mathcal{\widetilde{N}} \Leftrightarrow t = \widetilde{t}(\mu, v),
		 	\nonumber
		 	\end{equation}
		 	$\widetilde{\Phi}_{\mu}(\widetilde{t}(\mu, v)v) = \max_{t > 0}\widetilde{\Phi}_{\mu}(tv)$, $\widetilde{t} \in C((-\delta,0] \times W \backslash \{0\}, (0, \infty))$, where $\mathcal{\widetilde{N}} := \cup_{-\delta < \mu \leq 0}\mathcal{\widetilde{N}}_{\mu}$.
		 \end{lemma}
		
		 \textit{Proof.  } The proof is similar to the one of Lemma \ref{lem4.1}.		
		 \qed\vskip 5pt
		
		 \begin{lemma}
		 	\label{lem7.7}  Assume that $2 < p < q < 2_{s_{2}}^{\ast}$ and $(h_{2})$, $(h_{3})$ hold. Then $v_{\mu}$ is bounded in $H^{s_{2}}(\mathbb{R}^{N})$ when $\mu \rightarrow 0^{-}$.
		 \end{lemma}
		
		 \textit{Proof.  } Thanks to Lemma \ref{lem7.2},  we can follow the arguments of Step 1 in Theorem \ref{thmc.4} to show that
		 $$
		 \limsup_{\mu \rightarrow 0^{-}} \widetilde{h}(\mu) \leq \widetilde{h}(0).
		 $$
		 Mimicking the proof of Lemma \ref{lem4.2}, we can conclude that.
		 \qed\vskip 5pt
		
		 \begin{theorem}
		 	\label{thm7.8}  Assume that $2 < p < q < 2_{s_{2}}^{\ast}$ and $(h_{2})$, $(h_{3})$ hold. Then $v_{\mu} \rightarrow v^{\ast}$ in $H^{s_{2}}(\mathbb{R}^{N})$ as $\mu \rightarrow 0^{-}$ where $v^{\ast}$ is the unique and non-degenerate positive ground state in $H^{s_{2}}_{rad}(\mathbb{R}^{N})$ of
		 	$$
		 	(-\Delta)^{s_{2}}v = -v + h_{3}(0)|v|^{q-2}v.
		 	$$
		 \end{theorem}
		
		 \textit{Proof.  } Let $v_{n} = v_{\mu_{n}}$ when $\mu_{n} \rightarrow 0^{-}$. By Lemma \ref{lem7.7}, $v_{n}$ is bounded in $H^{s_{2}}(\mathbb{R}^{N})$. For $0 < s_{1} < s_{2} \leq 1$, we recall from \cite[Remark 1.4.1]{Ca} that $H^{s_{1}}(\mathbb{R}^{N}) \hookrightarrow H^{s_{2}}(\mathbb{R}^{N})$, and hence, $v_{n}$ is bounded in $H^{s_{1}}(\mathbb{R}^{N})$. Also, $v_{n}$ is bounded in $L^{p}(\mathbb{R}^{N})$ and $L^{q}(\mathbb{R}^{N})$. Up to a subsequence, we assume that $v_{n} \rightharpoonup v_{0}$ in $H^{s_{2}}(\mathbb{R}^{N})$ and $v_{n} \rightarrow v_{0}$ in $L^{p}(\mathbb{R}^{N})$ and $L^{q}(\mathbb{R}^{N})$.
		
		 Noticing that
		 \begin{equation} \label{eq7.11}
		 \lim_{n \rightarrow +\infty}\int_{\mathbb{R}^{N}}[h_{3}(|\mu|^{\frac{1}{2s_{2}}}|x|)-h_{3}(0)]|v_{n}|^{q}dx = 0
		 \end{equation}
		 by Lebesgue convergence theorem since $v_{n}$ is convergent in $L^{q}(\mathbb{R}^{N})$;
		 \begin{equation} \label{eq7.12}
		 \lim_{n \rightarrow +\infty}\int_{\mathbb{R}^{N}}|\mu_{n}|^{\frac{q-p}{q-2}}h_{2}(|\mu|^{\frac{1}{2s_{2}}}|x|)|v_{n}|^{p}dx = 0
		 \end{equation}
		 since $v_{n}$ is bounded in $L^{p}(\mathbb{R}^{N})$;
		 \begin{equation} \label{eq7.13}
		 \lim_{n \rightarrow +\infty}\int_{\mathbb{R}^{N}}|\mu_{n}|^{\frac{s_{2}-s_{1}}{s_{2}}}|(-\Delta)^{\frac{s_{1}}{2}}v_{n}|^{2}dx  = 0
		 \end{equation}
		 since $v_{n}$ is bounded in $H^{s_{1}}(\mathbb{R}^{N})$, we have
		 \begin{eqnarray} \label{eq7.14}
		 \widetilde{\Phi}_{0}(v_{n}) &=& \widetilde{\Phi}_{\mu_{n}}(v_{n}) - \frac{1}{2}\int_{\mathbb{R}^{N}}|\mu_{n}|^{\frac{s_{2}-s_{1}}{s_{2}}}|(-\Delta)^{\frac{s_{1}}{2}}v_{n}|^{2}dx \nonumber \\
		 &+& \frac{1}{q}\int_{\mathbb{R}^{N}}[h_{3}(|\mu|^{\frac{1}{2s_{2}}}|x|)-h_{3}(0)]|v_{n}|^{q}dx \nonumber \\
		 &+& \frac{1}{p}\int_{\mathbb{R}^{N}}|\mu_{n}|^{\frac{q-p}{q-2}}h_{2}(|\mu|^{\frac{1}{2s_{2}}}|x|)|v_{n}|^{p}dx \nonumber \\
		 &=& \widetilde{h}(\mu_{n}) + o_{n}(1).
		 \end{eqnarray}
		 Similarly,
		 $$
		 D_{v}\widetilde{\Phi}_{0}(v_{n}) \rightarrow 0.
		 $$
		 Then along the lines of the proof of Theorem \ref{thmc.4}, we can prove that $v_{n} \rightarrow v^{\ast}$ in $H^{s_{2}}(\mathbb{R}^{N})$ where $v^{\ast}$ is the unique and non-degenerate positive ground state in $H^{s_{2}}_{rad}(\mathbb{R}^{N})$ of
		 $$
		 (-\Delta)^{s_{2}}v = -v + h_{3}(0)|v|^{q-2}v.
		 $$
		 \qed\vskip 5pt
		
		 \begin{theorem}
		 	\label{thm7.9}  Assume that $2 < p < q < 2_{s_{2}}^{\ast}$ and $(h_{2})$, $(h_{3})$ hold. Then there exists $\Lambda_{2} < 0$ such that when $\lambda < \Lambda_{2}$, the ground state of (\ref{eq7.1}) is unique and non-degenerate. Furthermore,
		 	
		 	$(i)$ if $q > 2 + 4s_{2}/N$, then
		 	$$
		 	\lim_{\lambda \rightarrow -\infty}\int_{\mathbb{R}^{N}}u_{\lambda}^{2}dx \rightarrow 0,
		 	$$
		 	and thus there exists $c_{3} > 0$ small such that for any $c \in (0,c_{3})$, (\ref{eq7.1}) admits a ground state $u_{\lambda}$ with $\lambda < 0$ and $|\lambda|$ large such that $\int_{\mathbb{R}^{N}}u_{\lambda}^{2}dx = c$;
		 	
		 	$(ii)$ if $2 < q < 2 + 4s_{2}/N$, then
		 	$$
		 	\lim_{\lambda \rightarrow -\infty}\int_{\mathbb{R}^{N}}u_{\lambda}^{2}dx \rightarrow +\infty,
		 	$$
		 	and thus there exists $c_{4} > 0$ large such that for any $c > c_{4}$, (\ref{eq7.1}) admits a ground state $u_{\lambda}$ with $\lambda < 0$ and $|\lambda|$ large such that $\int_{\mathbb{R}^{N}}u_{\lambda}^{2}dx = c$;
		 	
		 	$(iii)$  if $q = 2 + 4s_{2}/N$, then
		 	$$
		 	\lim_{\lambda \rightarrow -\infty}\int_{\mathbb{R}^{N}}u_{\lambda}^{2}dx \rightarrow \int_{\mathbb{R}^{N}}(v^{\ast})^{2}dx.
		 	$$
		 \end{theorem}
		
		 \begin{remark}
		      There is a difference between cases when $\lambda \rightarrow 0^{-}$ and $\lambda \rightarrow -\infty$. In the latter case we get results when $2 < p < q < 2_{s_{2}}^{\ast}$ while get results when $2 < p < q < 2_{s_{1}}^{\ast}$ in the case of $\lambda \rightarrow 0^{-}$ ($2_{s_{1}}^{\ast} < 2_{s_{2}}^{\ast}$). If $\lambda$ is viewed as a Lagrange multiplier and unknown, this difference is difficult to distinguish.
		 \end{remark}
		
		 \begin{remark}
		 	Along the lines of Subsection 5.2, we can discuss the orbital stability or instability of ground states when $\lambda \rightarrow 0^{-}$ and $\lambda \rightarrow -\infty$. The situation is more complicated than the case of a single fractional Laplacian and details of the proof are left to readers.
		 \end{remark}
	
		 \begin{remark}
	        In addition to the applications in this article, we believe that our abstract framework applies to more equations. We are here to list some but not all of them.
	
	        $(i)$ Schr\"{o}dinger-Poisson-Slater equations:
	       	\begin{equation}
	       	\left\{
	       	\begin{array}{cc}
	       		-\Delta u + (|x|^{-1} \ast u^2) = \lambda u + f(|x|,u) \ in \ \mathbb{R}^{3},\\
	       	    u(x) \rightarrow 0 \ as \ |x| \rightarrow +\infty,
	       	\end{array}
	       	\right.
	       \end{equation}
	        where $f$ is $L^2$-supercritical and $H^1$-subcritical.
	
	        $(ii)$ Chern-Simons-Schr\"{o}dinger equations:
	        \begin{equation}
	        	\left\{
	        	\begin{array}{cc}
	        		-\Delta u + (\int_{|x|}^\infty\frac{h(s)}{s}u^2(s)ds)u + \frac{h^2(|x|)}{|x|^2}u = \lambda u + f(|x|,u) \ in \ \mathbb{R}^{2},\\
	        		u \in H_{rad}^1(\mathbb{R}^2),
	        	\end{array}
	        	\right.
	        \end{equation}
	        where $h(s) = \int_0^s\frac{l}{2}u^2(l)dl$, $f$ is $L^2$-supercritical and $H^1$-subcritical.
	
	        $(iii)$ NLS with a partial confinement:
	        \begin{equation}
	        	\left\{
	        	\begin{array}{cc}
	        		-\Delta u + |y|^k u = \lambda u + |u|^{p-2}u \ in \ \mathbb{R}^{N}, \\
	        		u(x) \rightarrow 0 \ as \ |x| \rightarrow +\infty,
	        	\end{array}
	        	\right.
	        \end{equation}
        where $x = (y,z) \in \mathbb{R}^{N-d}\times \mathbb{R}^{d}$, $k > 0$, $N \geq 2$, $1 \leq d \leq N-1$. Bellazzini, Boussaid, Jeanjean, and Visciglia, in \cite{BBJV}, proved the existence of orbitally stable ground states to NLS with a partial confinement together with qualitative and symmetry properties when the nonlinearity is $L^{2}$-subcritical with respect to the set of variables bearing no potential. We can cover the $L^{2}$-critical and $L^{2}$-supercritical cases, and to the best of our knowledge, this is the first method which can do this. Furthermore, we can give the information about uniqueness, non-degeneracy, multiplicity and bifurcation phenomena, allowing to  answer some conjectures in \cite{BBJV} (\cite[Remarks 1.8, 1.10]{BBJV}), see \cite{Song4}.
         \end{remark}

        \appendix

        \section{\textbf{Another method to prove Theorem \ref{thm4.9}}} \label{Appen1}

        The main aim of this appendix is to give a new method to prove Theorem \ref{thm4.9}. For the readers' convenience, we give a brief outline of our proof here. When $s = 1$, it is well known that positive solution of $-\Delta u = \lambda u + |u|^{q-2}u \ (\lambda < 0)$ is unique in $H^{1}(\mathbb{R}^{N})$. When $0 < s < 1$, \cite{FL, FLS} showed that positive solution with Morse index 1 of $(-\Delta)^{s}u = \lambda u + |u|^{q-2}u \ (\lambda < 0)$ is unique in $H^{s}(\mathbb{R}^{N})$. We want to find a suitable path connecting (\ref{eq4.14}) to the autonomous case. More precisely, consider
		\begin{eqnarray}
		G: H^{2s}_{rad}(\mathbb{R}^{N}) \times [0,1] \rightarrow L^{2}_{rad}(\mathbb{R}^{N}) \nonumber \\
		G(u,\zeta) = (-\Delta)^{s} u - \lambda u - h^{\zeta}(|x|)|u|^{q-2}u. \nonumber
		\end{eqnarray}
		Our first step is to show that all positive solutions of $G(u,\zeta) = 0$ with Morse index $1$ is non-degenerate in $H^{s}_{rad}(\mathbb{R}^{N})$ (the definition will be given in Definition \ref{def4.12}), see Lemma \ref{lem4.13}. Let $u^{\ast}$ be a positive ground state of $G(u,1) = 0$. Then applying the implicit function theorem when $\zeta = 1$ with initial value $u^{\ast}$, we construct a locally unique branch $u_{\zeta}$ nearby $u^{\ast}$ parameterized by $\zeta$ close to $1$ and satisfying $G(u_{\zeta},\zeta) = 0$, see Proposition \ref{prop4.14}. Next, we show that the local branch can be indeed globally continued to $\zeta = 0$. Let
		\begin{eqnarray}
		\zeta^{\ast} &=& \inf \{\tilde{\zeta} \in (0,1): u_{\zeta} \in C^{1}((\tilde{\zeta},1], H^{2s}_{rad}(\mathbb{R}^{N})), \nonumber \\
		&& u_{\zeta} \ satisfies \ the \ assumptions \ of \ Proposition \ \ref{prop4.14} \ for \ \zeta \in (\tilde{\zeta},1]\}. \nonumber
		\end{eqnarray}
		With the help of some uniform a priori estimates, see Lemma \ref{lem4.15}, we show that $\zeta^{\ast} = 0$ in Lemma \ref{lem4.18}. Finally, if $G(u,1) = 0$ admits two positive ground state solutions $u$ and $\tilde{u}$, then there are two branches $u_{\zeta}$ and $\tilde{u}_{\zeta}$ which do not coincide. Since $u_{0}$ and $\tilde{u}_{0}$ are positive and have Morse index $1$, $u_{0} = \tilde{u}_{0}$ by the uniqueness results of $G(u,0) = 0$. However, it is a contradiction since there are two different local branches with initial value $u_{0}$ parameterized by $\zeta$ close to $0$. Therefore, we know that $G(u,1) = 0$ admits at most one positive ground state solution.
		
		\begin{definition}
			\label{def4.11}(Morse index)  Let $u_{\zeta} \in H^{s}(\mathbb{R}^{N})$ be a solution of $G(u,\zeta) = 0$, and we define the linearized operator
			\begin{equation} \label{eq4.15}
			L_{+,\zeta} = (-\Delta)^{s} - \lambda - (q-1)h^{\zeta}(|x|)|u_{\zeta}|^{q-2}.
			\end{equation}
			$L_{+,\zeta}$ is a self-adjoint operator on $L^{2}(\mathbb{R}^{N})$ with quadratic-form domain $H^{s}(\mathbb{R}^{N})$ and operator domain $H^{2s}(\mathbb{R}^{N})$. The Morse index of $u_{\zeta}$ is defined as
			\[
			\mu(u_{\zeta}):= \sharp \{e < 0: e \ is \ an \ eigenvalue \ of \ L_{+,\zeta}\},
			\]
			where multiplicities of eigenvalues are taken into account. Noticing that $L_{rad}^{2}(\mathbb{R}^{N})$ is an invariant subspace of $L_{+,\zeta}$, we also introduce the Morse index of $u_{\zeta}$ in the sector of radial functions by defining
			\[
			\mu_{rad}(u_{\zeta}):= \sharp \{e < 0: e \ is \ an \ eigenvalue \ of \ L_{+,\zeta} \ restricted \ to \ L_{rad}^{2}(\mathbb{R}^{N})\}.
			\]
		\end{definition}
		
		If $u_{\zeta}$ is a ground state solution of $G(u,\zeta) = 0$, then $\mu(u_{\zeta}) = \mu_{rad}(u_{\zeta}) = 1$. On the one hand,
		$$
		\langle L_{+,\zeta}u_{\zeta},u_{\zeta}\rangle_{L^{2}} = (2-q)\int_{\mathbb{R}^{N}}h^{\zeta}(|x|)|u_{\zeta}|^{q}dx < 0.
		$$
		Thus, by the min-max principle, the operator $L_{+,\zeta}$ has at least one negative eigenvalue, i.e. $\mu(u_{\zeta}) \geq 1$. Since $u_{\zeta}$ is radial, then $\mu_{rad}(u_{\zeta}) \geq 1$ applying min-max principle in $L_{rad}^{2}(\mathbb{R}^{N})$. On the other hand, since $u_{\zeta}$ is a minimum in Nehari manifold whose codimension is $1$, $\mu_{rad}(u_{\zeta}) \leq \mu(u_{\zeta}) \leq 1$.
		
		\begin{definition}
			\label{def4.12}(non-degeneracy)  Let $u_{\zeta} \in H_{rad}^{s}(\mathbb{R}^{N})$ be a solution of $G(u,\zeta) = 0$. We say that $u_{\zeta}$ is non-degenerate in $H_{rad}^{s}(\mathbb{R}^{N})$ iff $\ker L_{+,\zeta|L_{rad}^{2}} = \{0\}$, where $L_{+,\zeta|L_{rad}^{2}}$ means $L_{+,\zeta}$ restricted to $L_{rad}^{2}(\mathbb{R}^{N})$. In other words, if $v \in H_{rad}^{s}(\mathbb{R}^{N})$ solves $(-\Delta)^{s}v = \lambda v + h^{\zeta}(|x|)|u_{\zeta}|^{q-2}v$, then $v = 0$.
		\end{definition}
		
		\begin{lemma}
			\label{lem4.13}  Let $u_{\zeta} \in H_{rad}^{s}(\mathbb{R}^{N})$ be a positive solution of $G(u,\zeta) = 0$ and $\mu_{rad}(u_{\zeta}) = 1$. Suppose that $(h_{1})$ holds true, then $u_{\zeta}$ is non-degenerate in $H_{rad}^{s}(\mathbb{R}^{N})$. In particular, all the ground states of $G(u,\zeta) = 0$ are non-degenerate in $H_{rad}^{s}(\mathbb{R}^{N})$.
		\end{lemma}
		
		\textit{Proof.  } Noticing that $h^{\zeta}$ satisfies $(h_{1})$ for all $\zeta \in [0,1)$ if $h$ satisfies $(h_{1})$, we merely give the proof when $\zeta = 1$. This is what has been proved in Lemma \ref{lem4.13b}. \qed\vskip 5pt
		
		As a next step, we will construct a local branch of solutions $u_{\zeta} \in H_{rad}^{2s}$, which is parameterized by $\zeta$ in some small interval. In fact, by a bootstrap argument, we see that $u_{\zeta} \in H_{rad}^{2s}$ holds if $u_{\zeta} \in H_{rad}^{s}$ satisfying $G(u_{\zeta},\zeta) = 0$.
		
		\begin{proposition}
			\label{prop4.14} Let $u^{\ast} \in H_{rad}^{2s}(\mathbb{R}^{N})$ be a positive solution of $G(u,1) = 0$ and $\mu_{rad}(u^{\ast}) = 1$. Suppose that $(h_{1})$ holds true, then for some $\delta > 0$, there exists a map $u(\zeta) \in C^{1}((1-\delta,1], H_{rad}^{2s})$ such that the following holds true, we denote $u_{\zeta} = u(\zeta)$ in the sequel.
			
			$(i)$ $u_{\zeta}$ solves $G(u,\zeta) = 0$ for all $\zeta \in (1-\delta,1]$.
			
			$(ii)$ There exists $\epsilon > 0$ such that $u_{\zeta}$ is the unique solution of $G(u,\zeta) = 0$ for $\zeta \in (1-\delta,1]$ in the neighborhood $\{u \in H_{rad}^{2s}: \|u - u^\ast\|_{H_{rad}^{2s}} < \epsilon\}$. In particular, we have that $u_{1} = u^{\ast}$ holds.
		\end{proposition}
		
		\textit{Proof.  } We use an implicit function argument for the map
		\begin{equation} \label{eq4.19}
		\Psi: H_{rad}^{2s} \times (1-\delta,1] \rightarrow H_{rad}^{2s}, \Psi(u,\zeta) = u - \frac{1}{(-\Delta)^{s}-\lambda}h^{\zeta}(|x|)u^{q-1}.
		\end{equation}
		It is not difficult to see that $\Psi$ is a well-defined map of class $C^{1}$. Furthermore, $\Psi(u^{\ast},1) = 0$ and $\Psi(u,\zeta) = 0$ if and only if $u \in H_{rad}^{2s}$ solves $G(u,\zeta) = 0$. Next, we consider the Fr\'{e}chet derivative
		\begin{equation} \label{eq4.20}
		\Psi_{u}(u,\zeta) = 1 - K, K = \frac{1}{(-\Delta)^{s}-\lambda}(q-1)h^{\zeta}(|x|)u^{q-2}.
		\end{equation}
		By Lemma \ref{lem4.13}, $0$ is not an eigenvalue of $L_{+|L^{2}_{rad}} = (-\Delta)^{s} - \lambda - (q-1)h^{\zeta}(|x|)u^{q-2}$. Note that the operator $K$ is compact on $L_{rad}^{2}$ and we have that $1 \notin \sigma(K)$. Since $h^{\zeta}(|x|)u^{q-2}\varphi \in L^{2}_{rad}$ if $\varphi \in H_{rad}^{2s}$, we have that $K$ maps $H_{rad}^{2s}$ to $H_{rad}^{2s}$. Thus the bounded inverse $(1 - K)^{-1}$ exists on $H_{rad}^{2s}$, i.e. the Fr\'{e}chet derivative $\Psi_{u}$ has a bounded inverse on $H_{rad}^{2s}$ at $(u^{\ast},1)$. By the implicit function theorem, we deduce that the assertions $(i)$ and $(ii)$ in Proposition \ref{prop4.14} hold for some $\delta > 0$ and $\epsilon > 0$ sufficiently small. \qed\vskip 5pt
		
		We now turn to the global extension of the locally unique branch $u_{\zeta}$ constructed in the above Proposition. First we will give some uniform a priori estimates. We write $a \lesssim b$ to denote that $a \leq Cb$ with some constant $C > 0$ that may change from line to line and only depends on $N$, $q$, $s$, $\lambda$, $h(r)$ and $u^{\ast}$. Furthermore, we use $a \thicksim b$ to denote that both $a \lesssim b$ and $b \lesssim a$ hold.
		
		\begin{lemma}
			\label{lem4.15} Suppose that $(h_{1})$ holds with $(N-2s)q \leq 2(N+\theta)$. Then we have the following power estimates
			\[
			\int_{\mathbb{R}^{N}}h^{\zeta}(|x|)|u_{\zeta}|^{q}dx \thicksim \int_{\mathbb{R}^{N}}u_{\zeta}^{2}dx \thicksim \int_{\mathbb{R}^{N}}|(-\Delta)^{\frac{s}{2}}u_{\zeta}|^{2}dx \thicksim 1,
			\]
			for all $\zeta \in (\zeta^{\ast},1]$, where $u_{\zeta}$ is given in Proposition \ref{prop4.14} and
			\begin{eqnarray}
			\zeta^{\ast} &=& \inf \{\tilde{\zeta} \in (0,1): u_{\zeta} \in C^{1}((\tilde{\zeta},1], H^{2s}_{rad}(\mathbb{R}^{N})), \nonumber \\
			&& u_{\zeta} \ satisfies \ the \ assumptions \ of \ Proposition \ \ref{prop4.14} \ for \ \zeta \in (\tilde{\zeta},1]\}. \nonumber
			\end{eqnarray}
		\end{lemma}
		
		\textit{Proof.  }  By integrating $G(u_{\zeta},\zeta) = 0$ with respect to $u_{\zeta}$, we obtain
		\begin{equation} \label{eq4.21}
		\int_{\mathbb{R}^{N}}|(-\Delta)^{s}u_{\zeta}|^{2}dx = \lambda\int_{\mathbb{R}^{N}}u_{\zeta}^{2}dx + \int_{\mathbb{R}^{N}}h^{\zeta}(|x|)|u_{\zeta}|^{q}dx.
		\end{equation}
		Similar to (\ref{eq4.6}) and (\ref{eq4.7}), Pohozaev identity and (\ref{eq4.21}) yield to
		\begin{equation} \label{eq4.22}
		-2s\lambda\int_{\mathbb{R}^{N}}u_{\zeta}^{2}dx = \int_{\mathbb{R}^{N}} [\frac{2N}{q} + \frac{2\zeta}{q}\frac{|x|h'(|x|)}{h(|x|)} - (N-2s)]h^{\zeta}(|x|)|u_{\zeta}|^{q}dx,
		\end{equation}
		\begin{equation} \label{eq4.23}
		2s\int_{\mathbb{R}^{N}}|(-\Delta)^{s}u_{\zeta}|^{2}dx = \int_{\mathbb{R}^{N}}(N - \frac{2N}{q} - \frac{2\zeta}{q}\frac{|x|h'(|x|)}{h(|x|)})h^{\zeta}(|x|)|u_{\zeta}|^{q}dx.
		\end{equation}
		By $(h_{1})$, $\theta \leq \frac{|x|h'(|x|)}{h(|x|)} \leq 0$. Therefore,
		\begin{eqnarray} \label{eq4.24}
		&& \int_{\mathbb{R}^{N}} [\frac{2(N+\theta)}{q} - (N-2s)]h^{\zeta}(|x|)|u_{\zeta}|^{q}dx \nonumber \\
		&\leq& -2s\lambda\int_{\mathbb{R}^{N}}u_{\zeta}^{2}dx \nonumber \\
		&\leq& \int_{\mathbb{R}^{N}} [\frac{2N}{q} - (N-2s)]h^{\zeta}(|x|)|u_{\zeta}|^{q}dx,
		\end{eqnarray}
		\begin{eqnarray} \label{eq4.25}
		\int_{\mathbb{R}^{N}}\frac{q-2}{q}Nh^{\zeta}(|x|)|u_{\zeta}|^{q}dx &\leq& 2s\int_{\mathbb{R}^{N}}|(-\Delta)^{s}u_{\zeta}|^{2}dx \nonumber \\
		&\leq& \int_{\mathbb{R}^{N}}(N - \frac{2(N+\theta)}{q})h^{\zeta}(|x|)|u_{\zeta}|^{q}dx.
		\end{eqnarray}
		Since $2 < q < \frac{2(N+\theta)}{N-2s}$, we have
		$$\int_{\mathbb{R}^{N}}h^{\zeta}(|x|)|u_{\zeta}|^{q}dx \thicksim \int_{\mathbb{R}^{N}}u_{\zeta}^{2}dx \thicksim \int_{\mathbb{R}^{N}}|(-\Delta)^{s}u_{\zeta}|^{2}dx.$$
		
		From the following fractional Gagliardo-Nirenberg-Sobolev inequality
		\begin{equation} \label{eq4.26}
		\int_{\mathbb{R}^{N}}|u|^{q}dx \leq C_{opt}(\int_{\mathbb{R}^{N}}|(-\Delta)^{\frac{s}{2}}u|^{2}dx)^{\frac{N(q-2)}{4s}} (\int_{\mathbb{R}^{N}}u^{2}dx)^{\frac{q}{2} - \frac{N(q-2)}{4s}}, \forall u \in H^{s},
		\end{equation}
		where $C_{opt} > 0$ denotes the sharp constant and depends on $s$, $N$, $q$, we can derive
		\begin{eqnarray} \label{eq4.27}
		&& \int_{\mathbb{R}^{N}}h^{\zeta}(|x|)|u_{\zeta}|^{q}dx \nonumber \\
		&\leq& C_{opt}\max\{\|h\|_{L^{\infty}},1\}(\int_{\mathbb{R}^{N}}|(-\Delta)^{s}u_{\zeta}|^{2}dx)^{\frac{N(q-2)}{4s}} (\int_{\mathbb{R}^{N}}u_{\zeta}^{2}dx)^{\frac{q}{2} - \frac{N(q-2)}{4s}} \nonumber \\
		&\lesssim& (\int_{\mathbb{R}^{N}}h^{\zeta}(|x|)|u_{\zeta}|^{q}dx)^{\frac{q}{2}}. \nonumber
		\end{eqnarray}
		Thus $\int_{\mathbb{R}^{N}}h^{\zeta}(|x|)|u_{\zeta}|^{q}dx \gtrsim 1$ for $\zeta \in (\zeta^{\ast},1]$.
		
		Next, we show that $\int_{\mathbb{R}^{N}}h^{\zeta}(|x|)|u_{\zeta}|^{q}dx \lesssim 1$ for $\zeta \in (\zeta^{\ast},1]$. Differentiating $G(u_{\zeta},\zeta) = 0$ with respect to $\zeta$, we get
		\begin{eqnarray} \label{eq4.28}
		&& (-\Delta)^{s}\frac{d}{d\zeta}u_{\zeta} = \nonumber \\
		&& \lambda \frac{d}{d\zeta}u_{\zeta} + (q-1)h^{\zeta}(|x|)|u_{\zeta}|^{q-2}\frac{d}{d\zeta}u_{\zeta} + \ln h(|x|)h^{\zeta}(|x|)|u_{\zeta}|^{q-2}u_{\zeta},
		\end{eqnarray}
		showing that $$\frac{d}{d\zeta}u_{\zeta} = L_{+,\zeta}^{-1}(\ln h(|x|)h^{\zeta}(|x|)|u_{\zeta}|^{q-2}u_{\zeta}).$$ Since $L_{+,\zeta}u_{\zeta} = (2-q)h^{\zeta}(|x|)|u_{\zeta}|^{q-2}u_{\zeta}$,
		\begin{eqnarray} \label{eq4.29}
		&& \langle h^{\zeta}(|x|)|u_{\zeta}|^{q-2}u_{\zeta}, \frac{d}{d\zeta}u_{\zeta}\rangle_{L^{2}} \nonumber \\
		&=& \langle L_{+,\zeta}^{-1}(h^{\zeta}(|x|)|u_{\zeta}|^{q-2}u_{\zeta}), \ln h(|x|)h^{\zeta}(|x|)|u_{\zeta}|^{q-2}u_{\zeta}\rangle_{L^{2}} \nonumber \\
		&=& \frac{1}{2-q}\int_{\mathbb{R}^{N}}\ln h(|x|)h^{\zeta}(|x|)|u_{\zeta}|^{q}dx.
		\end{eqnarray}
		We write $\phi(\zeta) = \int_{\mathbb{R}^{N}}h^{\zeta}(|x|)|u_{\zeta}|^{q}dx$, $\zeta \in (\zeta^{\ast},1]$. Then
		\begin{eqnarray} \label{eq4.30}
		\phi'(\zeta) &=& \int_{\mathbb{R}^{N}}\ln h(|x|)h^{\zeta}(|x|)|u_{\zeta}|^{q}dx + q\langle h^{\zeta}(|x|)|u_{\zeta}|^{q-2}u_{\zeta}, \frac{d}{d\zeta}u_{\zeta}\rangle_{L^{2}} \nonumber \\
		&=& -\frac{2}{q-2}\int_{\mathbb{R}^{N}}\ln h(|x|)h^{\zeta}(|x|)|u_{\zeta}|^{q}dx \nonumber \\
		&\geq& -\frac{2\ln h(0)}{q-2}\int_{\mathbb{R}^{N}}h^{\zeta}(|x|)|u_{\zeta}|^{q}dx \nonumber \\
		&=& -\frac{2\ln h(0)}{q-2}\phi(\zeta).
		\end{eqnarray}
		When $h(0) \leq 1$, $\phi'(\zeta) \geq 0$. Thus $\phi(\zeta) \leq \phi(1)$ for any $\zeta \in (\zeta^{\ast},1]$. When $h(0) > 1$, let $\chi(\zeta) = e^{-\frac{2\ln h(0)}{q-2}\zeta}$. Then
		$$
		\chi'(\zeta) = -\frac{2\ln h(0)}{q-2}e^{-\frac{2\ln h(0)}{q-2}\zeta} = -\frac{2\ln h(0)}{q-2}\chi(\zeta)
		$$
		and
		$$
		(\frac{\phi}{\chi})' = \frac{\chi\phi' - \chi'\phi}{\chi^{2}} \geq 0.
		$$
		Therefore, $\frac{\phi(\zeta)}{\chi(\zeta)} \leq \frac{\phi(1)}{\chi(1)}$, $\forall \zeta \in (\zeta^{\ast},1]$, showing
		$$
		\phi(\zeta) \leq \frac{\phi(1)}{\chi(1)}\chi(\zeta) \leq e^{\frac{2\ln h(0)}{q-2}}\phi(1), \forall \zeta \in (\zeta^{\ast},1].
		$$
		Thus we have shown that $\int_{\mathbb{R}^{N}}h^{\zeta}(|x|)|u|^{q}dx \lesssim 1$ for $\zeta \in (\zeta^{\ast},1]$.
		\qed\vskip 5pt
		
		\begin{lemma}
			\label{lem4.16} Suppose that $(h_{1})$ holds with $(N-2s)q \leq 2(N+\theta)$. $u_{\zeta}$ is given in Proposition \ref{prop4.14} and $u^{\ast} > 0$. Then for all $\zeta \in (\zeta^{\ast},1]$, we have
			
			$(i)$ $u_{\zeta}(x) > 0$ for all $x \in \mathbb{R}^{N}$;
			
			$(ii)$ (uniform decay estimate) $u_{\zeta}(x) \leq C|x|^{-(N+2s)}$ for $|x| > R$, where $R$ is some constant independent of $\zeta$;
			
			$(iii)$ $\mu_{rad}(u_{\zeta}) = 1$.
		\end{lemma}
		
		\textit{Proof.  } Define the linearized operator
		\begin{equation} \label{eq4.31}
		L_{-,\zeta} = (-\Delta)^{s} - \lambda - h^{\zeta}(|x|)|u_{\zeta}|^{q-2}.
		\end{equation}
		Obviously, if $\zeta \rightarrow \tilde{\zeta} \in (\zeta^{\ast},1]$, then $L_{\pm,\zeta} \rightarrow L_{\pm,\tilde{\zeta}}$ in norm-resolvent sense.
		
		The proof of $(i)$ is divided into two steps. First, we show that if $u_{\tilde{\zeta}} > 0$ for some $\tilde{\zeta} \in (\zeta^{\ast},1]$, then $u_{\zeta} > 0$ for $\zeta$ close to $\tilde{\zeta}$. Note that $L_{-,\zeta}u_{\zeta} = 0$, i.e. $0$ is an eigenvalue of $L_{-,\zeta}$ and $u_{\zeta}$ is an eigenfunction with respect to $0$. Since $\sigma_{ess}(L_{-,\zeta}) = [-\lambda,+\infty)$, we know that $\inf \sigma(L_{-,\zeta})$ is an eigenvalue. \cite[Theorem XIII.48]{RS} yields that the first eigenvalue is simple and the corresponding eigenfunction can be chosen strictly positive when $s = 1$ and similar result when $0 < s < 1$ can be found in \cite[Lemma C.4]{FLS}. Thus $\lambda_{1}(L_{-,\tilde{\zeta}}) = 0$ is simple and we choose $\psi_{1,\tilde{\zeta}} = u_{\tilde{\zeta}} > 0$ as the corresponding eigenfunction, where $\lambda_{1}(L_{-,\tilde{\zeta}})$ means the first eigenvalue of $L_{-,\tilde{\zeta}}$. Since $L_{-,\zeta} \rightarrow L_{-,\tilde{\zeta}}$ in norm-resolvent sense when $\zeta \rightarrow \tilde{\zeta} \in (\zeta^{\ast},1]$, we derive that $\lambda_{1}(L_{-,\zeta}) \rightarrow 0$. Since $\lambda_{1}(L_{-,\zeta})$ is simple and isolate, we may assume that $\sigma(L_{-,\zeta}) \cap (-\epsilon,\epsilon) = \{\lambda_{1}(L_{-,\zeta})\}$ for small $\epsilon$ and $\zeta$ close enough to $\tilde{\zeta}$. Furthermore, $0$ is an eigenvalue of $L_{-,\zeta}$, showing $\lambda_{1}(L_{-,\zeta}) = 0$. Since $u_{\zeta} \rightarrow u_{\tilde{\zeta}} > 0$ in $H^{2s}$, we obtain $u_{\zeta} > 0$ for close enough to $\tilde{\zeta}$.
		
		The second step is to show that if $\{\zeta_{n}\}_{n=1}^{+\infty} \subset (\hat{\zeta},1]$ is a sequence with $\zeta_{n} \rightarrow \hat{\zeta}$ such that $u_{\zeta_{n}} > 0$ for some $\hat{\zeta} \in (\zeta^{\ast},1]$, then $u_{\hat{\zeta}} > 0$. Obviously, $u_{\hat{\zeta}} \geq 0$ since $u_{\zeta_{n}} \rightarrow u_{\hat{\zeta}}$ in $H^{2s}$. Thus $0$ is the first eigenvalue of $L_{-,\hat{\zeta}}$ and $u_{\hat{\zeta}} > 0$. The proof of $(i)$ is complete.
		
		Next, we prove $(ii)$. Since $u_{\zeta}$ is monotonically decreasing with respect to $r = |x|$, we have $|x|^{N}u_{\zeta}^{2}(x) \lesssim u_{\zeta}^{2}(x)\int_{B_{|x|}}dy \leq \int_{B_{|x|}}u_{\zeta}^{2}(y)dy \lesssim 1$ (we use the uniform a priori estimate shown in Lemma \ref{lem4.15} in the final inequality). Thus $u_{\zeta}(x) \lesssim |x|^{-\frac{N}{2}}$.
		
		We write $u_{\zeta} = \frac{1}{(-\Delta)^{s} + \omega}z_{\zeta}$ where $0 < \omega < -\lambda$ and $z_{\zeta} = (\omega + \lambda)u_{\zeta} + h^{\zeta}u_{\zeta}^{q-1}$. Since $u_{\zeta}(x) \lesssim |x|^{-\frac{N}{2}}$, we can take $R$ large enough s.t. $supp z_{\zeta}^{+} \subset B_{\frac{R}{2}}$ and $R$ is independent of $\zeta$. Suppose $G_{s,\omega}(x)$ is the kernel of the resolvent $((-\Delta)^{s} + \omega)^{-1}$ on $\mathbb{R}^{N}$ with $\omega > 0$. In other words, $G_{s,\omega}(x)$ denotes the Fourier transform of $(|\xi|^{2s} + \omega)^{-1}$. \cite[Lemma C.1]{FLS} shows that $0 < G_{s,\omega}(x) \leq C|x|^{-(N+2s)}$ for $|x| > 0$ where $C > 0$ depends on $s$, $\omega$ and $N$ when $0 < s < 1$ and this result also hold when $s = 1$. Notice that
		\begin{eqnarray} \label{eq4.32}
		\int_{B_{\frac{R}{2}}}z_{\zeta}^{+}(y)dy &\leq& \int_{B_{\frac{R}{2}} \cap supp z_{\zeta}^{+}}(h^{\zeta})^{\frac{1}{q}}[(h^{\zeta})^{\frac{1}{q}}u_{\zeta}]^{q-1}dy \nonumber \\
		&\leq& (\int_{B_{\frac{R}{2}} \cap supp z_{\zeta}^{+}}h^{\zeta}dy)^{\frac{1}{q}}(\int_{B_{\frac{R}{2}} \cap supp z_{\zeta}^{+}}h^{\zeta}u_{\zeta}^{q}dy)^{\frac{q-1}{q}} \nonumber \\
		&\lesssim& R^{\frac{N}{q}},
		\end{eqnarray}
		we can conclude that: For all $|x| > R$,
		\begin{eqnarray} \label{eq4.33}
		u_{\zeta}(x) &=& \int_{\mathbb{R}^{N}}G_{s,\omega}(x-y)z_{\zeta}(y)dy \nonumber \\
		&\leq& \int_{B_{\frac{R}{2}}}G_{s,\omega}(x-y)z_{\zeta}^{+}(y)dy \nonumber \\
		&\lesssim& \int_{B_{\frac{R}{2}}}|x-y|^{-(N+2s)}z_{\zeta}^{+}(y)dy \nonumber \\
		&\lesssim& \int_{B_{\frac{R}{2}}}|x|^{-(N+2s)}z_{\zeta}^{+}(y)dy \nonumber \\
		&\lesssim& |x|^{-(N+2s)}.
		\end{eqnarray}
		
		Finally, we prove $(iii)$. Arguing by contradiction, we may assume that there exists some $\tilde{\zeta} \in (\zeta^{\ast},1)$ such that $0$ is an eigenvalue of $L_{+,\tilde{\zeta}}$ restricted on $L^{2}_{rad}$ and $\mu_{rad}(u_{\tilde{\zeta}}) = 1$ by the continuity of the eigenvalues of $L_{+,\zeta}$. Since $(i)$ shows that $u_{\tilde{\zeta}}$ is strictly positive, $u_{\tilde{\zeta}}$ is non-degenerate in $H_{rad}^{s}$ by Lemma \ref{lem4.13}, which is a contradiction. The proof now is complete. \qed\vskip 5pt
		
		\begin{lemma}
			\label{lem4.17}  Let $\{\zeta_{n}\}_{n=1}^{+\infty} \subset (\zeta^{\ast},1]$ be a sequence such that $\zeta_{n} \rightarrow \zeta^{\ast}$, $u_{\zeta_{n}} = u_{\zeta_{n}}(|x|) > 0$ and $\mu_{rad}(u_{\zeta_{n}}) = 1$. Then, after possibly passing to a subsequence, we have $u_{\zeta_{n}} \rightarrow u_{\zeta^{\ast}}$ in $H^{2s}_{rad}$ where $\mu_{rad}(u_{\zeta^{\ast}}) = 1$ and $u_{\zeta^{\ast}} = u_{\zeta^{\ast}}(|x|) > 0$ satisfies
			\begin{equation} \label{eq4.34}
			(-\Delta)^{s}u_{\zeta^{\ast}} = \lambda u_{\zeta^{\ast}} + h^{\zeta^{\ast}}(|x|)u_{\zeta^{\ast}}^{q-1} \ in \ \mathbb{R}^{N}.
			\end{equation}
		\end{lemma}
		
		\textit{Proof.  } By Lemma \ref{lem4.15}, $u_{\zeta_{n}}$ is bounded in $H^{s}$. Passing to a subsequence if necessary, we assume that $u_{\zeta_{n}} \rightharpoonup u_{\zeta^{\ast}}$ in $H^{s}$, $u_{\zeta_{n}} \rightarrow u_{\zeta^{\ast}}$ in $L^{p}_{loc}$ for $1 \leq p < 2_{s}^{\ast}$, $u_{\zeta_{n}} \rightarrow u_{\zeta^{\ast}}$ a.e. on $\mathbb{R}^{N}$. With the help of the uniform decay estimate shown in Lemma \ref{lem4.16} $(ii)$, we obtain that $u_{\zeta_{n}} \rightarrow u_{\zeta^{\ast}}$ in $L^{p}$ for $1 \leq p < 2_{s}^{\ast}$. Bootstrap arguments implies that $u_{\zeta_{n}} \rightarrow u_{\zeta^{\ast}}$ in $H^{2s}$.
		
		Obviously, $u_{\zeta^{\ast}} \geq 0$ satisfies equation (\ref{eq4.34}). Thus $0$ is the first eigenvalue of $L_{-,\zeta^{\ast}}$, showing $u_{\zeta^{\ast}} > 0$.
		
		Finally, we show that $\mu_{rad}(u_{\zeta^{\ast}}) = 1$. On the one hand, $\mu_{rad}(u_{\zeta^{\ast}}) \geq 1$ since
		\[
		\langle L_{+,\zeta^{\ast}}u_{\zeta^{\ast}},u_{\zeta^{\ast}}\rangle_{L^{2}} = -(q-2)\int_{\mathbb{R}^{N}}u_{\zeta^{\ast}}^{q}dx < 0.
		\]
		One the other hand, since the Morse index is lower semi-continuous with respect to the norm-resolvent topology, we conclude that
		\[
		\mu_{rad}(u_{\zeta^{\ast}}) \leq \inf_{n \rightarrow +\infty}\mu_{rad}(u_{\zeta_{n}}) = 1,
		\]
		which completes the proof. \qed\vskip 5pt
		
		\begin{lemma}
			\label{lem4.18}  Let $u^{\ast} \in H_{rad}^{2s}(\mathbb{R}^{N})$ be a positive solution of $G(u,1) = 0$ and $\mu_{rad}(u^{\ast}) = 1$. Suppose $(h_{1})$ holds with $(N-2s)q \leq 2(N+\theta)$, then $\zeta^{\ast} = 0$ where $\zeta^{\ast}$ is given in Lemma \ref{lem4.15}.
		\end{lemma}
		
		\textit{Proof.  } Suppose on the contrary that $\zeta^{\ast} > 0$. By Lemma \ref{lem4.16}, for all $\zeta \in (\zeta^{\ast},1]$, $u_{\zeta}(x) > 0$ and $\mu_{rad}(u_{\zeta}) = 1$. Then by Lemma \ref{lem4.17}, there exists $u_{\zeta^{\ast}}$ in $H^{2s}_{rad}$ satisfying (\ref{eq4.34}) where $\mu_{rad}(u_{\zeta^{\ast}}) = 1$ and $u_{\zeta^{\ast}} = u_{\zeta^{\ast}}(|x|) > 0$. By Lemma \ref{lem4.13}, $u_{\zeta^{\ast}}$ is non-degenerate in $H_{rad}^{s}(\mathbb{R}^{N})$. Similar to the proof of Proposition \ref{prop4.14}, applying the implicit function argument when $\zeta = \zeta^{\ast}$, we find that $\zeta < \zeta^{\ast}$ such that $u_{\zeta}$ satisfies the assumptions of Proposition \ref{prop4.14}, contradicting the definition of $\zeta^{\ast}$. \qed\vskip 5pt
		
		Finally, we give the proof of Theorem \ref{thm4.9}.
		
		\textbf{Proof of Theorem \ref{thm4.9}:  } Note that all ground states of (\ref{eq4.14}) are in $H_{rad}^{s}$. We show the uniqueness in $H_{rad}^{s}$ here. Arguing by contradiction, (\ref{eq4.14}) admits two positive ground states $u_{1}$ and $\tilde{u}_{1}$ in $H_{rad}^{s}$. $\mu_{rad}(u_{1}) = \mu_{rad}(\tilde{u}_{1}) = 1$. By Proposition \ref{prop4.14} and Lemma \ref{lem4.18}, we get two global branches $u_{\zeta}$ and $\tilde{u}_{\zeta}$, $\zeta \in [0,1]$, starting from $u_{1}$ and $\tilde{u}_{1}$ respectively. By Lemma \ref{lem4.17}, both $u_{0}$ and $\tilde{u}_{0}$ are positive solutions of $G(u,0) = 0$ where $\mu_{rad}(u_{0}) = \mu_{rad}(\tilde{u}_{0}) = 1$, showing $u_{0} = \tilde{u}_{0}$. By Lemma \ref{lem4.13}, $u_{0}$ is non-degenerate in $H_{rad}^{s}(\mathbb{R}^{N})$. However, there are two different local branches $u_{\zeta}$ and $\tilde{u}_{\zeta}$, $\zeta \in [0,\delta)$ stemming form $u_{0}$, contradicting the non-degeneracy of $u_{0}$. \qed\vskip 5pt
		
		\section{\textbf{Proof of Theorem \ref{thm4.27}}} \label{Appen2}
		
		Along the lines of Appendix \ref{Appen1}, we can show the uniqueness of positive ground states for  $(\ref{eq4.36})$ when $\lambda \leq 0$ and give the proof of Theorem \ref{thm4.27}.
		
		\textit{Proof of Theorem \ref{thm4.27}.  } We apply implicit function argument in $H^{2s}_{rad}$. Consider
		\begin{eqnarray}
		K: H^{2s}_{rad} \times (0,1] \rightarrow L^{2}_{rad}(\mathbb{R}^{N}) \nonumber \\
		K(u,\eta) = (-\Delta)^{s} u + [\eta V + (1-\eta)\lambda_{1}]u - \lambda u + h(|x|)|u|^{q-2}u. \nonumber
		\end{eqnarray}
		Note that $\inf \sigma ((-\Delta)^{s} + \eta V + (1-\eta)\lambda_{1}) \geq \lambda_{1} > 0$. In fact,
		\begin{eqnarray}
		&& \inf \sigma ((-\Delta)^{s} + \eta V + (1-\eta)\lambda_{1}) \\
		&=& \inf_{u \in H^{s}}\frac{\int_{\mathbb{R}^{N}}[|(-\Delta)^{\frac{s}{2}}u|^{2} + \eta Vu^{2}]dx}{\int_{\mathbb{R}^{N}}u^{2}dx} + (1-\eta)\lambda_{1} \nonumber \\
		&\geq& \eta\inf_{u \in H^{s}}\frac{\int_{\mathbb{R}^{N}}[|(-\Delta)^{\frac{s}{2}}u|^{2} + Vu^{2}]dx}{\int_{\mathbb{R}^{N}}u^{2}dx} + (1-\eta)\lambda_{1} \nonumber \\
		&=& \eta \lambda_{1} + (1-\eta)\lambda_{1} = \lambda_{1}. \nonumber
		\end{eqnarray}
		Denote $\eta V + (1-\eta)\lambda_{1}$ by $V_{\eta}$. $(V_{1})$ holds, then $V_{\eta}$ and $V_{\eta} + 2V'_{\eta} - \lambda \geq 0$ are non-decreasing for any $\eta \in [0,1]$. By the proof of Lemma \ref{lem4.26}, $u_{\eta}$ is non-degenerate in $H_{rad}^{s}$ whenever $u_{\eta} \in H_{rad}^{s}$ solves $K(u,\eta) = 0$, $u_{\eta}$ is positive and $\mu_{rad}(u_{\eta}) = 1$. Let $u^{\ast} \in H_{rad}^{2s}$ be a positive solution of $K(u,1) = 0$ and $\mu_{rad}(u^{\ast}) = 1$. Applying implicit function theorem to $K(u,\eta)$ at $(u^{\ast},1)$, we construct a local branch of solutions $u_{\eta} \in H_{rad}^{2s}$, which is parameterized by $\eta$ in some small interval, i.e. for some $\delta > 0$, there exists a map $u(\eta) \in C^{1}((1-\delta,1], H_{rad}^{2s})$ such that the following holds, where we denote $u_{\eta} = u(\eta)$ in the sequel.
		
		$(i)$ $u_{\eta}$ solves $K(u,\eta) = 0$ for all $\eta \in (1-\delta,1]$.
		
		$(ii)$ There exists $\epsilon > 0$ such that $u_{\eta}$ is the unique solution of $K(u,\eta) = 0$ for $\eta \in (1-\delta,1]$ in the neighborhood $\{u \in H_{rad}^{2s}: \|u - u_{\eta}\|_{H_{rad}^{2s}} < \epsilon\}$. In particular, we have that $u_{1} = u^{\ast}$ holds.
		
		Let
		\begin{eqnarray}
		\eta^{\ast} &=& \inf \{\tilde{\eta} \in (0,1): u_{\eta} \in C^{1}((\tilde{\eta},1], H^{2s}_{rad}(\mathbb{R}^{N})), \nonumber \\
		&& u_{\eta} \ satisfies \ the \ above \ assumptions \ (i), (ii) \ for \ \eta \in (\tilde{\eta},1]\}. \nonumber
		\end{eqnarray}
		With the help of Lemma \ref{lem4.29} shown later, we can derive $\eta^{\ast} = 0$ along the lines of Lemmas \ref{lem4.15} - \ref{lem4.18}. Then like the proof of Theorem \ref{thm4.9}, we can complete the proof. \qed\vskip 5pt
		
		\begin{lemma}
			\label{lem4.29} Suppose $(h_{1})$ and $(V_{1})$ hold with $(N-2s)q \leq 2(N+\theta)$. Then there is a uniform a priori estimate:
			\[
			\int_{\mathbb{R}^{N}}(V_{\eta}+1)u_{\eta}^{2}dx \thicksim \int_{\mathbb{R}^{N}}h(|x|)|u_{\eta}|^{q}dx \thicksim \int_{\mathbb{R}^{N}}|(-\Delta)^{s}u_{\eta}|^{2}dx \thicksim 1,
			\]
			for all $\eta \in (\eta^{\ast},1]$, where $u_{\eta}$ and $\eta^{\ast}$ are given by the proof of Theorem \ref{thm4.27}.
			
		\end{lemma}
		
		\textit{Proof.  }  By integrating $K(u_{\eta},\eta) = 0$ with respect to $u_{\eta}$, we obtain that:
		\begin{equation} \label{eq4.40}
		\int_{\mathbb{R}^{N}}(|(-\Delta)^{s}u_{\eta}|^{2} + V_{\eta}u_{\eta}^{2})dx = \lambda\int_{\mathbb{R}^{N}}u_{\eta}^{2}dx + \int_{\mathbb{R}^{N}}h(|x|)|u_{\eta}|^{q})dx.
		\end{equation}
		Pohozaev identity and (\ref{eq4.40}) yield to
		\begin{eqnarray} \label{eq4.41}
		&& 2s\int_{\mathbb{R}^{N}}V_{\eta}u_{\eta}^{2}dx + \int_{\mathbb{R}^{N}}|x|V_{\eta}'(|x|)u_{\eta}^{2}dx - 2s\lambda\int_{\mathbb{R}^{N}}u_{\eta}^{2}dx \nonumber \\
		&=& \int_{\mathbb{R}^{N}} [\frac{2N}{q} + \frac{2|x|h'(|x|)}{qh(|x|)} - (N-2s)]h(|x|)|u_{\eta}|^{q}dx,
		\end{eqnarray}
		\begin{eqnarray} \label{eq4.42}
		&& 2s\int_{\mathbb{R}^{N}}|(-\Delta)^{s}u_{\eta}|^{2}dx - \int_{\mathbb{R}^{N}}|x|V_{\eta}'(|x|)u_{\eta}^{2}dx \nonumber \\
		&=& \int_{\mathbb{R}^{N}}(N - \frac{2N}{q} - \frac{2|x|h'(|x|)}{qh(|x|)})h(|x|)|u_{\eta}|^{q}dx,
		\end{eqnarray}
		\begin{eqnarray} \label{eq4.43}
		&& (\gamma + 2s)\int_{\mathbb{R}^{N}}|(-\Delta)^{s}u_{\eta}|^{2}dx + \int_{\mathbb{R}^{N}}(\gamma V_{\eta} - |x|V_{\eta}'(|x|))u_{\eta}^{2}dx \nonumber \\
		&=& \gamma\lambda\int_{\mathbb{R}^{N}}u_{\eta}^{2}dx + \int_{\mathbb{R}^{N}}(N + \gamma - \frac{2N}{q} - \frac{2|x|h'(|x|)}{qh(|x|)})h(|x|)|u_{\eta}|^{q}dx.
		\end{eqnarray}
		(\ref{eq4.42}) implies that
		\begin{equation} \label{eq4.44}
		2s\int_{\mathbb{R}^{N}}|(-\Delta)^{s}u_{\eta}|^{2}dx \geq
		\frac{q-2}{q}N\int_{\mathbb{R}^{N}}h(|x|)|u_{\eta}|^{q}dx.
		\end{equation}
		By (\ref{eq4.43}), we have that
		\begin{equation} \label{eq4.45}
		(\gamma+2s)\int_{\mathbb{R}^{N}}|(-\Delta)^{s}u_{\eta}|^{2}dx \leq
		[N + \gamma - \frac{2(N+\theta)}{q}]\int_{\mathbb{R}^{N}}h(|x|)|u_{\eta}|^{q}dx.
		\end{equation}
		(\ref{eq4.41}) yields
		\begin{equation} \label{eq4.46}
		2s\int_{\mathbb{R}^{N}}V_{\eta}u_{\eta}^{2}dx - 2s\lambda\int_{\mathbb{R}^{N}}u_{\eta}^{2}dx \leq [\frac{2N}{q} - (N-2s)]\int_{\mathbb{R}^{N}}h(|x|)|u_{\eta}|^{q}dx,
		\end{equation}
		and
		\begin{equation} \label{eq4.47}
		(\gamma+2s)\int_{\mathbb{R}^{N}}V_{\eta}u_{\eta}^{2}dx - 2s\lambda\int_{\mathbb{R}^{N}}u_{\eta}^{2}dx \geq [\frac{2(N+\theta)}{q} - (N-2s)]\int_{\mathbb{R}^{N}} h(|x|)|u_{\eta}|^{q}dx.
		\end{equation}
		(\ref{eq4.44}) - (\ref{eq4.47}) show that
		\[
		\int_{\mathbb{R}^{N}}(V_{\eta}+1)u_{\eta}^{2}dx \thicksim \int_{\mathbb{R}^{N}}h(|x|)|u_{\eta}|^{q}dx \thicksim \int_{\mathbb{R}^{N}}|(-\Delta)^{s}u_{\eta}|^{2}dx.
		\]
		
		From the fractional Gagliardo-Nirenberg-Sobolev inequality, we can derive
		\begin{eqnarray} \label{eq4.48}
		\int_{\mathbb{R}^{N}}h(|x|)|u_{\eta}|^{q}dx &\lesssim& (\int_{\mathbb{R}^{N}}|(-\Delta)^{s}u_{\eta}|^{2}dx)^{\frac{N(q-2)}{4s}} (\int_{\mathbb{R}^{N}}u_{\eta}^{2}dx)^{\frac{q}{2} - \frac{N(q-2)}{4s}} \nonumber \\
		&\lesssim& (\int_{\mathbb{R}^{N}}|(-\Delta)^{s}u_{\eta}|^{2}dx)^{\frac{N(q-2)}{4s}} (\int_{\mathbb{R}^{N}}(V_{\eta}+1)u_{\eta}^{2}dx)^{\frac{q}{2} - \frac{N(q-2)}{4s}} \nonumber \\
		&\lesssim& (\int_{\mathbb{R}^{N}}h(|x|)|u_{\eta}|^{q}dx)^{\frac{q}{2}}.
		\end{eqnarray}
		Thus $\int_{\mathbb{R}^{N}}h(|x|)|u_{\eta}|^{q}dx \gtrsim 1$ for $\eta \in (\eta^{\ast},1]$.
		
		Next, we show that $\int_{\mathbb{R}^{N}}h(|x|)|u_{\eta}|^{q}dx \lesssim 1$ for $\eta \in (\eta^{\ast},1]$. Differentiating $K(u_{\eta},\eta) = 0$ with respect to $\eta$, we get
		\begin{equation} \label{eq4.49}
		(-\Delta)^{s}\frac{d}{d\eta}u_{\eta} + V_{\eta}\frac{d}{d\eta}u_{\eta} + (V - \lambda_{1})u_{\eta} = \lambda \frac{d}{d\eta}u_{\eta} + (q-1)h(|x|)|u_{\eta}|^{q-2}\frac{d}{d\eta}u_{\eta},
		\end{equation}
		showing that $\frac{d}{d\eta}u_{\eta} = L_{+,\eta}^{-1}(\lambda_{1} - V)u_{\eta}$ where $L_{+,\eta} = (-\Delta)^{s} + V_{\eta} - \lambda - (q-1)h(|x|)|u_{\eta}|^{q-2}$. Since $L_{+,\eta}u_{\eta} = (2-q)h(|x|)|u_{\eta}|^{q-2}u_{\eta}$,
		\begin{eqnarray} \label{eq4.50}
		&& \langle h(|x|)|u_{\eta}|^{q-2}u_{\eta}, \frac{d}{d\eta}u_{\eta}\rangle_{L^{2}} \nonumber \\
		&=& \langle L_{+,V_{\eta},\eta}^{-1}(h(|x|)|u_{\eta}|^{q-2}u_{\eta}),(\lambda_{1} - V)u_{\eta}\rangle_{L^{2}} \nonumber \\
		&=& \frac{1}{2-q}\int_{\mathbb{R}^{N}}(\lambda_{1} - V)u_{\eta}^{2}dx.
		\end{eqnarray}
		We write $\phi(\eta) = \int_{\mathbb{R}^{N}}h(|x|)|u_{\eta}|^{q}dx$, $\eta \in (\eta^{\ast},1]$. Then
		\begin{eqnarray} \label{eq4.51}
		\phi'(\eta) &=& q\langle h(|x|)|u_{\eta}|^{q-2}u_{\eta}, \frac{d}{d\eta}u_{\eta}\rangle_{L^{2}} \nonumber \\
		&=& \frac{2}{q-2}\int_{\mathbb{R}^{N}}(V - \lambda_{1})u_{\eta}^{2}dx \nonumber \\
		&\geq& -\frac{2\lambda_{1}}{q-2}\int_{\mathbb{R}^{N}}h(|x|)|u_{\eta}|^{q}dx \nonumber \\
		&=& -\frac{2\lambda_{1}}{q-2}\phi(\eta).
		\end{eqnarray}
		Similar to the proof of Lemma \ref{lem4.15}, we have $\phi(\eta) \leq e^{\frac{2\lambda_{1}}{q-2}}\phi(1)$, $\forall \eta \in (\eta^{\ast},1]$. Thus we have shown that $\int_{\mathbb{R}^{N}}h(|x|)|u|^{q}dx \lesssim 1$ for $\eta \in (\eta^{\ast},1]$. \qed
\end{sloppypar}		
	\end{document}